\documentclass[12pt]{article}
\usepackage{amsmath}
\usepackage{graphicx}
\usepackage{enumerate}
\usepackage{natbib}
\usepackage{url} % not crucial - just used below for the URL 

\usepackage{mathrsfs}
\usepackage{psfrag,epsf}
\usepackage{overpic}
\usepackage{float}
\usepackage{subfigure} 
\usepackage{natbib}
\usepackage{amsfonts}
\usepackage[hypertexnames=false,hidelinks]{hyperref}
\usepackage{amssymb}
\usepackage{stmaryrd}
\usepackage{amsthm}
\usepackage{scalerel}
\usepackage{pifont}
\usepackage{cases}
\usepackage{enumitem} 
\usepackage{booktabs}

\newcommand{\blind}{0}

\newtheorem{lemma}{Lemma}
\newtheorem{theorem}{Theorem}
\newtheorem{assumption}{Assumption}
\newtheorem{proposition}{Proposition}

\newcommand{\R}{\mathbb{R}}
\renewcommand{\hat}[1]{\widehat{#1}}

\renewcommand{\P}{\mathbf{P}}
\newcommand{\E}{\mathbf{E}}
\newcommand{\ve}{\varepsilon}
\newcommand{\cov}{\textup{cov}}
\newcommand{\tr}{\textup{tr}}
\newcommand{\ttop}{^{\top}}
\newcommand{\ts}{\textstyle}
\newcommand{\var}{\operatorname{var}}
\renewcommand{\u}{\mathbf{u}}
\newcommand{\x}{\mathbf{x}}

\newcommand{\sps}{\scaleto{\,\mathcal{S}}{4.5pt}}
\newcommand{\SIGMA}{\Sigma}
\newcommand{\TILDEKAPPA}{\tilde\kappa}
\newcommand{\z}{\mathbf{z}}
\newcommand{\e}{\epsilon}

\begin{document}

\def\spacingset#1{\renewcommand{\baselinestretch}%
{#1}\small\normalsize} \spacingset{1}

%%%%%%%%%%%%%%%%%%%%%%%%%%%%%%%%%%%%%%%%%%%%%%%%%%%%%%%%%%%%%%%%%%%%%%%%%%%%%%

%\iffalse
{
  \title{\bf Testing Elliptical Models in High Dimensions}
  \author{Siyao Wang
  %\thanks{
    %\textit{}}
    \hspace{.2cm}\\
    Department of Statistics, University of California, Davis\\
    and \\
    Miles E. Lopes \\
    Department of Statistics, University of California, Davis}
    
    \date{}
  \maketitle
} 
%\fi

%\iffalse
\if0\blind
{
 % \bigskip
 % \bigskip
 % \bigskip
 % \begin{center}
%    {\Large \bf Testing Elliptical Models in High Dimensions}
%\end{center}
%  \medskip
} \fi
%\fi

\bigskip
\begin{abstract}
 Due to the broad applications of elliptical models, there is a long line of research on goodness-of-fit tests for empirically validating them. However, the existing literature on this topic is generally confined to low-dimensional settings, and to the best of our knowledge, there are no established goodness-of-fit tests for elliptical models that are supported by theoretical guarantees in high dimensions. In this paper, we propose a new goodness-of-fit test for this problem, and our main result shows that the test is asymptotically valid when the dimension and sample size diverge proportionally. Remarkably, it also turns out that the asymptotic validity of the test requires \emph{no assumptions} on the population covariance matrix. With regard to numerical performance, we confirm that the empirical level of the test is close to the nominal level across a range of conditions, and that the test is able to reliably detect non-elliptical distributions. Moreover, 
 when the proposed test is specialized to the problem of testing normality in high dimensions, we show that it compares favorably with a state-of-the-art method, and hence, this way of using the proposed test is of independent interest.
\end{abstract}

\noindent%
{\it Keywords: elliptical models, goodness-of-fit, high-dimensional statistics, Isserlis' Theorem} 
%\vfill

\spacingset{1.9} % DON'T change the spacing!
\section{Introduction}
\addtolength{\textheight}{.5in}%

A centered random vector in $\R^p$ with covariance matrix $\Sigma$ is said to have an elliptical distribution if it can be represented in the form $\xi\Sigma^{1/2}\u$, where $\u$ is a uniformly distributed random vector on the $p$-dimensional unit sphere, and $\xi$ is a non-negative scalar random variable that is independent of $\u$ and normalized to satisfy $\E(\xi^2)=p$~\citep{cambanis1981theory,muirhead2009aspects}.
Elliptical models
play an essential role in multivariate analysis and high-dimensional statistics, as they possess many attractive theoretical properties of normal models, and are also flexible enough to capture important non-normal effects that arise in applications. For instance, elliptical models have been widely adopted throughout the domains of finance~\citep{Embrechts:2011,gupta2013elliptically} and signal processing~\citep{Tyler_survey,mahmood2015modeling,hung2022robust,li2022mises}. Furthermore, elliptical models provide the foundation for a plethora of methods that address fundamental statistical problems, including covariance estimation~\citep{tyler1987distribution,sun2016robust,Nadler:2020}, regression~\citep{li1991sliced,lin2019sparse}, and inference in graphical models~\citep{vogel2011elliptical,liu2012transelliptical,drton2017structure}.

Based on the extensive use of elliptical models, the problem of validating them with goodness-of-fit tests has become a long-standing topic of research~\citep[e.g.][]{beran1979testing,koltchinskii2000testing,manzotti2002statistic,schott2002testing,
huffer2007test,albisetti2020testing,ducharme2020goodness,babic2021optimal,Li:2023}. However, this literature is generally confined to low-dimensional settings---often due to the use of the inverse sample covariance matrix in the construction of goodness-of-fit tests. Likewise, several authors have recently drawn attention to the need for new approaches to goodness-of-fit testing that can handle elliptical and related models in high dimensions~\citep{Drton:2023,chen2023normality,Li:2023,Lin:2023}. Yet, to the best of our knowledge, there are no established goodness-of-fit tests for elliptical models that are theoretically supported in high dimensions. Indeed, even in the simpler case of normal models, theoretical guarantees for goodness-of-fit tests in high dimensions are at the cutting edge of research~\citep{chen2023normality}.

In the current paper, we seek to rectify this situation by developing a new goodness-of-fit test for elliptical models and establishing its asymptotic validity when the dimension and sample size diverge proportionally. To outline our formulation of the problem, we deal with a set of centered i.i.d.~observations $\x_1,\dots,\x_n\in\R^p$ drawn from an unknown distribution $\P$, and the goal is to use the observations to test
\begin{equation}\label{eqn:gof}
    \mathsf{H}_0: \P\in\mathscr{P} \ \ \ \ \ \text{versus}  \ \ \ \ \ \mathsf{H}_1: \P\not\in\mathscr{P},
\end{equation}
where $\mathscr{P}$ is a class of elliptical distributions. To account for the role of high dimensionality and allow $p$ to increase as a function of $n$, the problem~\eqref{eqn:gof} should be understood more precisely as representing a sequence of hypothesis testing problems implicitly indexed by $n$, in which $\P$ and $\mathscr{P}$ both depend on $n$. The exact definition of the sequence of classes $\mathscr{P}$ will be given later in Assumption~\ref{Data generating model}, but in brief, the associated elliptical distributions   will only be subject to some simple moment conditions and will be allowed to have \emph{unrestricted covariance matrices}.

The idea of the proposed test is based on a special property of the moment structure of elliptical distributions. Specifically, if $\x_1=(x_{11},\dots,x_{1p})$ is drawn from an elliptical distribution for which each covariate is non-degenerate and has finite kurtosis, then all of the covariates must have the same kurtosis, denoted by $\kappa$.
In particular, when $\x_1$ is centered, this means $\kappa=\E(x_{1j}^4)/(\E(x_{1j}^2))^2$ for all $j=1,\dots,p$, as recorded in Lemma~\ref{lem:equal_kurt} with a short accompanying proof. 
It also turns out that $\kappa$ can be independently calculated in a ``coordinate-free'' way (based on Lemma~\ref{lem:quadform}), which leads to the following system of $p$ simultaneous equations that hold for any centered elliptical distribution with non-degenerate covariates and $\kappa<\infty$,
\begin{equation}\label{eqn:kurt}
    \kappa=\frac{\E(x_{1j}^4)}{(\E(x_{1j}^2))^2}=\frac{3(\var(\|\x_1\|_2^2)+\tr(\Sigma)^2)}{2\tr(\Sigma^2)+\tr(\Sigma)^2} \text{ \ \ \ \  for all  \ \ \ \ $j=1,\dots,p$},
\end{equation}
where $\Sigma$ denotes the covariance matrix of $\x_1$. The equations above will lead us to construct two distinct kurtosis estimates, denoted $\tilde\kappa$ and $\check{\kappa}$, and our strategy is to ``play them off of each other'' in designing a goodness-of-fit test. That is, our proposed test statistic will measure the discrepancy between $\tilde\kappa$ and $\check{\kappa}$, and the null hypothesis will be rejected when this discrepancy is sufficiently large. So, in essence, the proposed test seeks to detect violations of  equations~\eqref{eqn:kurt}, and this approach is well suited to the high-dimensional setting because larger numbers of covariates provide more opportunities to detect violations. Furthermore, \emph{the test can even detect non-elliptical distributions when all of the covariates have the same kurtosis}. This is because one of the kurtosis estimates is based on a structural property that generally does not hold for non-elliptical distributions, and our simulation results for Model (4) in Section~\ref{sec:power} demonstrate the power of our test in such cases.
%Lastly, one more merit of our approach is that it does not rely on the inverse sample covariance matrix, and so unlike many existing goodness-of-fit tests for elliptical models, our approach be readily applied to data for which $p>n$.

From a theoretical standpoint, our main result (Theorem~\ref{thm:main}) is a central limit theorem for the proposed test statistic under the null hypothesis. One of the main virtues of this result is that it allows the population covariance matrix to vary in an unrestricted manner as $n$ and $p$ diverge. An interesting theoretical consequence of this flexibility is that the difference $\tilde\kappa-\check{\kappa}$ may fluctuate over many different scales, depending on the correlation structure of the data. This is in contrast to more conventional high-dimensional central limit theorems, in which an unnormalized statistic often has fluctuations of the same order (such as $1/\sqrt n$ or $1/p$) for all of the allowed choices of the population covariance matrix. Accordingly, constructing an accurate estimate for the standard deviation of $\tilde\kappa-\check{\kappa}$ is especially crucial in our method, and proving the consistency of this estimate is a substantial undertaking (see Appendix C).

A key technical obstacle in proving Theorem~\ref{thm:main} is bounding sums that have the form 
    $\sum_{j_1,j_2,j_3,j_4=1}^p \E\big( x_{1j_1}^{k_1} x_{1j_2}^{k_2} x_{1j_3}^{k_3}  x_{1j_4}^{k_4}  \big)$ for exponents $k_1,\dots,k_4 \in \{0,2,4\}$. In order to analyze such sums, we develop a systematic approach based on \emph{Isserlis' theorem} \citep{janson1997gaussian}, also known as Wick's theorem, which is a powerful tool for calculating mixed moments of normal random vectors that is based on the combinatorial framework of pair partitions. Because the statement of Isserlis' theorem is only applicable to normal random vectors, an important merit of our analysis is that we are able to successfully adapt it to the setting of elliptical random vectors. Furthermore, even in the special case of normal data, Isserlis' theorem produces a formula for each mixed moment $\E\big( x_{1j_1}^{k_1} x_{1j_2}^{k_2} x_{1j_3}^{k_3}  x_{1j_4}^{k_4}  \big)$ that may involve an extremely large number of terms. To deal with this issue, we introduce an equivalence relation among the terms that allows us to decompose the moment calculations across a tractable number of equivalence classes. (See the proof of Proposition~\ref{prop:combinationofhighmoments} in Appendix~\ref{app:T}.)

With regard to numerical performance, we show through simulations that the proposed test has an empirical level that nearly matches the nominal level when the data are generated from elliptical distributions in 60 distinct parameter settings. Also, the test is able to reject the null hypothesis even when the data-generating distribution differs from elliptical distributions by small perturbations. For comparison, we show that the proposed test has substantially higher power than a recently developed normality test~\citep{chen2023normality}, which seems to be the only such test whose asymptotic validity has been established in high dimensions. (Such a comparison of power is justified because both tests are intended to reject in the presence of non-elliptical data.) Lastly, we present some examples involving natural datasets, illustrating that for commonly used significance levels, the proposed test is able to reject, whereas the normality test is not. Thus, our approach is of independent interest in connection with the topic of high-dimensional normality testing.\\[-0.3cm]

\noindent\textbf{Notation and conventions. }
For any real number $q\geq 1$, the entrywise $\ell_q$ norm of a real matrix $A$ is  $\|A\|_q = (\sum_{i,j } |A_{ij}|^q)^{1/q}$, and the $L^q$ norm of a scalar random variable $\zeta$ is $\|\zeta\|_{L^q}= (\E(|\zeta|^q))^{1/q}$.
For two sequences of real numbers $\{a_n\}$ and $\{b_n\}$, the relation $a_n=o(b_n)$ means $a_n/b_n\to 0$ as $n\to\infty$. When $\{a_n\}$ and $\{b_n\}$ are non-negative, we write $a_n\lesssim b_n$ if there is a constant $C>0$ not depending on $n$ such that  $a_n\leq Cb_n$ holds for all large $n$. When both of the relations $a_n\lesssim b_n$ and $b_n\lesssim a_n$ hold, we write $a_n\asymp b_n$. The maximum and minimum of a pair of real numbers $a$ and $b$ are sometimes written as $\max\{a,b\}=a\vee b$ and $\min\{a,b\}=a\wedge b$.
Convergence in probability and convergence in distribution are respectively denoted by $\xrightarrow{\P}$ and $\Rightarrow$.
Lastly, certain formulas in our testing procedure involve fractions of random variables in which the denominator may vanish.
To make the procedure well defined in any finite-sample situation, a fraction with a zero denominator is arbitrarily defined to be equal to one. This convention is of minor importance, because it is only relevant for events that occur with asymptotically negligible probability, and hence, does not affect the limiting null distribution of the proposed test statistic.

\section{Method}
\label{sec:meth}
Recall from our discussion of~\eqref{eqn:kurt} that the proposed test seeks to detect violations of the null hypothesis by measuring the discrepancy between two distinct estimates $\tilde\kappa$ and $\check{\kappa}$ for the kurtosis parameter $\kappa$. These two estimates will be constructed using separate halves of the data, $\x_1 \ldots, \x_{n/2}$ and $\x_{n/2 + 1} \ldots, \x_{n}$, and we assume henceforth that $n$ is even for simplicity. Since the left side of~\eqref{eqn:kurt} asserts that $\kappa$ is equal to the entrywise kurtosis $\E(x_{1j}^4)/\E(x_{1j}^2)^2$ for all $j=1,\dots,p$, it is natural to construct the first estimate $\tilde\kappa$ by averaging the entrywise kurtosis estimates $\tilde \kappa_j  =  \frac{1}{n/2} \sum_{i=1}^{n/2} x_{ij}^4\big/(\frac{1}{n/2} \sum_{i=1}^{n/2}  x_{ij}^2)^2 $
over $j=1,\dots,p$, which gives
\begin{equation*}
% \label{eq:definationofbarkappa}
    \TILDEKAPPA \ = \ \frac{1}{p}\sum_{j=1}^p \tilde \kappa_j.
\end{equation*}

To define the second estimate $\check{\kappa}$, it is helpful to introduce the parameters $ \varsigma^2 =  \var(\|\mathbf{x}_1\|_2^2) $ and $\nu_k  =  \tr(\SIGMA^k)$ for $k\geq 1$,
so that the right side of~\eqref{eqn:kurt} can be rewritten as
 \begin{align}\label{eq:definitionofkappa1}
    \kappa 
    \ = \ & \frac{3(\varsigma^2 + \nu_1^2)}{\nu_1^2 + 2\nu_2}.
\end{align}
From this equation, we construct $\check{\kappa}$ by substituting in estimates of $\varsigma^2$, $\nu_1$ and $\nu_2$. For estimating the latter two parameters, let $\check{\Sigma}=\frac{1}{n/2}\sum_{i=n/2+1}^{n}\x_i\x_i\ttop$, as well as $\check{\nu}_k=\tr(\check{\Sigma}^k)$ for $k\geq 1$. Specifically, we estimate $\nu_1$ using $\check{\nu}_1$, and we estimate $\nu_2$ using $\check{\nu}_2-\frac{1}{n/2}\check{\nu}_1^2$, which is a bias-corrected version of $\check{\nu}_2$. Next, we estimate $\varsigma^2$ using the sample variance of the squared norms $\|\x_i\|_2^2$ for $i=n/2+1,\dots,n$, which is 
\begin{align*}
&\check{\varsigma}^2 \ = \ \frac{1}{2\binom{n/2}{2}}\sum_{n/2<i<i' \leq n}\big(\|\x_i\|_2^2-\|\x_{i'}\|_2^2\big)^2.
\end{align*}
Hence, the second estimate of $\kappa$ is defined by
\begin{equation*}
    \check{\kappa}  \ = \  \frac{3(\check{\varsigma}^2+\check{\nu}_1^2)}{\check{\nu}_1^2 + 2(\check{\nu}_2 - \frac{2}{n} \check{\nu}_1^2)}.
\end{equation*}

To compare $\tilde\kappa$ and $\check{\kappa}$, it will be technically convenient to work with a shifted and rescaled version of the difference $\tilde\kappa-\check{\kappa}$, denoted
\begin{equation*}
    T_n \ = \ \sqrt{\frac{pn}{2}}\Big(\frac{\tilde\kappa - \check{\kappa}}{3} + \frac{4}{n}\Big),
\end{equation*}
where the term $4/n$ merely ensures that $T_n$ is asymptotically centered. Finally, the statistic $T_n$ must be normalized by an estimate of its standard deviation, which is denoted by $\hat\sigma_n$ and will be defined below in Section~\ref{sec:variance}. So altogether, the proposed test statistic is $T_n/\hat\sigma_n$, and our main theoretical result in Theorem~\ref{thm:main} will show that under the null hypothesis, it satisfies
%\begin{equation}\label{eqn:prelim}
$T_n/\hat\sigma_n  \Rightarrow  N(0,1)$
%\end{equation}
as $n\to\infty$ with $p\asymp n$. Accordingly, the proposed testing procedure rejects at a nominal level of $\alpha\in(0,1)$ if the event 
\begin{equation*}
% \label{eq:reject}
   \mathsf{R}_n(\alpha) \ = \  \Big\{|T_n|>\hat\sigma_n \,\Phi^{-1}(1-\alpha/2)\Big\}
\end{equation*}
occurs, where $\Phi$ is the standard normal distribution function.

\subsection{Estimating the variance of $T_n$}\label{sec:variance}

We will show in Proposition~\ref{prop:T} of Appendix~\ref{app:T} that the variance of $T_n$ is asymptotically equivalent to a parameter $\sigma_n^2 = \sigma_{n,1}^2+\sigma_{n,2}^2$ defined by
\begin{align}\label{eq:definitionofsigma1}
       \sigma_{n,1}^2  & \ = \ \frac{8\| R\|_4^4}{3p} \text{ \ \ \ \  \ \ \ \ \ and \ \ \ \ \ \ \ \ \ }
       \sigma_{n,2}^2   \ = \  \frac{8p(2\nu_4 +\nu_2^2)}{(\nu_1^2 +2\nu_2)^2},
\end{align}
where $R$ is the correlation matrix of $\x_1$.
To write down an estimate for $\sigma_n^2$, let $\hat\Sigma=\frac{1}{n}\sum_{i=1}^n \x_i\x_i\ttop$, as well as $\hat{\nu}_k = \tr(\hat \SIGMA^k)$ for $k\geq 1$, and $\hat R_{jj'}=\hat \Sigma_{jj'}/(\hat \Sigma_{jj}\hat \Sigma_{j'j'})^{1/2}$ for $1\leq j,j'\leq p$. Then, we define $\hat\sigma_n^2=\hat\sigma_{n,1}^2+\hat\sigma_{n,2}^2$ according to
\begin{align}\label{eqn:hatsigma1def}
   \  \ \  \hat\sigma_{n,1}^2 & 
   \ = \ \frac{8}{3p}\Big(\|\hat R\|_4^4 -\hat c\Big) \text{  \ \ \ \  \ \ \ \ \ and \ \ \ \ \ \ \ \ \  }
    \hat\sigma_{n,2}^2  \ = \ \displaystyle\frac{8p(2\hat{\nu}_4 +(\hat{\nu}_2- \frac{1}{n} \hat{\nu}_1^2)^2)}{ (\hat \nu_1^2 +2(\hat{\nu}_2- \frac{1}{n} \hat{\nu}_1^2))^2},
\end{align}
where $\hat c$ is a correction term whose definition requires some additional notation: For any real number $x$ and positive number $t$, define the threshold function $\llbracket x\rrbracket_t=\text{sgn}(x)(|x|\wedge t)$.
Also, for any fixed symmetric matrix $M\in\R^{p\times p}$ and integer $k\geq 1$, define the function $g_k(M)=\E((\z_1\ttop M\z_1)^k)$, where $\z_1$ is a standard normal vector in $\R^p$. (Explicit formulas for this function are provided in Lemma \ref{lem:momentsofz}.)
In this notation, the correction term $\hat c$ is given by
 \begin{equation}\label{eqn:hatcdef}
 \hat c \ = \  \hat\beta\|\hat R\|_4^4 - 3 \big\llbracket1-\hat\beta+\hat\gamma\big\rrbracket_{t_p}\|\hat R\|_{2}^2,
 \end{equation}
where $t_p=p^{-3/4}\log(p)$, and
\begin{align}
    \hat \beta \ = \ &1- \frac{\frac{1}{n}\sum_{i=1}^n \|\mathbf{x}_{i}\|_2^{8}}{g_4(\hat \SIGMA)  - \ts \frac{12}{n}(\hat{\nu}_1^4 +2\hat{\nu}_1^2\hat{\nu}_2 - \ts \frac{1}{n}\hat{\nu}_1^4    )} \label{eq:definitionofhatbeta}\\[0.2cm]
    \hat \gamma \ = \ & 1 + \frac{\hat{\varsigma}^2 - 2(\hat{\nu}_2- \ts\frac{1}{n} \hat{\nu}_1^2)}{g_2(\hat \SIGMA)  - \ts\frac{2}{n} \hat{\nu}_1^2} - \frac{\frac{1}{n}\sum_{i=1}^n \|\mathbf{x}_{i}\|_2^{6}}{\frac{1}{2}g_3(\hat \SIGMA) - \ts\frac{3}{n} \hat{\nu}_1^3}.
\end{align}
The definition of the correction term in~\eqref{eqn:hatcdef} is motivated by an asymptotic formula for $\sigma_{n,1}^2$ that results from~\eqref{eqn:correction_motivation} and~\eqref{eqn:sig1equiv} in the proof of Lemma~\ref{lem:NullI_n}. In particular, it can be shown that the term $ 3 \big\llbracket1-\hat\beta+\hat\gamma\big\rrbracket_{t_p}\|\hat R\|_{2}^2$ is asymptotically negligible, but it turns out that retaining this term is beneficial at a finite-sample level.
Finally, when interpreting the formulas above, recall that in the unlikely event where the denominator of a fraction vanishes,  we arbitrarily define the fraction to be equal to one.

\section{Theory}
Our analysis will be done with a sequence of models implicitly embedded in a triangular array, where the rows are indexed by $n$, and all parameters may vary with $n$ except when stated otherwise. The following assumption specifies the sequence of null hypotheses in~\eqref{eqn:gof} corresponding to different values of $n$. 

\begin{assumption}[]
\label{Data generating model}
As $n\to\infty$, the dimension $p$ grows such that $p\asymp n$. For each $i=1,\dots,n$, the observation $\x_i\in\R^p$ has the form
  \begin{equation*}
      \x_i \ = \ \xi_i\SIGMA^{1/2}\u_i,
  \end{equation*}
  where $\SIGMA\in\R^{p\times p}$ is deterministic and positive semidefinite with positive diagonal entries, and $(\u_1,\xi_1),\dots,$ $(\u_n,\xi_n)$ are i.i.d.~random vectors in $\R^p\times [0,\infty)$ satisfying the following conditions: The vector $\u_1$  is drawn from the uniform distribution on the unit sphere of $\R^p$, and is independent of $\xi_1$. In addition, the random variable $\xi_1^2$ satisfies $\E(\xi_1^2)=p$, as well as the conditions
 \begin{align} \label{eqn:moment assumption}
   \var\!\bigg(\frac{\xi_1^2-p}{\sqrt p}\bigg)   \ = \   \tau +o(1) \ \ \ \ \ \  \text{ and} \ \ \ \ \ \  \  
   \bigg\|\frac{\xi_1^2-p}{\sqrt p}\bigg\|_{L^8}\, = \, o(p^{1/4})
  \end{align}
as $n\to\infty$, for some constant $\tau \geq 0$ that is fixed with respect to $n$.
\end{assumption} 
\noindent\textbf{Remarks.} Perhaps the most striking feature of Assumption~\ref{Data generating model} is that the population covariance matrix $\SIGMA$ is unrestricted. (The condition that $\SIGMA$ has positive diagonal entries does not impose any real restriction, because if $\SIGMA_{jj}=0$, then the $j$th covariate is zero almost surely and is hence irrelevant.) Also, the model only places conditions on a finite number of moments. Similar moment assumptions have been considered elsewhere, such as in the papers~\citep{hu2019aos,hu2019ieee}, which analyze eigenvalue statistics of sample covariance matrices in high-dimensional elliptical models. In essence, the purpose of the two conditions in~\eqref{eqn:moment assumption} is to make it possible to precisely determine the magnitudes of higher order remainder terms in the proofs of certain CLTs underlying Theorem~\ref{thm:main} below. It should be noted too that the second condition in \eqref{eqn:moment assumption} allows for the eighth moment of $(\xi_1^2-p)/\sqrt p$ to diverge with the dimension, i.e.~$\|(\xi_1^2-p)/\sqrt p\|_{L^8}\to\infty$ as $p\to\infty$. 
In fact, this condition can be weakened to allow for even faster divergence as $p\to\infty$, but it would involve a somewhat more complicated statement of Assumption~\ref{Data generating model}, and so we have retained the version above to keep the statement as concise as possible.

\noindent\textbf{Examples.} Below, we provide an assortment of examples for the distribution of $\xi_1^2$ that conform with Assumption~\ref{Data generating model}. In each case, the distribution is parameterized so that the basic normalization condition $\E(\xi_1^2)=p$ holds.
\begin{itemize}
    \item 
    Chi-Squared on $p$ degrees of freedom,
    %\\[-0.2cm] %, \quad $\tau = 2$ 
    \item Poisson($p$),
    %\\[-0.2cm]%, \quad $\tau = 1$
    \item $(1-\tau)$Negative-Binomial$(p,1-\tau)$, \ for any fixed $\tau \in (0,1)$,
    %\\[-0.2cm]
    % 
    \item $(p+2b)$\,Beta$(p/2,b)$, \ for any fixed $b>0$,
    \item Gamma$(p/\tau,1/\tau)$, \ for any fixed $\tau>0$,
    \item Beta-Prime$\left(\frac{p(1+p+\tau)}{\tau},\frac{1+p+2\tau}{\tau}\right)$, \ for any fixed $\tau >0$,
    \item Log-Normal$\left(\log(p)-\frac{1}{2}\log\big(1+\frac{\tau}{p}\big), \log\big(1+\frac{\tau}{p}\big) \right)$, \ for any fixed $\tau > 0$.
\end{itemize}
Several of these examples for $\xi_1^2$ can be subsumed within a more general non-parametric class of distributions that satisfy Assumption~\ref{Data generating model}. This class corresponds to choices of $\xi_1^2$ that can be represented in the form $\xi_1^2= \sum_{j=1}^p \eta_j$ for some independent random variables $\eta_1,\dots,\eta_p\geq 0$ satisfying 
 \begin{equation*}
 % \label{eqn:suff}
 \small
         \ts\frac{1}{p}\sum_{j=1}^p \E(\eta_j)=1, \quad \ts\frac{1}{p}\sum_{j=1}^p\var(\eta_j) = \tau+o(1), \text{\quad and \quad } \displaystyle\max_{1\leq j\leq p}\E(\eta_j^8)\,=\, o(p^{5}),
 \end{equation*}
again with $\tau\geq 0$ being a constant not depending on $n$. (This can be verified using Rosenthal's inequality, as given in Lemma~\ref{lem:Rosenthal} of the supplement.) It is also notable that the moments $\E(\eta_j^8)$ are allowed to grow quite rapidly as $p\to\infty$.

The following theorem is our main theoretical result.
\begin{theorem}
\label{thm:main}
    If Assumption \ref{Data generating model} holds, then as $n\to\infty$ we have
    \begin{equation*}
         \frac{T_n}{\hat\sigma_n} \ \Rightarrow \ N(0,1),
    \end{equation*}
    and hence, for any fixed $\alpha\in (0,1)$, the level of the test satisfies
    \begin{equation*}
        \P(\mathsf{R}_n(\alpha)) \ \to \ \alpha.
    \end{equation*}
\end{theorem}
\noindent\textbf{Remarks.} At a high level, the proof of Theorem~\ref{thm:main} involves several main parts. The first two are given in Propositions~\ref{prop:T} and~\ref{prop:rho} of the supplement, which respectively show that $T_n/\sigma_n\Rightarrow N(0,1)$ and $\hat\sigma_n/\sigma_n\xrightarrow{\P}1$ as $n\to\infty$. The problem of showing that $T_n/\sigma_n$ is asymptotically normal is then divided into two more substantial parts, corresponding to separate CLTs for standardized versions of the estimates $\tilde\kappa$ and $\check{\kappa}$, which utilize different sets of theoretical tools. 

The CLT for $\tilde\kappa$ involves analyzing a sum of i.i.d.~random variables $\sum_{i=1}^{n/2} w_i$, where $w_i=\frac{1}{\sqrt{p}}\sum_{j=1}^p \ts\frac{1}{3}(x_{ij}^4 - \E(x_{ij}^4)) - 2 (x_{ij}^2-\E(x_{ij}^2))$. We establish asymptotic normality by showing that a Lyapunov condition of the form  $\E(w_1^4)=o(n\var(w_1)^2)$ holds, which in turn, requires bounding sums of the form $\sum_{j_1,j_2,j_3,j_4=1}^p \E\big( x_{1j_1}^{k_1} x_{1j_2}^{k_2} x_{1j_3}^{k_3}  x_{1j_4}^{k_4}  \big)$ for exponents $k_1,\dots,k_4 \in \{0,2,4\}$. 
%(Note that the control of such moments is the primary reason for the second condition in~\eqref{eqn:moment assumption}.) 
Although the random vector $\x_1$ is not necessarily normal under the null hypothesis, we are able to control its mixed moments by adapting Isserlis' theorem~\citep{janson1997gaussian}, which states that the following formula holds for any dimension $d$ and any centered normal vector $(g_1,\dots,g_d)$,
\begin{align}\label{eqn:firstisserlis}
    \E(g_1\cdots g_d) \ = \ \sum_{\pi \in \mathcal{P}_d} \prod_{\{j,j'\}\in \pi} \E(g_jg_{j'}),
\end{align}
where $\mathcal{P}_d$ is the collection of all partitions $\pi$ of $\{1,\dots,d\}$ into distinct pairs $\{j,j'\}$. Here, it is important to note that even if the vector $\x_1$ was assumed to be normal, a brute-force application of Isserlis' theorem would not provide a workable solution to establishing the Lyapunov condition, because of the enormous number of terms on the right side of~\eqref{eqn:firstisserlis}. For instance, in the special case when $\x_1$ is normal and $k_1=\cdots=k_4=4$, there are more than 2 million terms in the formula for $\E( x_{1j_1}^{k_1} x_{1j_2}^{k_2} x_{1j_3}^{k_3}  x_{1j_4}^{k_4} )$, since $d=\sum_{l=1}^4 k_l$ and $\textup{card}(\mathcal{P}_{16})=15\cdot 13 \cdot 11 \cdots 1$. Nevertheless, we are able to overcome this difficulty by taking advantage of symmetries among the partitions. To explain the how this works, let $t_{\pi}(j_1,j_2,j_3,j_4)$ denote the term corresponding to the partition $\pi$ in the formula for $\E( x_{1j_1}^{4} x_{1j_2}^{4} x_{1j_3}^{4}  x_{1j_4}^{4} )$ so that interchanging the order of summation gives
\begin{equation*}
    \sum_{j_1,j_2,j_3,j_4=1}^p \E\big( x_{1j_1}^{4} x_{1j_2}^{4} x_{1j_3}^{4}  x_{1j_4}^{4}  \big) \ = \ \sum_{\pi\in\mathcal{P}_{16}} \sum_{j_1,j_3,j_3,j_4=1}^p t_{\pi}(j_1,j_2,j_3,j_4).
\end{equation*}
If we define an equivalence relation on $\mathcal{P}_{16}$ by saying that $\pi$ and $\pi'$ are equivalent if $t_{\pi}(j_1,j_2,j_3,j_4)$ and $t_{\pi'}(j_1,j_2,j_3,j_4)$ produce the same value when summed over all $j_1,j_2,j_3,j_4$, then remarkably, it turns out that there is only a handful of equivalence classes. Hence, the sum over $\mathcal{P}_{16}$ that contains more than 2 million terms can be reduced to a tractable calculation, making it possible to establish the asymptotic normality of $\tilde \kappa$. (Strictly speaking, the technical aspects of this strategy are implemented somewhat differently in the proof of Proposition~\ref{prop:combinationofhighmoments}, and there are additional considerations required to find the sizes of the different equivalence classes, but this explanation conveys the main idea.) 

The asymptotic normality of the second estimate $\check{\kappa}$ is proven in an entirely different way, by adapting classical results on U-statistics to the high-dimensional setting. For this purpose, 
Propositions~\ref{prop:T2smallr} and~\ref{prop:T2larger} show that the behavior of $\check{\kappa}$ is largely determined by the statistic $U=\check{\varsigma}^2 - 2(\check{\nu}_2 - \frac{2}{n} \check{\nu}_1^2)$, which is shown to be a U-statistic in the proof of Lemma~\ref{lem:UCLT}. There are two aspects of the classical theory of U-statistics that break down if they are used to prove a CLT for $U$ in our setup. The first is in showing that $U$ is asymptotically equivalent to its H\'ajek projection $\hat U$, and the second is in showing that $\hat U$ is asymptotically normal. In essence, the issue in the high-dimensional setting is that we must account for the fact that $\var(U)$ and $\var(\hat U)$ may vary in their orders of magnitude as $p\to\infty$, and this is accomplished in Lemmas~\ref{lem:UCLT}-\ref{lem:psi1} through a delicate analysis of the covariances of random quadratic forms.

\section{Numerical experiments}

In this section, we study the level and power of the proposed test under a substantial range of simulation settings. In addition, we offer comparisons with a high-dimensional normality test recently developed in~\citep{chen2023normality}, as well as some examples involving natural datasets.

\subsection{Assessment of level}\label{sec:level}
\noindent\textbf{Design of experiments.} For the simulations under the null hypothesis, we generated data from elliptical models based on five choices of the distribution of the squared random variable $\xi_1^2$, and four choices of the covariance matrix $\Sigma$. The choices for $\xi_1^2$ satisfy $\E(\xi_1^2)=p$ and are given by:

\begin{enumerate}[label=(\roman*).]
    \item Chi-Squared distribution with $p$ degrees of freedom, 
    \item Beta-Prime$\left(\frac{p(1+p+3)}{3},\frac{1+p+6}{3}\right)$,
    \item $(p+4)\textup{Beta}(p/2,2)$,
    \item $\textup{Gamma}(p/5, 1/5)$,
    \item $\frac{1}{p+1}\textup{Gamma}(p,1)^2$.
\end{enumerate}
\noindent The choices for the covariance matrix $\Sigma$ are:
\begin{enumerate}[leftmargin=*,label=\arabic*.]
    \item[(1).] (Spiked eigenvalues with generic eigenvectors). The eigenvalues of $\SIGMA$ are $\lambda_{j}(\SIGMA)=5$ for $j=1,\dots,5$ and $\lambda_j(\SIGMA)=1$ for $j=6,\dots,p$. The $p\times p$ matrix of eigenvectors of $\SIGMA$ is drawn from the uniform distribution on orthogonal matrices.
    \item[(2).] (Toeplitz). The matrix $\SIGMA$ has entries of the form $\SIGMA_{ij} = \rho^{|i-j|}$ with $\rho=0.1$.
    \item[(3).] (Decaying eigenvalues with generic eigenvectors). The eigenvalues of $\SIGMA$ are  \smash{$\lambda_j(\SIGMA)=j^{-1/4}$,} and the eigenvectors are the same as in case (1).
    \item[(4).] The covariance matrix $\SIGMA$ is the identity matrix.
\end{enumerate}
The dimension was selected from $p \in\{ 200, 400, 600\}$, and the sample size was chosen to be $n=400$.
In each possible setting of the tuple $(\xi_1^2,\Sigma,n,p)$, we generated 10000 datasets, and we applied the proposed test at a 5\% nominal level to each dataset. Hence, in total, the experiments presented here cover 60 different settings of the null hypothesis. Furthermore, in Appendix~\ref{app:dim}, we present additional results in cases where the dimension is even larger, for values of the ratio $p/n$ as large as 10.

\noindent\textbf{Discussion of results.} In Table~\ref{table:level}, we report the rejection rate (i.e.~empirical level) based on the 10000 trials under each elliptical distribution. For example, the entry in the lower left corner of Table~\ref{table:level} shows that in settings (v) and (1) with $p/n=0.5$, the null hypothesis was rejected in 4.52\% of the trials.

\begin{table}[H]
% \small
\centering
\setlength\tabcolsep{6.5pt}
\caption{Results on empirical level for the proposed test (5\% nominal level).}
\begin{tabular}{lllllllllllllll}
\hline
      &\multicolumn{4}{c}{$p/n=0.5$} &  & \multicolumn{4}{c}{$p/n=1$} &  & \multicolumn{4}{c}{$p/n=1.5$} \\ \cline{2-5} \cline{7-10} \cline{12-15} 
      \specialrule{0em}{1pt}{1pt}
 & (1)   & (2)   & (3)   & (4) &  & (1)   & (2)   & (3)   & (4) &  & (1)   & (2)   & (3)   & (4) \\
\hline
\specialrule{0em}{1pt}{1pt}
(i) & 4.06&4.05 &3.94 &4.00 &  &3.92&4.05 &4.12 &4.04 &  &3.91&4.05 &3.91 &3.99        \\
 \specialrule{0em}{1pt}{1pt}
 \hline
\specialrule{0em}{1pt}{1pt}
(ii) & 4.26&3.96 &3.96 &4.01 &  &4.04&4.04 &4.16 &3.93 &  &3.78&3.79 &4.00 &3.82        \\
 \specialrule{0em}{1pt}{1pt}
 \hline
\specialrule{0em}{1pt}{1pt}
(iii) & 4.20&4.16 &4.31 &4.17 &  &4.49&4.23 &4.38 &4.32 &  &4.11&4.16 &3.96 &4.17 \\
 \specialrule{0em}{1pt}{1pt}
 \hline
\specialrule{0em}{1pt}{1pt}
(iv) & 5.27&4.79 &4.78 &4.61 &  &4.19&4.25 &4.28 &4.28 &  &4.11&4.15 &3.92 &4.14        \\
 \specialrule{0em}{1pt}{1pt}
 \hline
\specialrule{0em}{1pt}{1pt}
(v) & 4.52&4.43 &4.71 &4.25 &  &4.03&4.10 &4.22 &4.11 &  &3.95&4.12 &3.99 &3.98         \\
 \specialrule{0em}{1pt}{1pt}
 \hline
\end{tabular}
\label{table:level}
\end{table}

% \begin{table}[H]
% \centering
% \setlength\tabcolsep{10pt}
% \caption{Results on empirical level for the proposed test (5\% nominal level)}
% \begin{tabular}{lllllllllllll}
% %
% \hline
%       &  & \multicolumn{3}{c}{$p/n=0.5$} &  & \multicolumn{3}{c}{$p/n=1$} &  & \multicolumn{3}{c}{$p/n=1.5$} \\ \cline{3-5} \cline{7-9} \cline{11-13} 
%       \specialrule{0em}{1pt}{1pt}
%  &  & (1)   & (2)   & (3) &  & (1)   & (2)   & (3) &  & (1)   & (2)   & (3) \\
% \hline
% \specialrule{0em}{1pt}{1pt}
%  %
% (i) & &3.85&4.05 &4.05 &  &4.09&4.05 &4.12 &  &3.9&4.05 &3.91         \\
% %
%  \specialrule{0em}{1pt}{1pt}
%  \hline
% \specialrule{0em}{1pt}{1pt}
%  %
% (ii) & &4.54&3.96 &4.05 &  &3.87&4.04 &4.16 &  &3.8&3.79 &4.00        \\
%  %
%  \specialrule{0em}{1pt}{1pt}
%  \hline
% \specialrule{0em}{1pt}{1pt}
%  % 
% (iii) & &4.43&4.16 &4.31 &  &4.08&4.23 &4.38 &  &3.86&4.16 &3.96 \\
%  \specialrule{0em}{1pt}{1pt}
%  \hline
% \specialrule{0em}{1pt}{1pt}
%  %
% (iv) & &4.99&4.79 &4.78 &  &4.21&4.25 &4.28 &  &3.77&4.15 &3.92         \\
%  %
%  \specialrule{0em}{1pt}{1pt}
%  \hline
% \specialrule{0em}{1pt}{1pt}
%  %   
% (v) & &4.64&4.43 &4.71 &  &4.13&4.1 &4.09 &  &3.91&4.12 &3.99           \\
%  %
%  \specialrule{0em}{1pt}{1pt}
%  \hline
% \end{tabular}
% \label{table:level}
% \end{table}
% %
In all settings, the empirical level is within approximately 1\% of the nominal level. It is also notable that the difference between the empirical and nominal levels is reliably in the conservative direction, which is often preferable to a situation where the direction of disagreement is unpredictable or anti-conservative. Lastly, a comparison of the empirical level of our method and that of~\citep{chen2023normality} is implicitly contained in the rejection rate curves at $h=0$ in Figures~\ref{fig:lap&0.5}-\ref{fig:beta&1.5}, and this comparison is discussed below.

\subsection{Assessment of power} \label{sec:power}
\noindent\textbf{Design of experiments.} To study the power of the proposed test, we generated data under the alternative hypothesis using models of the form
$\x_1 = \SIGMA^{1/2}\mathbf{s}_1$, where $\mathbf{s}_1\in\R^p$ is a random vector whose entries are i.i.d random variables satisfying $\E(s_{11}) = 0$ and $\E(s_{11}^2) = 1$. For the choice of $\Sigma$, we used the same menu of options (1)-(4) described in Section~\ref{sec:level}. The ``distance'' between the null and alternative hypotheses is controlled through a parameter $h\in[0,1]$, such that the entries of $\mathbf{s}_1$ are given by
\begin{align*}
    s_{11}  \ = \  \sqrt{1- h} \, z_{11} +\sqrt{h}\, y_{11},
\end{align*}
where $z_{11}$ is a standard normal random variable and $y_{11}$ is a standardized (a) Laplace or (b) Beta random variable. In detail, the choices of $y_{11}$ correspond to
\begin{enumerate}
    \item[(a).] $y_{11} \sim \frac{\textup{Laplace}(0,1)}{\sqrt{2}}$,
    % $k = 3h$;
    \item[(b).] $y_{11} \sim  \frac{\textup{Beta}(2,3/2) - 4/7}{\sqrt{8/147}}$.
\end{enumerate}
Note that if $h=0$, then $\x_1$ has a $N(0,\Sigma)$ distribution, which is elliptical, and so the null hypothesis holds. By contrast, if $0<h\leq 1$, then the alternative hypothesis holds. Another reason for generating the data in this manner is that it allows for convenient comparisons of both level and power with respect to the high-dimensional normality test proposed in~\citep{chen2023normality}. Such a comparison is legitimate, because both tests are intended to accept the  null hypothesis when $h=0$ and reject it when $h>0$.

As in Section~\ref{sec:level}, the dimension was selected from $p\in\{200,400,600\}$ and the sample size was chosen to be $n=400$. For each possible choice of $(y_{11},\Sigma,n,p)$, we generated 1000 datasets corresponding to several values of $h$ in a grid within $[0,1]$, and we performed both the proposed test and the normality test at a 5\% nominal level on every dataset. (The normality test also involves two  parameters, denoted $L$ and $B$, which we set to the values $L=1$ and $B=500$ suggested in~\citep{chen2023normality}.) Finally, at each setting of $(y_{11},\Sigma,n,p,h)$, we recorded the rejection rate of each test over the 1000 datasets, allowing the results to be visualized as curves depending on $h$, such that their heights at $h=0$ represent the empirical level, and their heights at $h\in(0,1]$ represent the empirical power.

\begin{figure}[H]
\setlength{\abovecaptionskip}{20pt}
\centering 
\subfigure{
\begin{overpic}[width=0.21\textwidth]{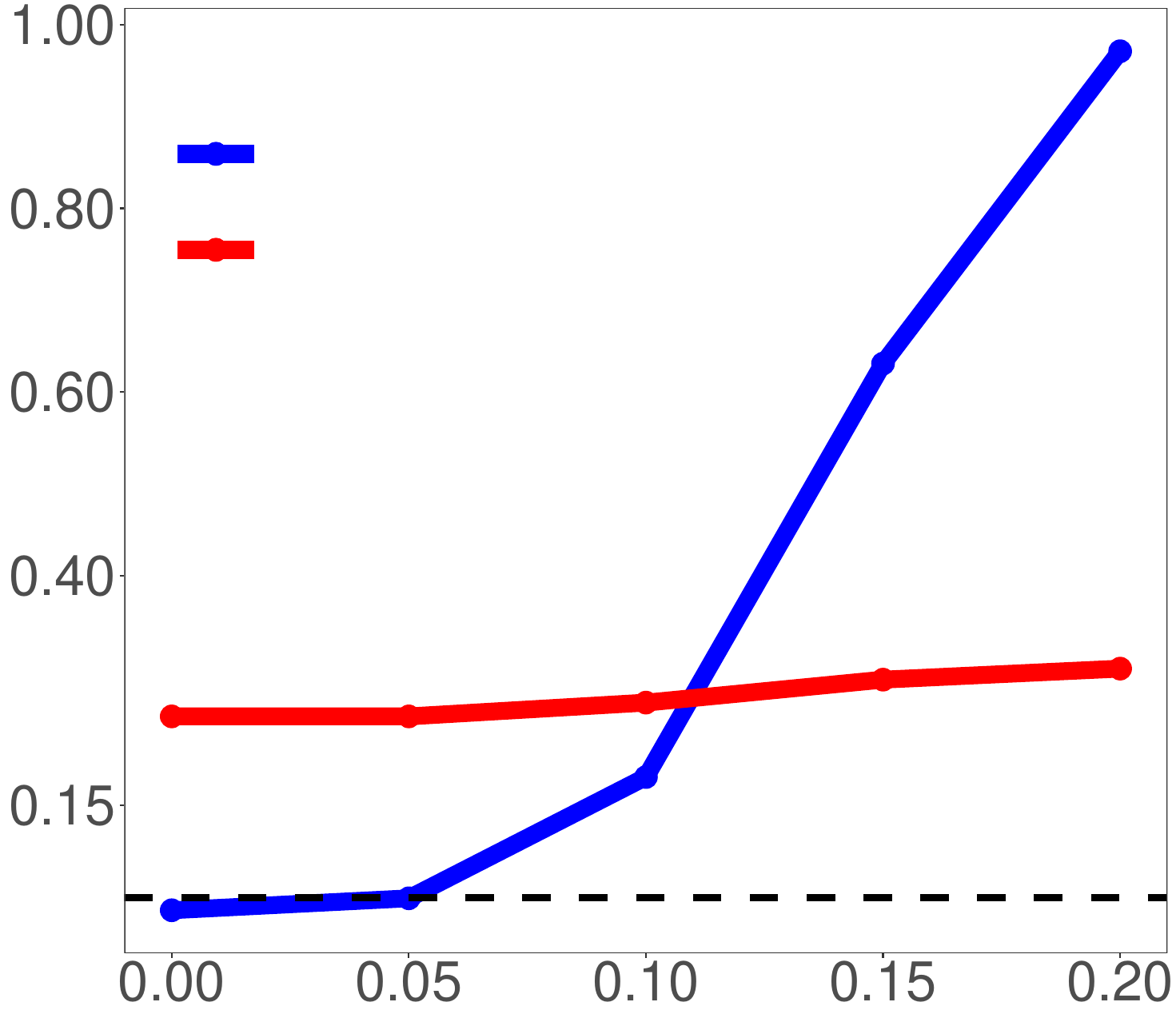}
		\put(-10,-4){\rotatebox{90}{ {  \ \ \ \footnotesize rejection rate  \ \ }}}
\put(50,-9){ \footnotesize $h$  }
\put(27,-22){ Model (1) }
\put(20,63){ \scriptsize normality test}
\put(20,73){ \scriptsize proposed test}
	\end{overpic}}
\hspace{4mm}
\subfigure{
\begin{overpic}[width=0.21\textwidth]{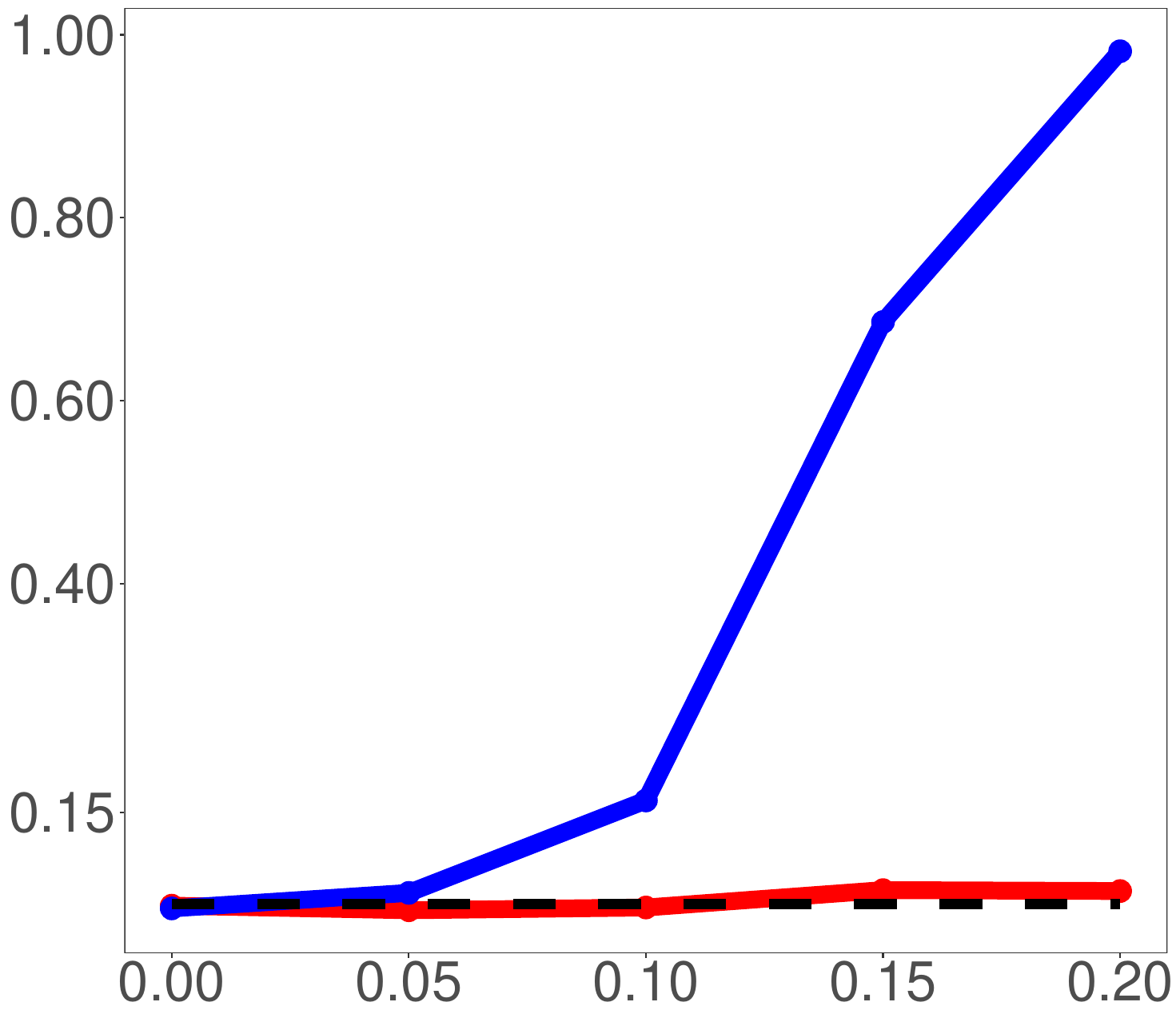}
		\put(-10,-4){\rotatebox{90}{ {  \ \ \ \footnotesize rejection rate  \ \ }}}
\put(50,-9){ \footnotesize $h$  }
\put(27,-22){ Model (2) }
	\end{overpic}}
\hspace{4mm}
\subfigure{
\begin{overpic}[width=0.21\textwidth]{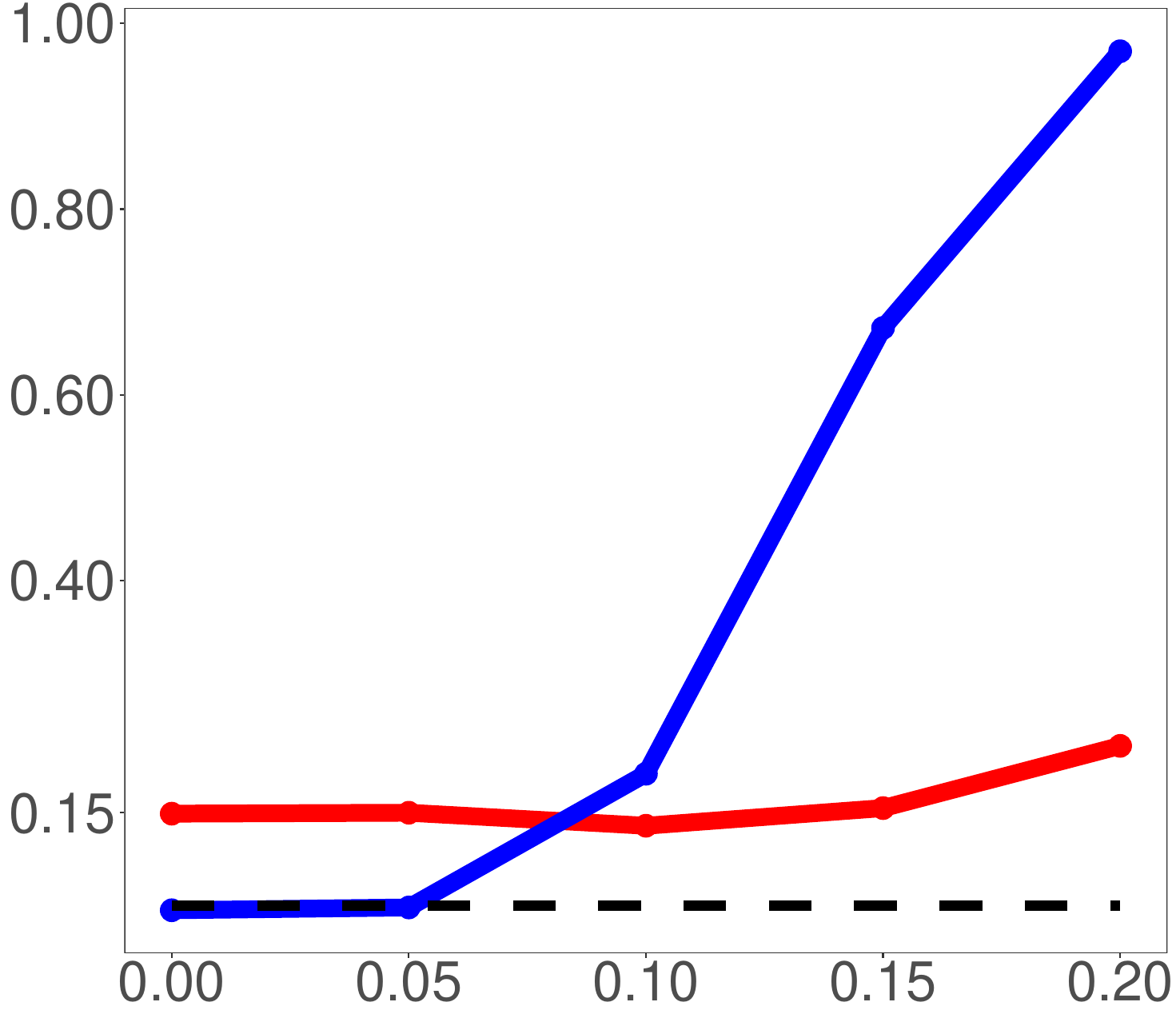}
		\put(-10,-4){\rotatebox{90}{ {  \ \ \ \footnotesize rejection rate  \ \ }}}
\put(50,-9){ \footnotesize $h$  }
\put(27,-22){ Model (3) }
	\end{overpic}}
 \hspace{4mm}
\subfigure{
\begin{overpic}[width=0.21\textwidth]{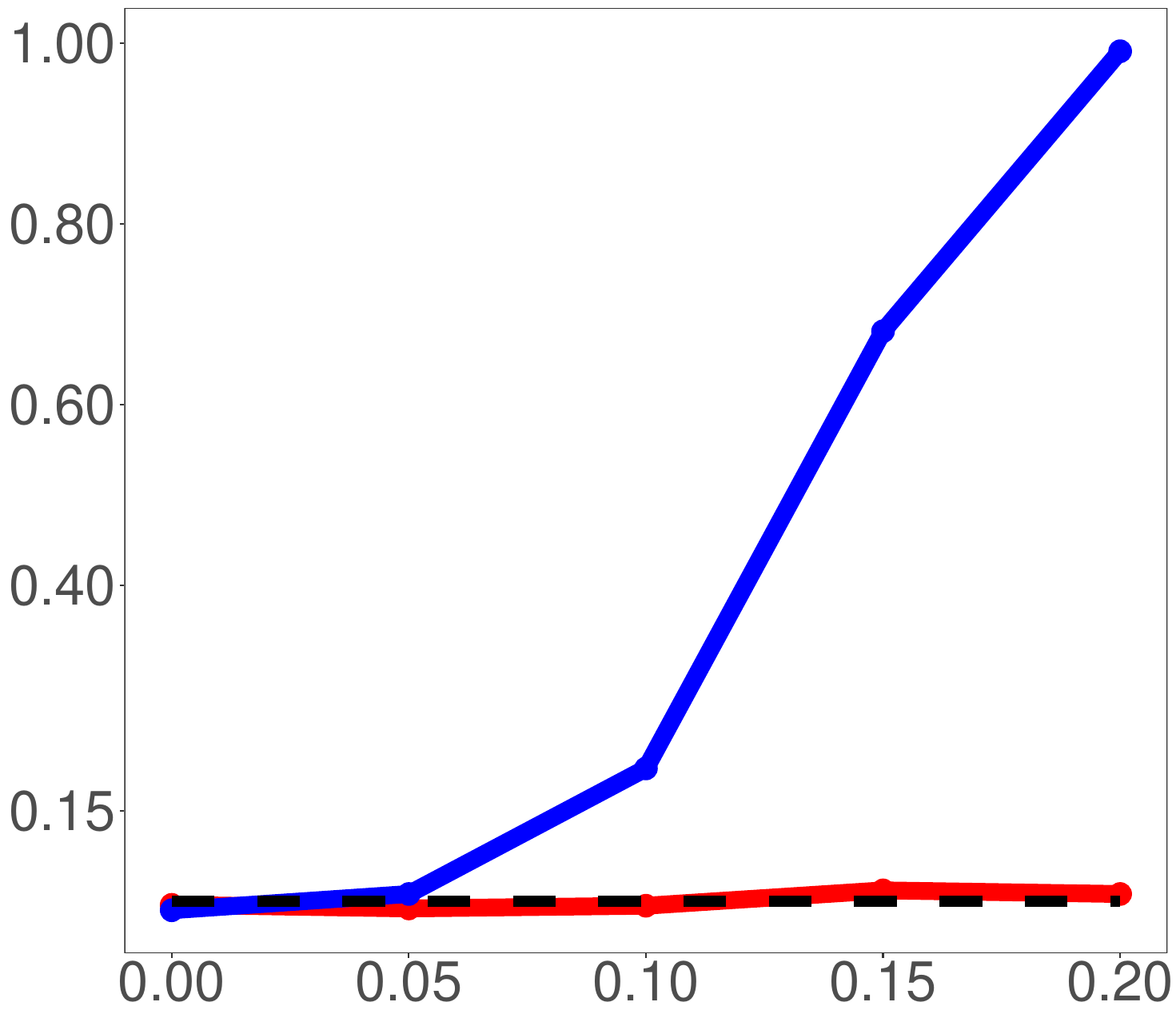}
            \put(-10,-4){\rotatebox{90}{ {  \ \ \ \footnotesize rejection rate  \ \ }}}
\put(50,-9){ \footnotesize $h$  }
\put(27,-22){ Model (4) }
	\end{overpic}}
 % \vspace{0.2cm}
\caption{ Results for $y_{11}\sim \frac{\textup{Laplace}(0,1)}{\sqrt{2}}$ when $p/n = 0.5$ and $p=200$.}
\label{fig:lap&0.5}
\end{figure}

\begin{figure}[H]
\setlength{\abovecaptionskip}{20pt}
\centering 
\subfigure{
\begin{overpic}[width=0.21\textwidth]{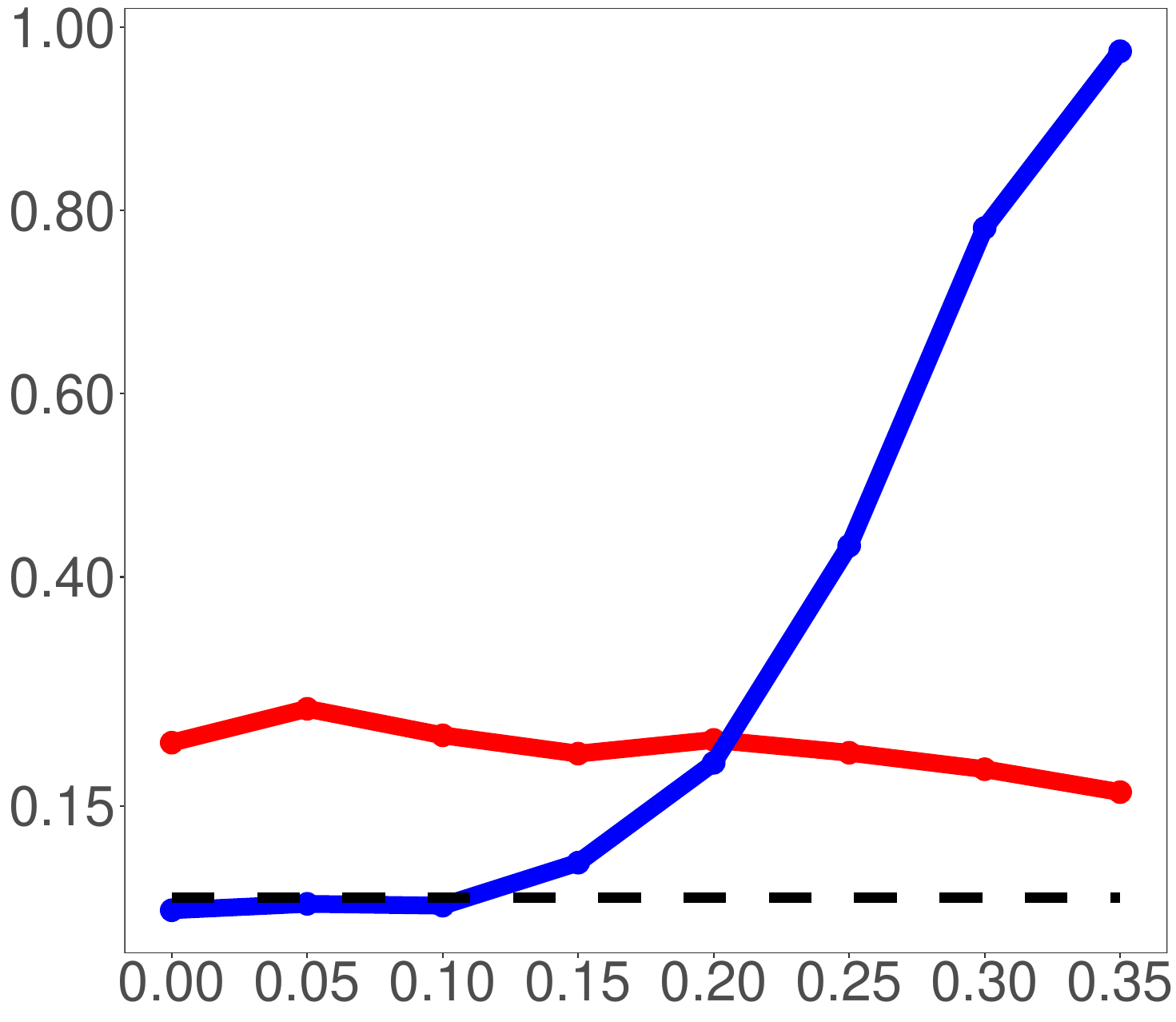}
		\put(-10,-4){\rotatebox{90}{ {  \ \ \ \footnotesize rejection rate  \ \ }}}
\put(50,-9){ \footnotesize $h$  }
\put(27,-22){ Model (1) }
	\end{overpic}}
\hspace{4mm}
\subfigure{
\begin{overpic}[width=0.21\textwidth]{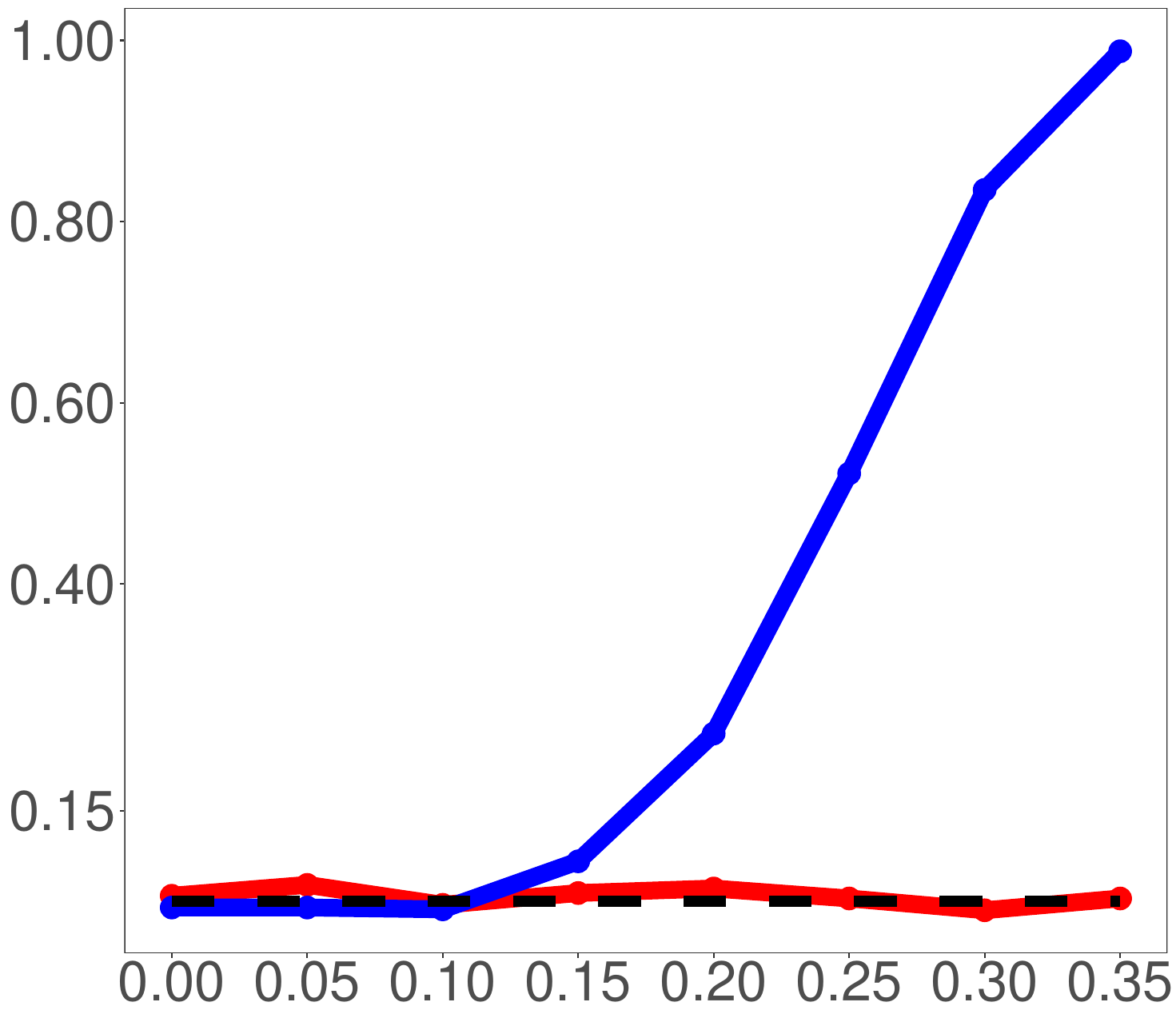}
		\put(-10,-4){\rotatebox{90}{ {  \ \ \ \footnotesize rejection rate  \ \ }}}
\put(50,-9){ \footnotesize $h$  }
\put(27,-22){ Model (2) }
	\end{overpic}}
\hspace{4mm}
\subfigure{
\begin{overpic}[width=0.21\textwidth]{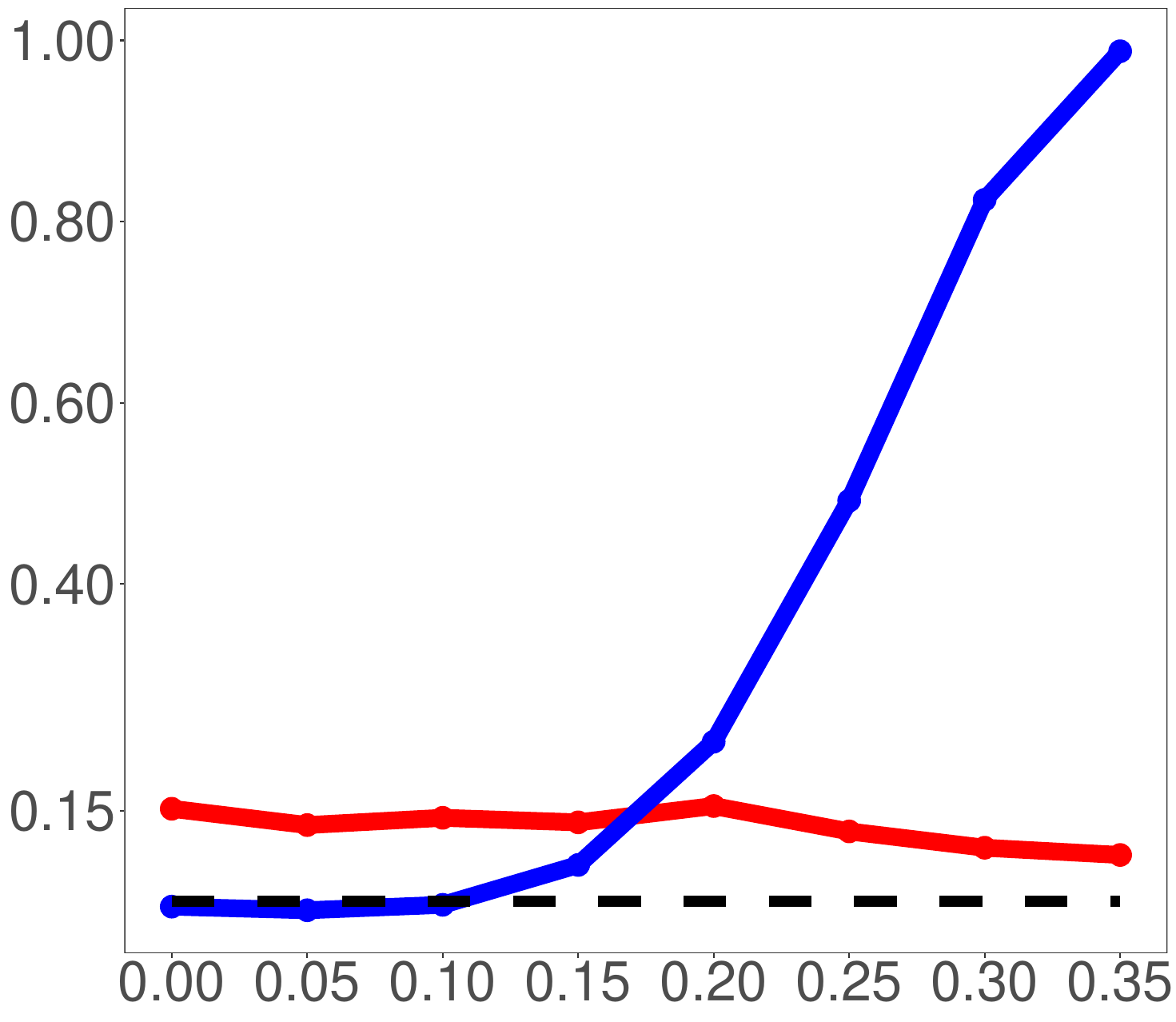}
		\put(-10,-4){\rotatebox{90}{ {  \ \ \ \footnotesize rejection rate  \ \ }}}
\put(50,-9){ \footnotesize $h$  }
\put(27,-22){ Model (3) }
	\end{overpic}}
\hspace{4mm}
\subfigure{
\begin{overpic}[width=0.21\textwidth]{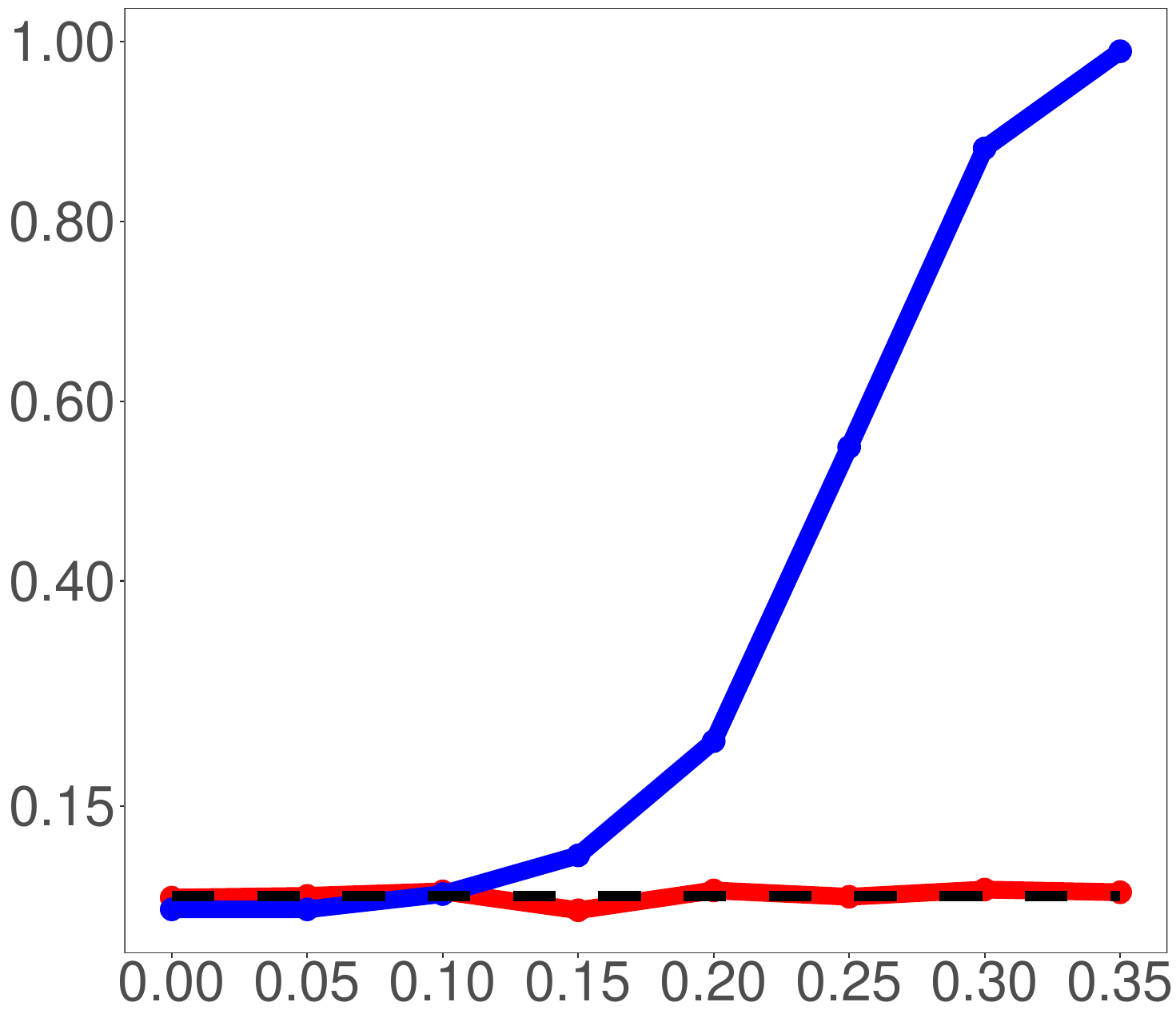}
		\put(-10,-4){\rotatebox{90}{ {  \ \ \ \footnotesize rejection rate  \ \ }}}
\put(50,-9){ \footnotesize $h$  }
\put(27,-22){ Model (4) }
	\end{overpic}} 
 % \vspace{0.2cm}
\caption{Results for $y_{11}\sim \frac{\textup{Beta}(2,1.5) - 4/7}{\sqrt{8/147}}$ when $p/n = 0.5$ and $p=200$.}
\label{fig:beta&0.5}
\end{figure}

\begin{figure}[H]
\setlength{\abovecaptionskip}{20pt}
\centering 
\subfigure{
\begin{overpic}[width=0.21\textwidth]{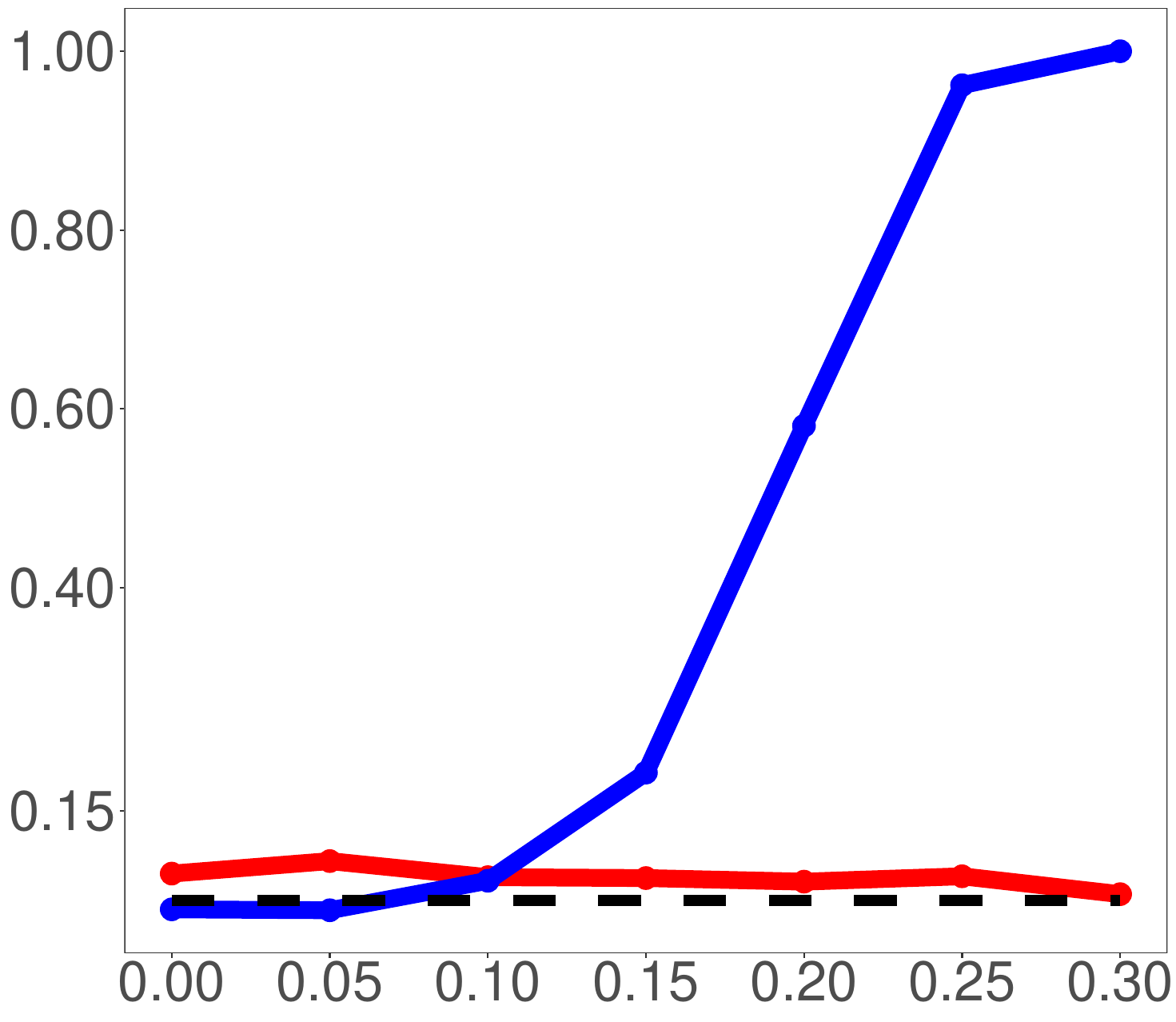}
		\put(-10,-4){\rotatebox{90}{ {  \ \ \ \footnotesize rejection rate  \ \ }}}
\put(50,-9){ \footnotesize $h$  }
\put(27,-22){ Model (1) }
	\end{overpic}}
\hspace{4mm}
\subfigure{
\begin{overpic}[width=0.21\textwidth]{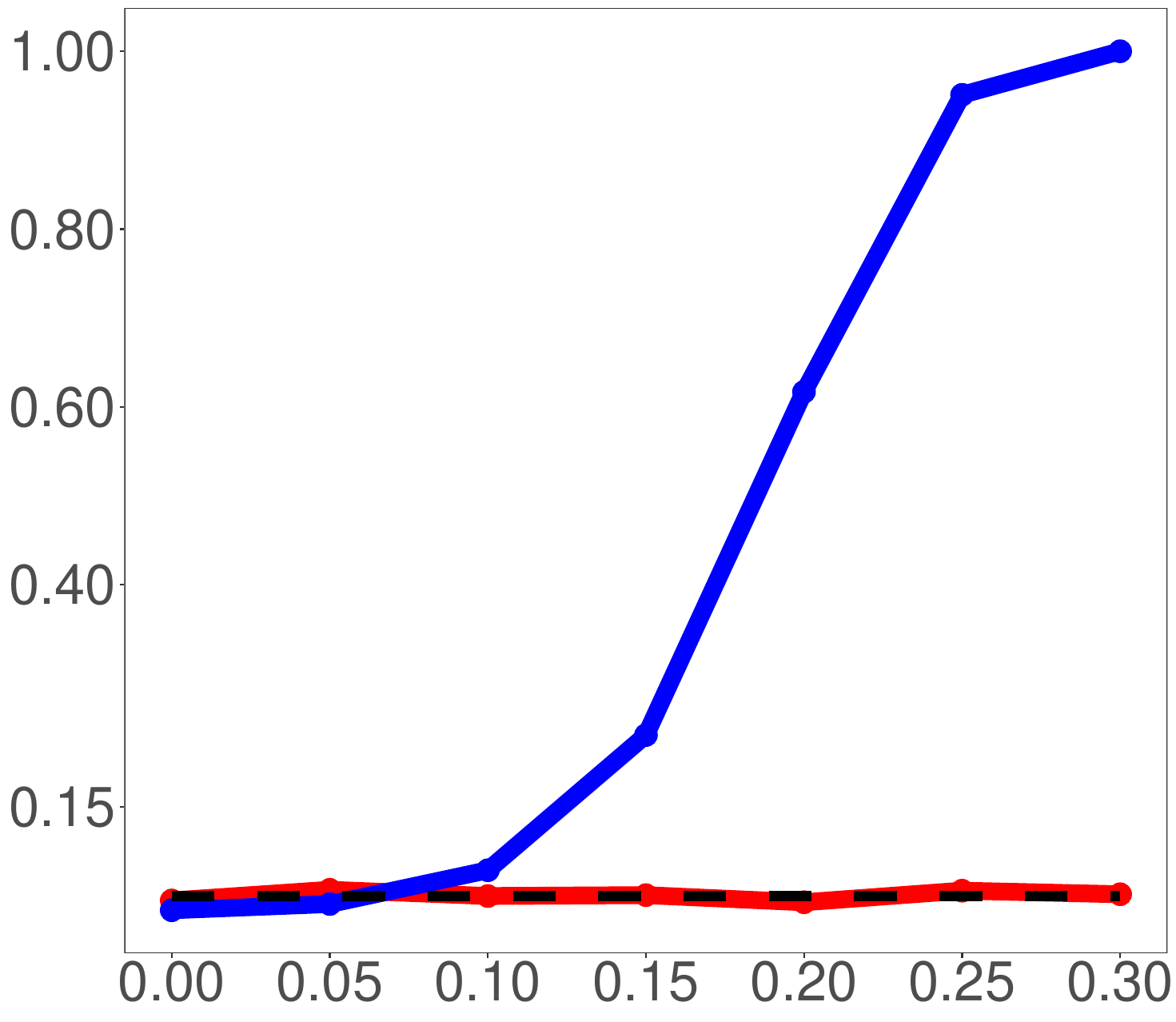}
		\put(-10,-4){\rotatebox{90}{ {  \ \ \ \footnotesize rejection rate  \ \ }}}
\put(50,-9){ \footnotesize $h$  }
\put(27,-22){ Model (2) }
	\end{overpic}}
\hspace{4mm}
\subfigure{
\begin{overpic}[width=0.21\textwidth]{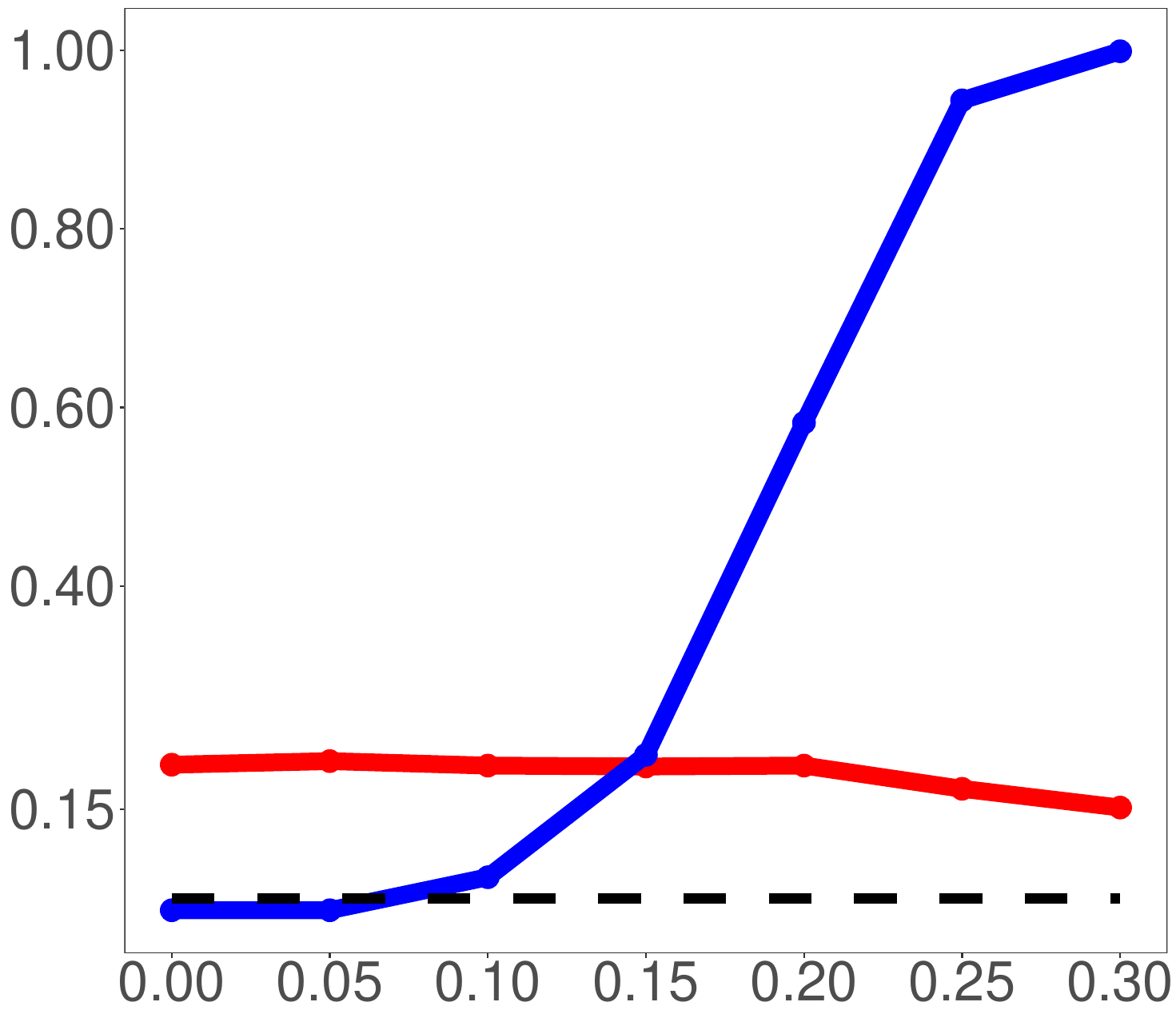}
		\put(-10,-4){\rotatebox{90}{ {  \ \ \ \footnotesize rejection rate  \ \ }}}
\put(50,-9){ \footnotesize $h$  }
\put(27,-22){ Model (3) }
	\end{overpic}}
\hspace{4mm}
\subfigure{
\begin{overpic}[width=0.21\textwidth]{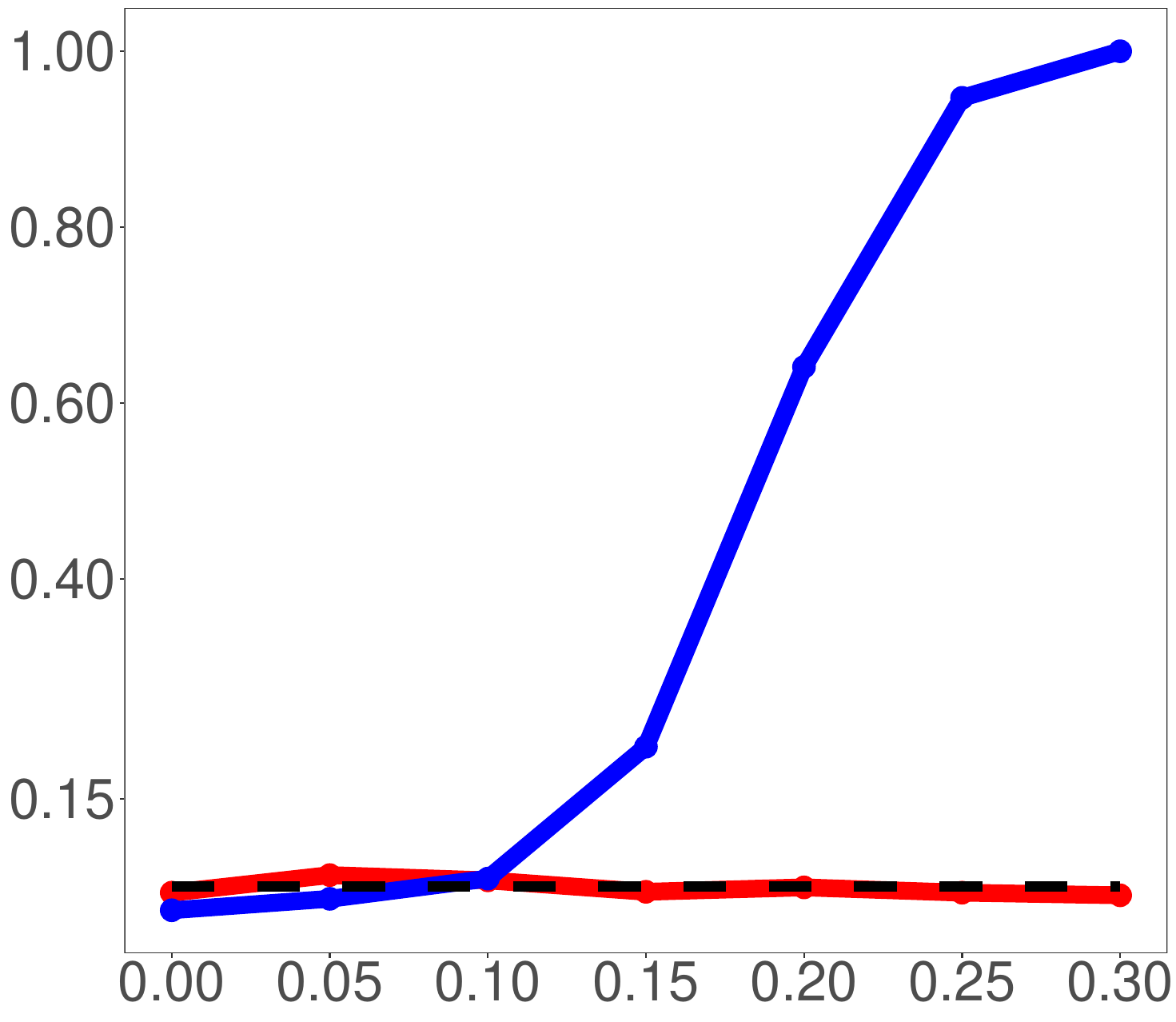}
		\put(-10,-4){\rotatebox{90}{ {  \ \ \ \footnotesize rejection rate  \ \ }}}
\put(50,-9){ \footnotesize $h$  }
\put(27,-22){ Model (4) }
	\end{overpic}}
 % \vspace{0.2cm}
\caption{Results for $y_{11}\sim \frac{\textup{Beta}(2,1.5) - 4/7}{\sqrt{8/147}}$ when $p/n = 1.5$ and $p=600$.}
\label{fig:beta&1.5}
\end{figure}
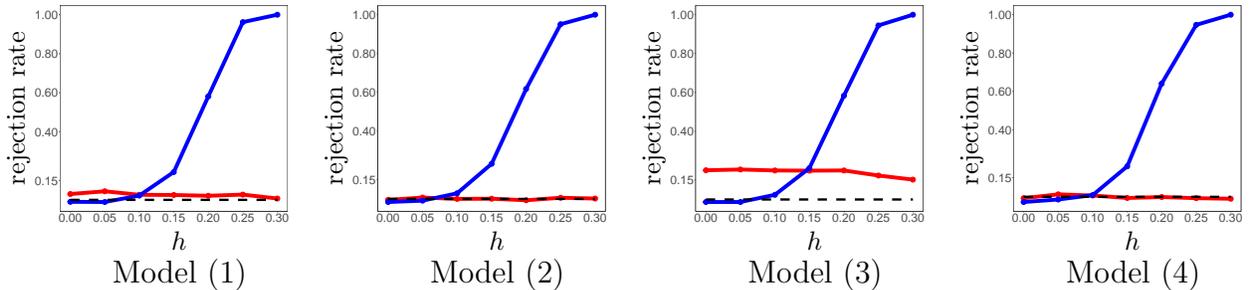

\noindent\textbf{Discussion of results.} The results are displayed in Figures~\ref{fig:lap&0.5}-\ref{fig:beta&1.5}, which contain 24 panels corresponding to the different settings of $(y_{11},\Sigma,n,p)$. In each panel, the rejection rates of the proposed test and the normality test are respectively plotted with blue and red curves. Overall, the most apparent difference between the tests is that they tend to produce curves with different shapes. The power of the proposed test has the intuitive property that it increases toward 1 as $h$ increases, whereas the power of the normality test stays relatively flat as a function of $h$. (Recall that larger values of $h$ correspond to data-generating distributions that are farther from being elliptical.) To fully appreciate the favorable power characteristics of the proposed test, it is worth emphasizing that normal models correspond to a strict subclass of elliptical models. So, in a heuristic sense, normality is an easier hypothesis to reject---and yet, the proposed test still manages to outperform.

In addition to the power of the tests, it is important to compare how well they respect the nominal level of 5\%, which is marked with a dashed black line in every panel. In all settings, the empirical level of the proposed test at $h=0$ is close to the desired value of 5\%, whereas the empirical level of the normality is often substantially inflated, and in several cases it is nearly 20\%. Due to this issue, the fact that the red curves are sometimes higher than the blue curves for small values of $h$ is not a disadvantage of the proposed test. In other words, the higher values of the red curves for small $h$ occur simply because those curves start too high when $h=0$.

\noindent\textbf{Remark on covariates with equal kurtosis.} Apart from the comparisons of power and level just discussed, a distinct point that deserves special emphasis is the power of the proposed test to detect non-elliptical distributions whose covariates have the same kurtosis. Considering that the test is motivated by the identity~\eqref{eqn:kurt}, which implies that all covariates in an elliptical model have the same kurtosis, it might seem surprising that the test can still be powerful against alternatives for which all covariates have the same kurtosis. Nevertheless, the results for Model (4) in Figures~\ref{fig:lap&0.5}-\ref{fig:beta&1.5} show that the test achieves substantially higher power in such cases compared to the normality test. This behavior becomes more understandable once it is recognized that the kurtosis estimate $\breve \kappa$ used to construct the proposed test is designed based on a distinct structural property of elliptical models. Hence, it is less surprising that $\breve \kappa$ disagrees with the other kurtosis estimate $\tilde \kappa$ when data are generated from non-elliptical distributions whose covariates have the same kurtosis.

\begin{figure}[H]
\setlength{\abovecaptionskip}{20pt}
\centering 
\subfigure{
\begin{overpic}[width=0.21\textwidth]{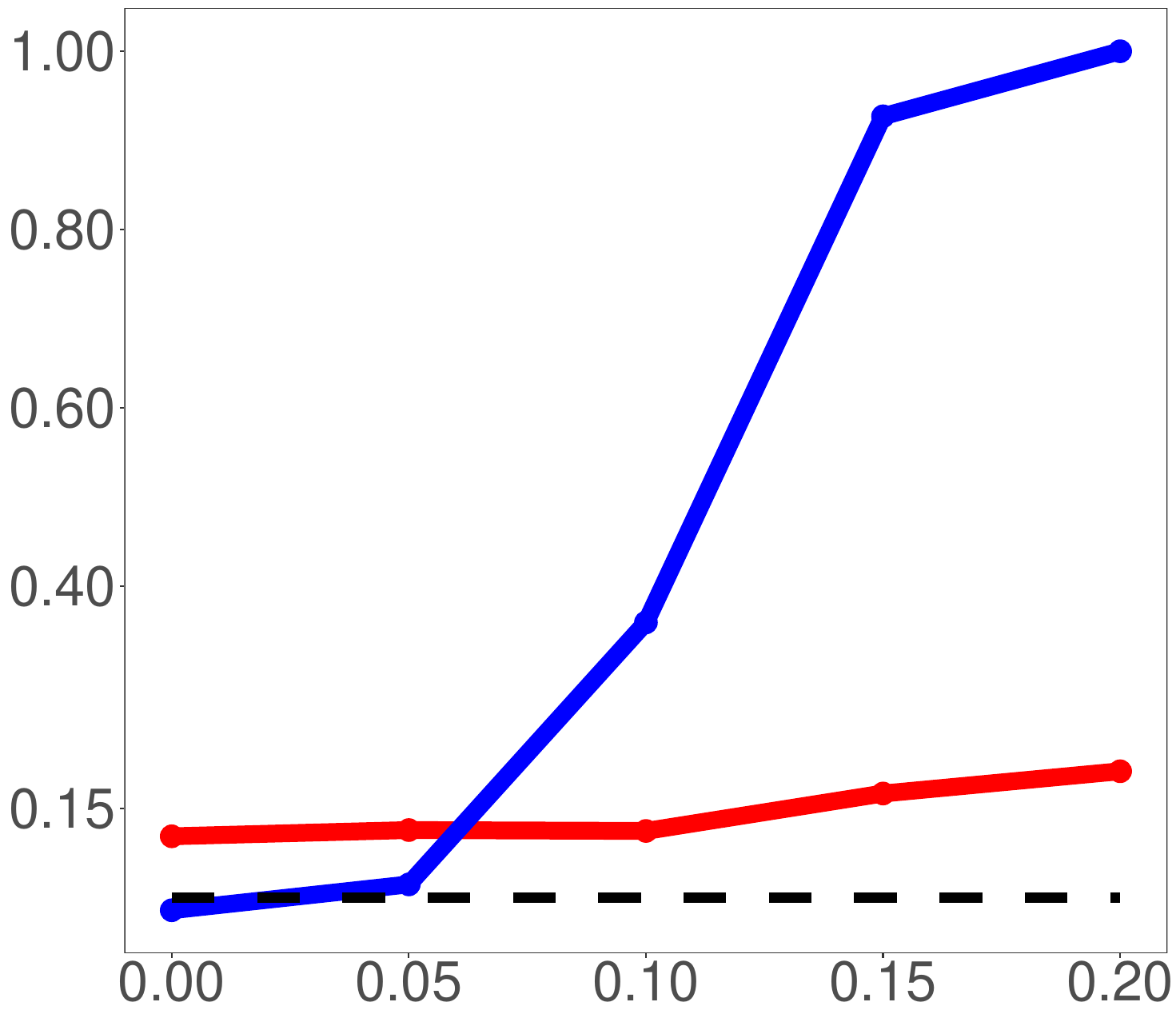}
		\put(-10,-4){\rotatebox{90}{ {  \ \ \ \footnotesize rejection rate  \ \ }}}
\put(50,-9){ \footnotesize $h$  }
\put(27,-22){ Model (1) }
	\end{overpic}}
\hspace{4mm}
\subfigure{
\begin{overpic}[width=0.21\textwidth]{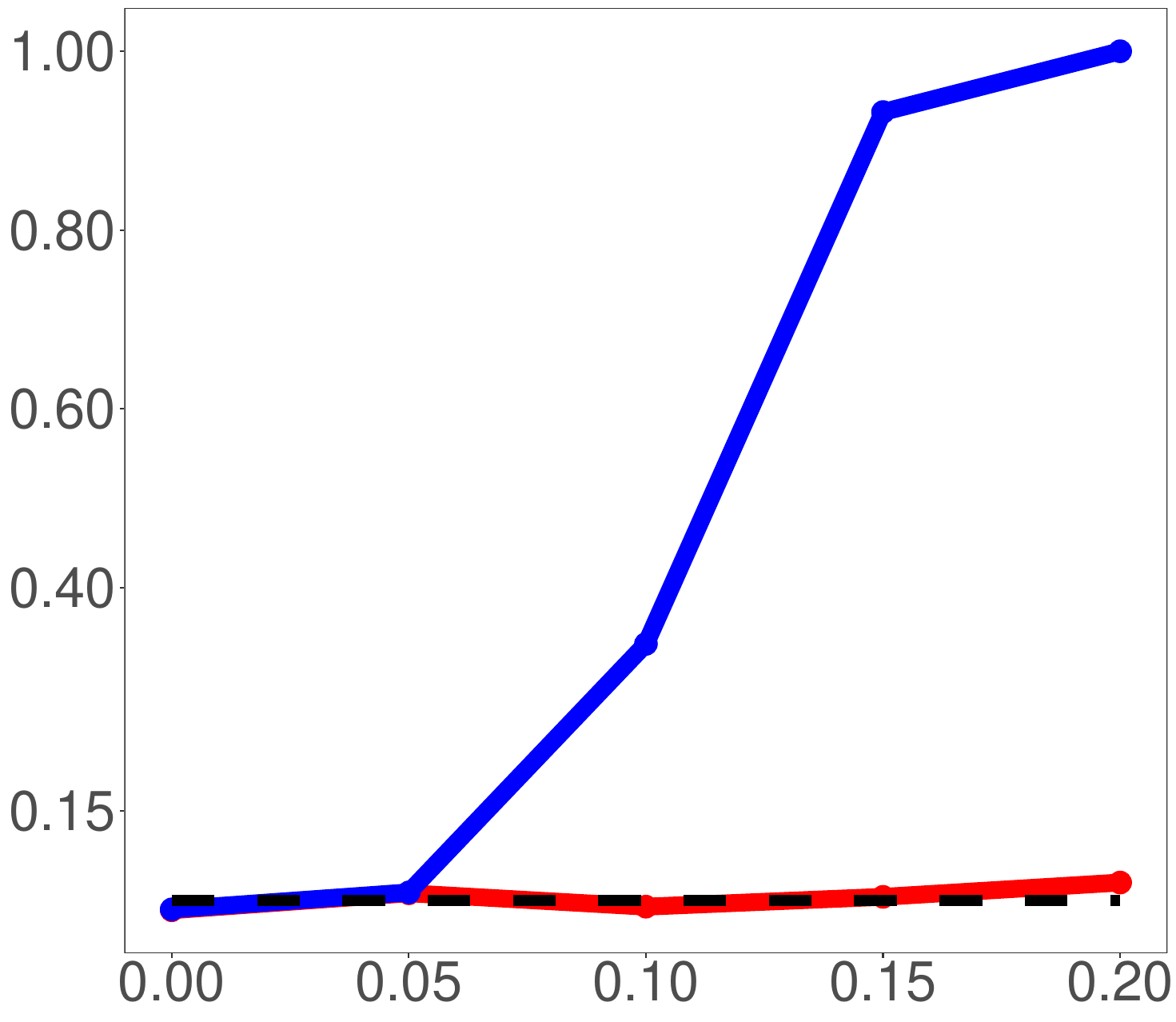}
		\put(-10,-4){\rotatebox{90}{ {  \ \ \ \footnotesize rejection rate  \ \ }}}
\put(50,-9){ \footnotesize $h$  }
\put(27,-22){ Model (2) }
	\end{overpic}}
\hspace{4mm}
\subfigure{
\begin{overpic}[width=0.21\textwidth]{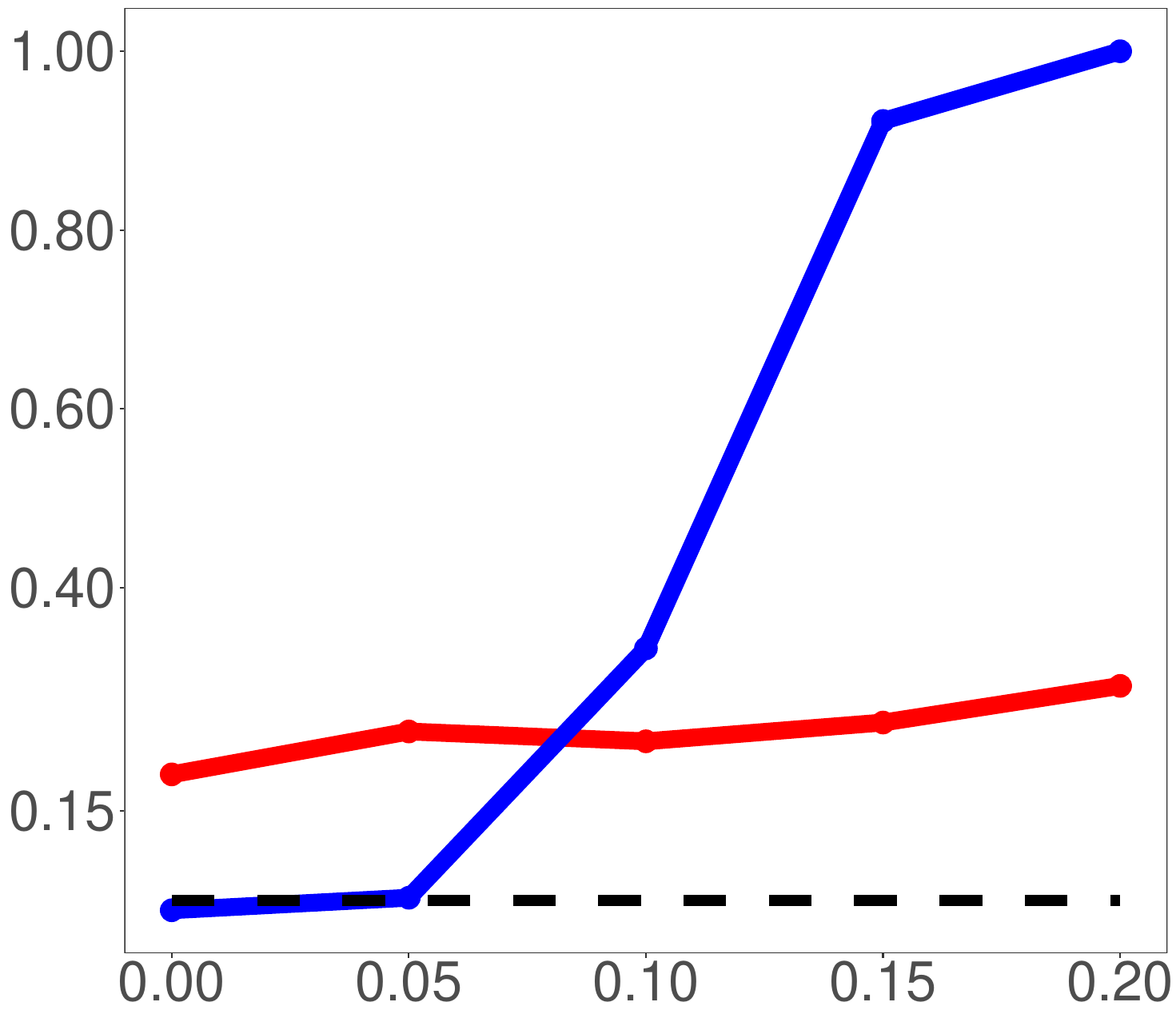}
		\put(-10,-4){\rotatebox{90}{ {  \ \ \ \footnotesize rejection rate  \ \ }}}
\put(50,-9){ \footnotesize $h$  }
\put(27,-22){ Model (3) }
	\end{overpic}}
 \hspace{4mm}
\subfigure{
\begin{overpic}[width=0.21\textwidth]{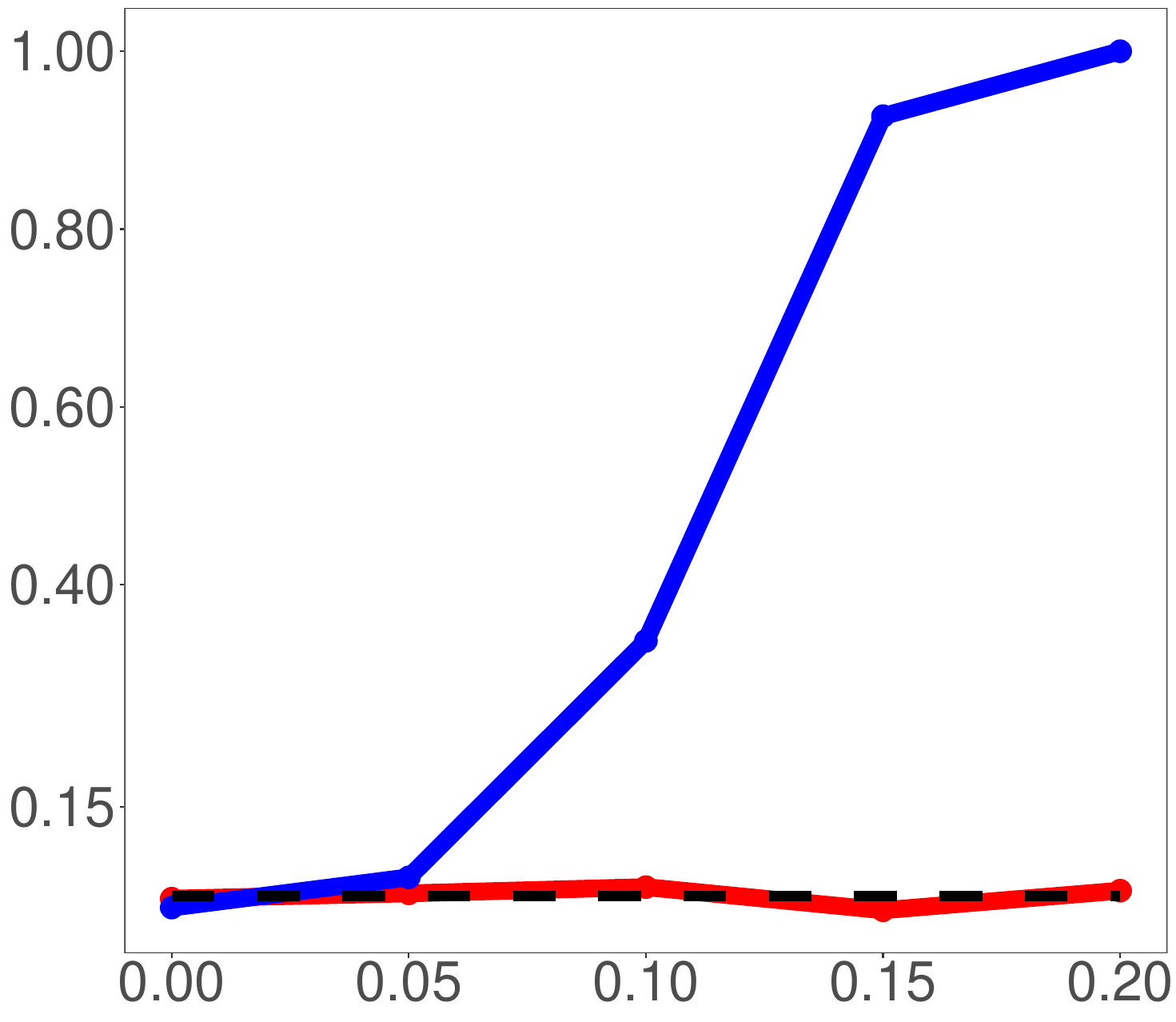}
		\put(-10,-4){\rotatebox{90}{ {  \ \ \ \footnotesize rejection rate  \ \ }}}
\put(50,-9){ \footnotesize $h$  }
\put(27,-22){ Model (4) }
	\end{overpic}}
 % \vspace{0.2cm}
\caption{Results for $y_{11}\sim \frac{\textup{Laplace}(0,1)}{\sqrt{2}}$ when $p/n = 1$ and $p = 400$.}
\end{figure}

\begin{figure}[H]
\setlength{\abovecaptionskip}{20pt}
\centering 
\subfigure{
\begin{overpic}[width=0.21\textwidth]{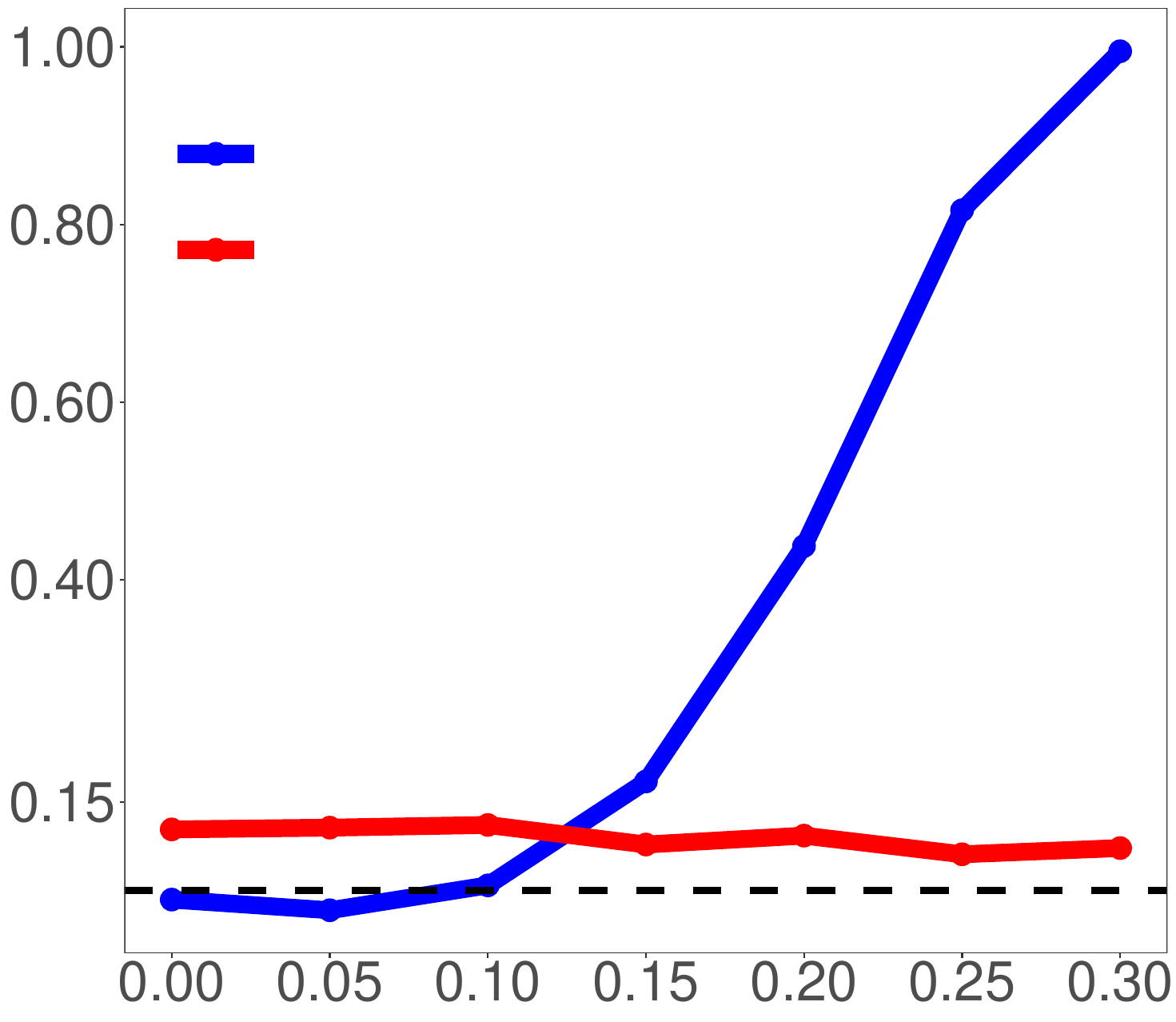}
		\put(-10,-4){\rotatebox{90}{ {  \ \ \ \footnotesize rejection rate  \ \ }}}
\put(50,-9){ \footnotesize $h$  }
\put(27,-22){ Model (1) }
\put(20,63){ \scriptsize normality test}
\put(20,73){ \scriptsize proposed test}
	\end{overpic}}
\hspace{4mm}
\subfigure{
\begin{overpic}[width=0.21\textwidth]{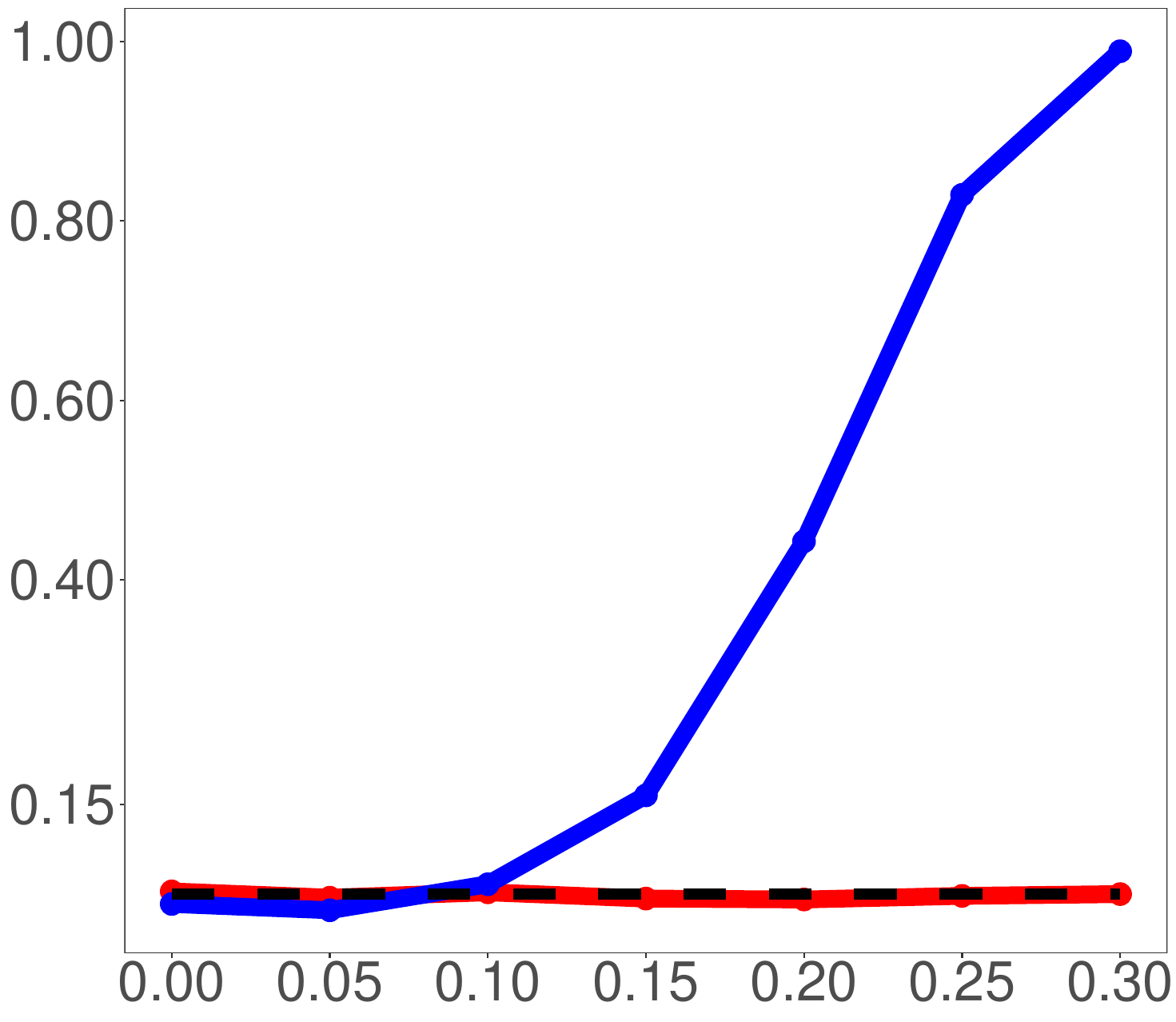}
		\put(-10,-4){\rotatebox{90}{ {  \ \ \ \footnotesize rejection rate  \ \ }}}
\put(50,-9){ \footnotesize $h$  }
\put(27,-22){ Model (2) }
	\end{overpic}}
\hspace{4mm}
\subfigure{
\begin{overpic}[width=0.21\textwidth]{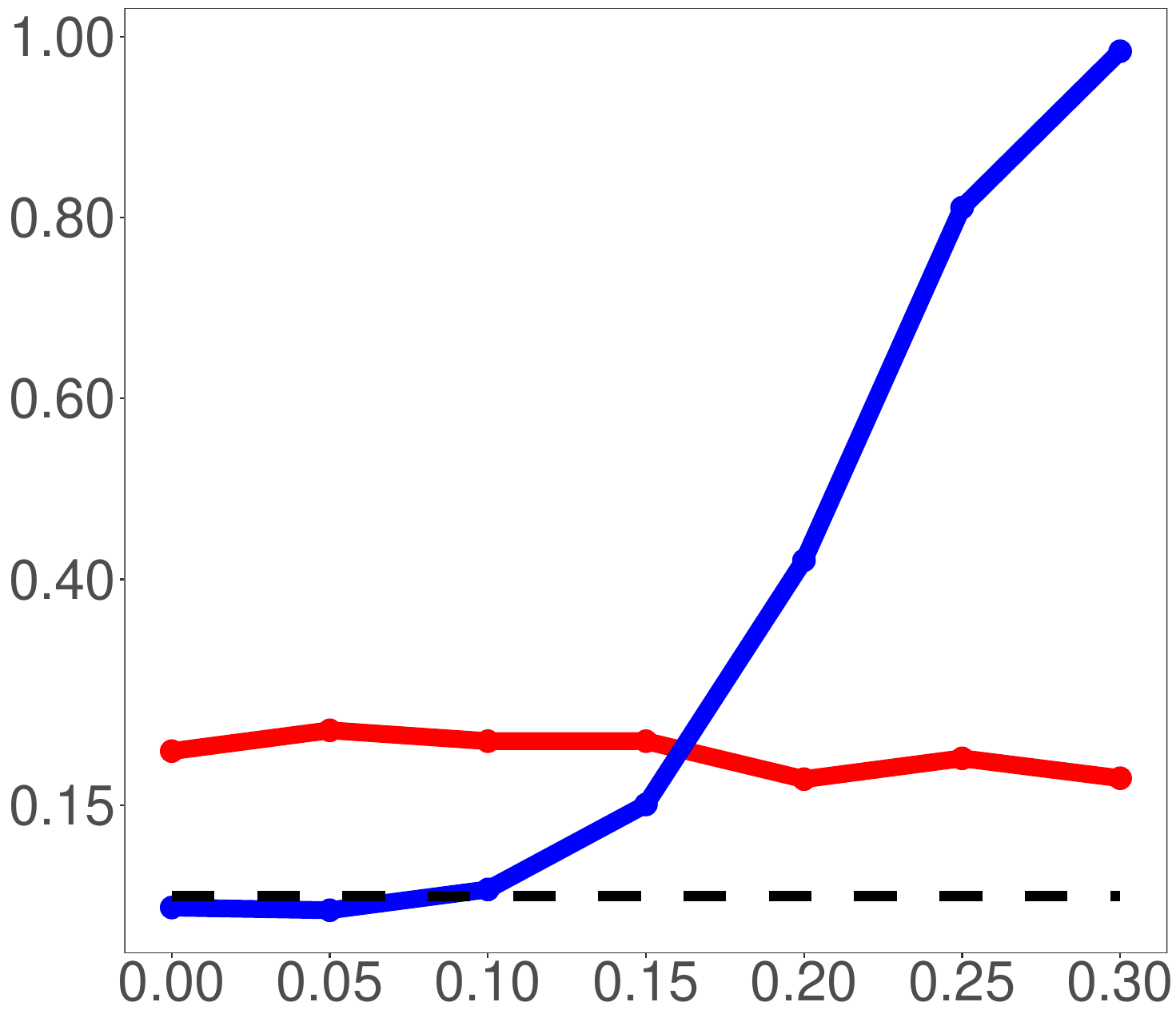}
		\put(-10,-4){\rotatebox{90}{ {  \ \ \ \footnotesize rejection rate  \ \ }}}
\put(50,-9){ \footnotesize $h$  }
\put(27,-22){ Model (3) }
	\end{overpic}}
\hspace{4mm}
\subfigure{
\begin{overpic}[width=0.21\textwidth]{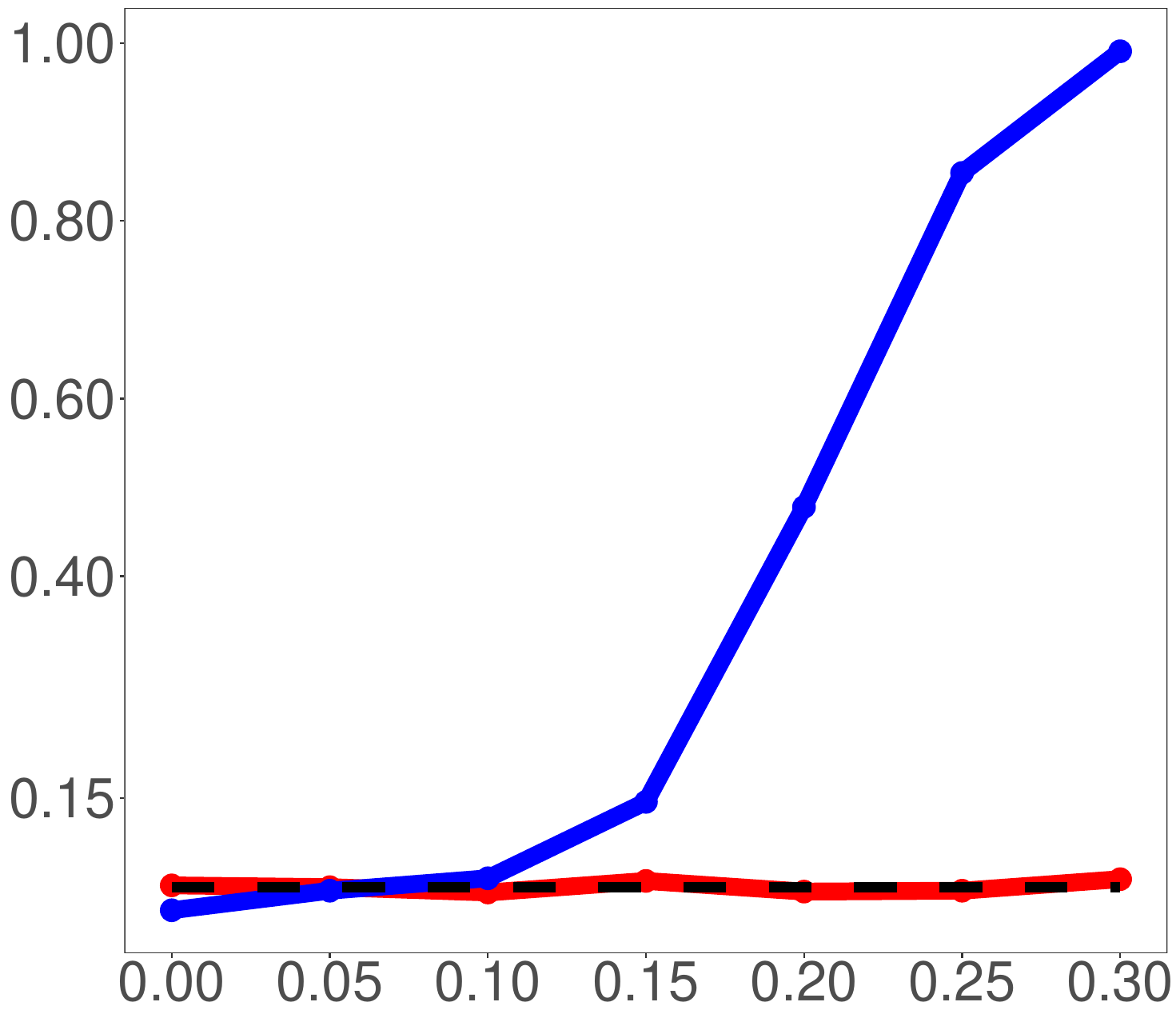}
		\put(-10,-4){\rotatebox{90}{ {  \ \ \ \footnotesize rejection rate  \ \ }}}
\put(50,-9){ \footnotesize $h$  }
\put(27,-22){ Model (4) }
	\end{overpic}}
 % \vspace{0.2cm}
\caption{Results for $y_{11}\sim \frac{\textup{Beta}(2,1.5) - 4/7}{\sqrt{8/147}}$ when $p/n = 1$ and $p=400$.}
\end{figure}

\begin{figure}[H]
\setlength{\abovecaptionskip}{20pt}
\centering 
\subfigure{
\begin{overpic}[width=0.21\textwidth]{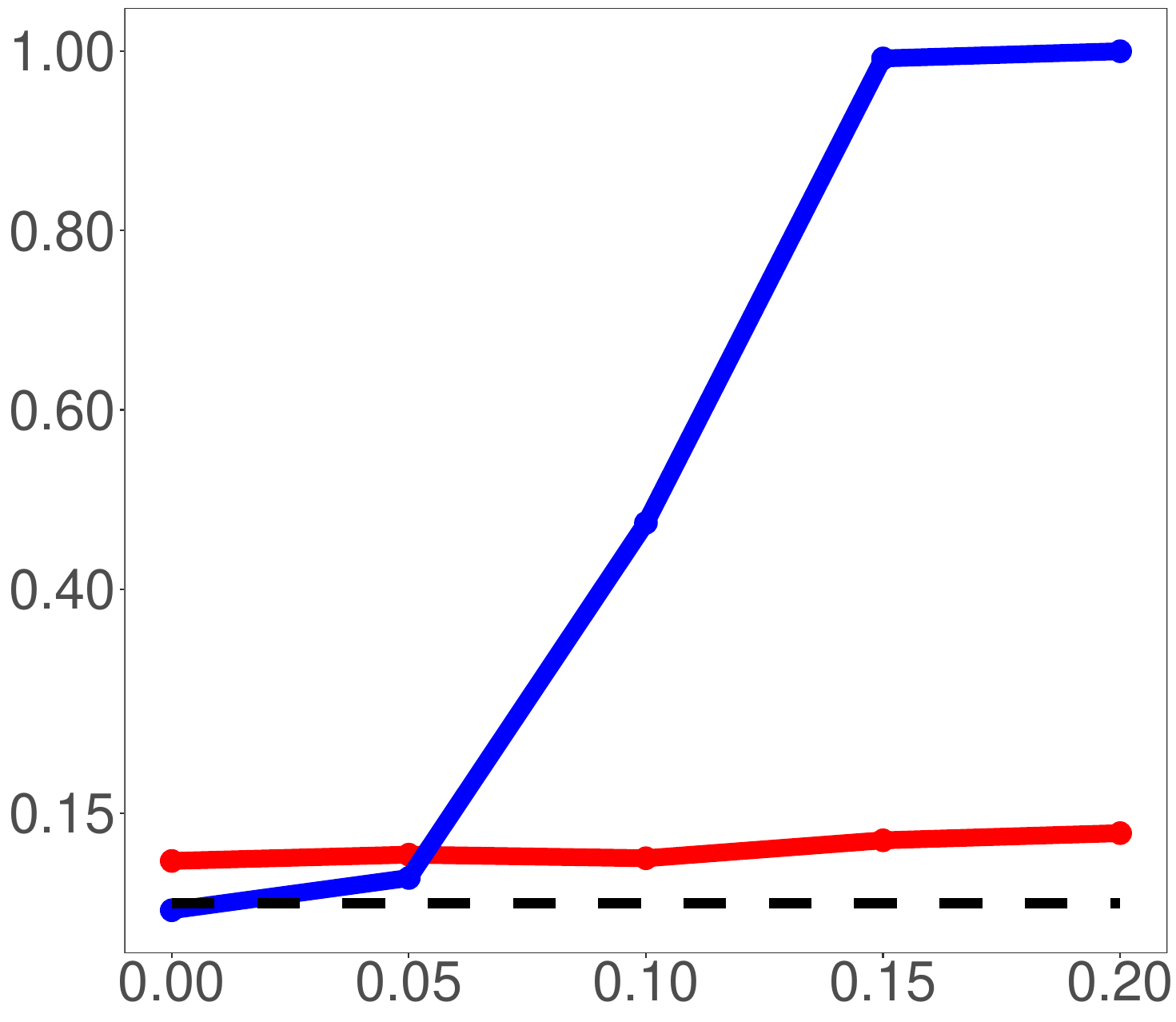}
		\put(-10,-4){\rotatebox{90}{ {  \ \ \ \footnotesize rejection rate  \ \ }}}
\put(50,-9){ \footnotesize $h$  }
\put(27,-22){ Model (1) }
	\end{overpic}}
\hspace{4mm}
\subfigure{
\begin{overpic}[width=0.21\textwidth]{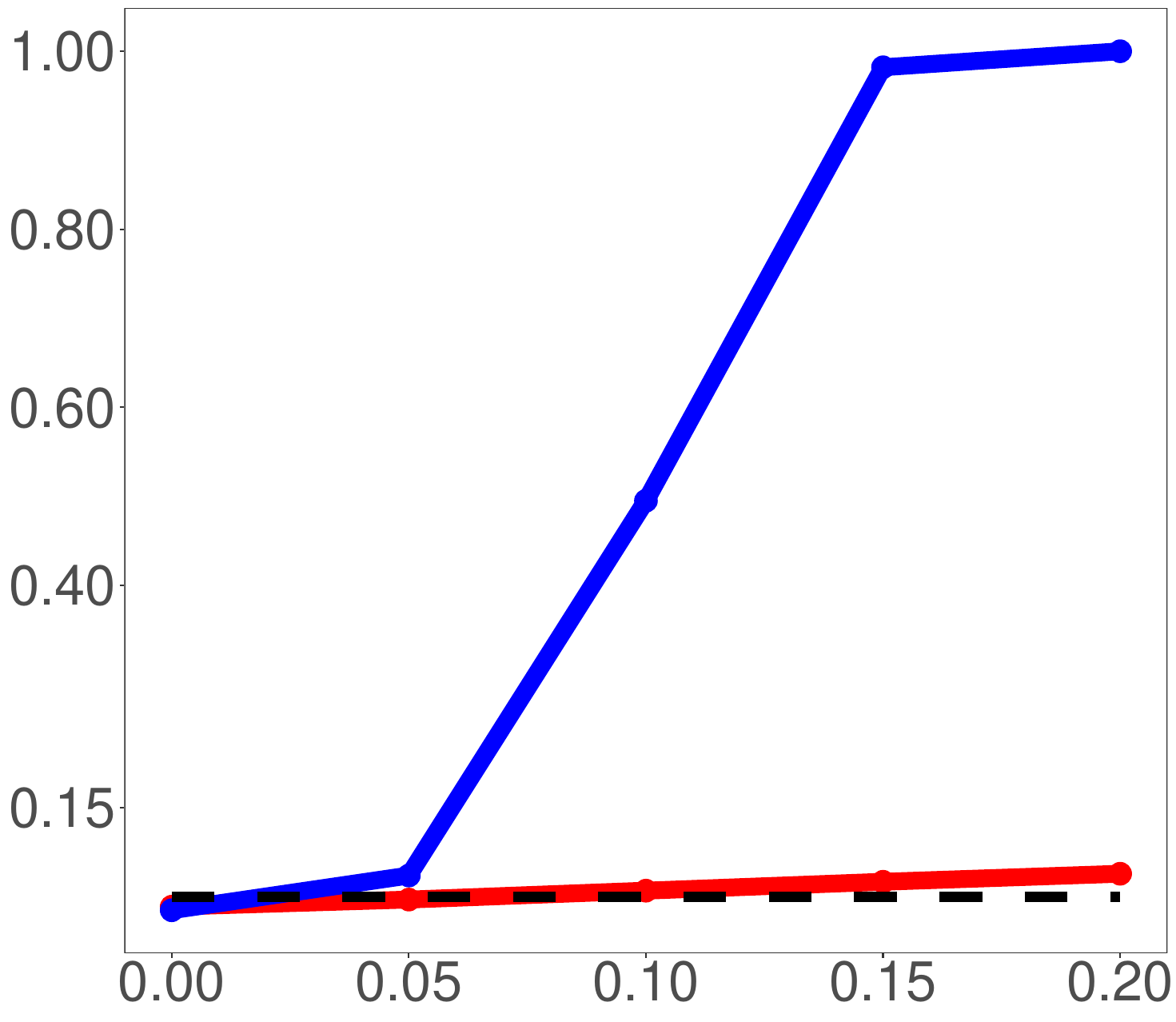}
		\put(-10,-4){\rotatebox{90}{ {  \ \ \ \footnotesize rejection rate  \ \ }}}
\put(50,-9){ \footnotesize $h$  }
\put(27,-22){ Model (2) }
	\end{overpic}}
\hspace{4mm}
\subfigure{
\begin{overpic}[width=0.21\textwidth]{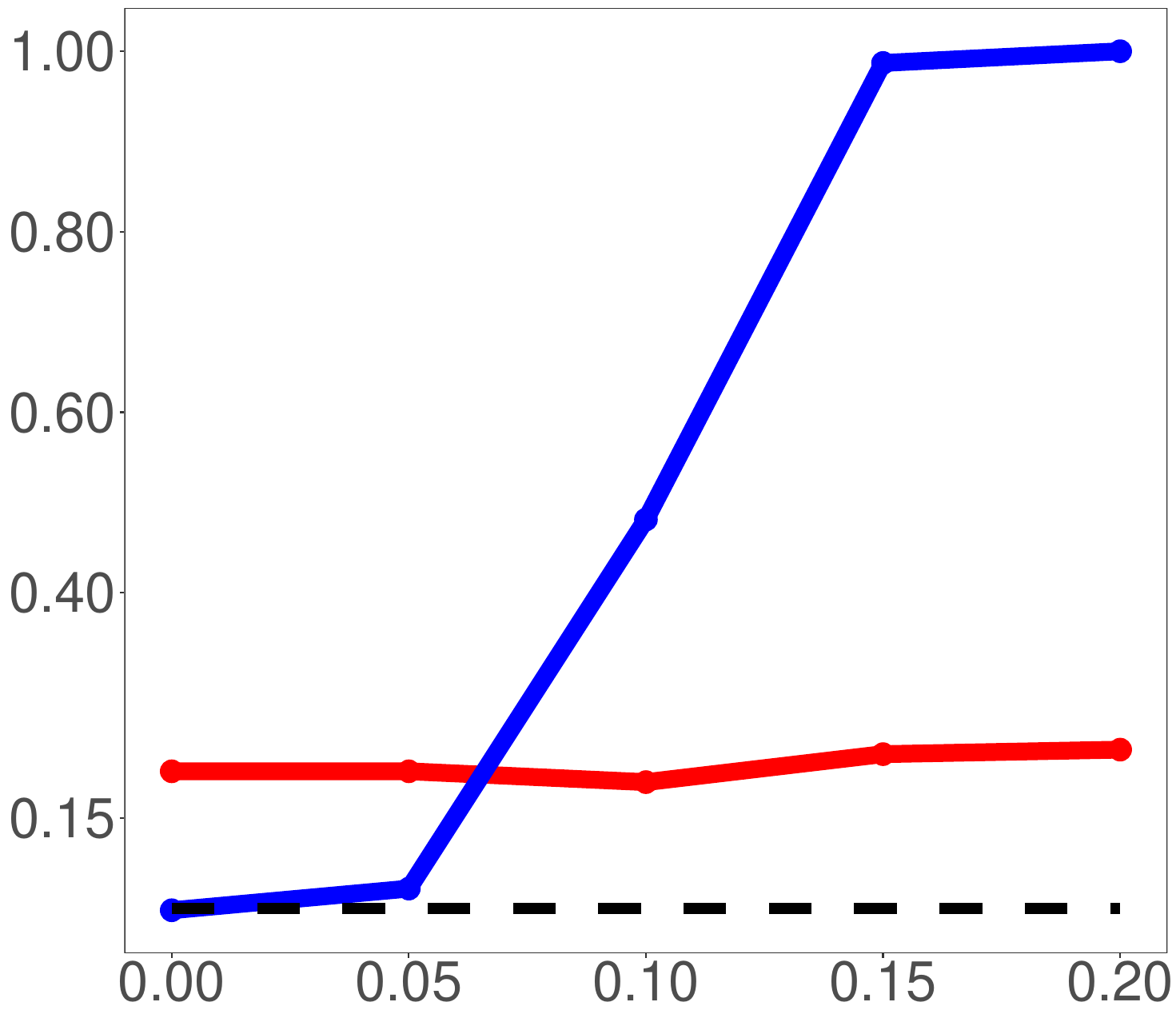}
		\put(-10,-4){\rotatebox{90}{ {  \ \ \ \footnotesize rejection rate  \ \ }}}
\put(50,-9){ \footnotesize $h$  }
\put(27,-22){ Model (3) }
	\end{overpic}}
\hspace{4mm}
\subfigure{
\begin{overpic}[width=0.21\textwidth]{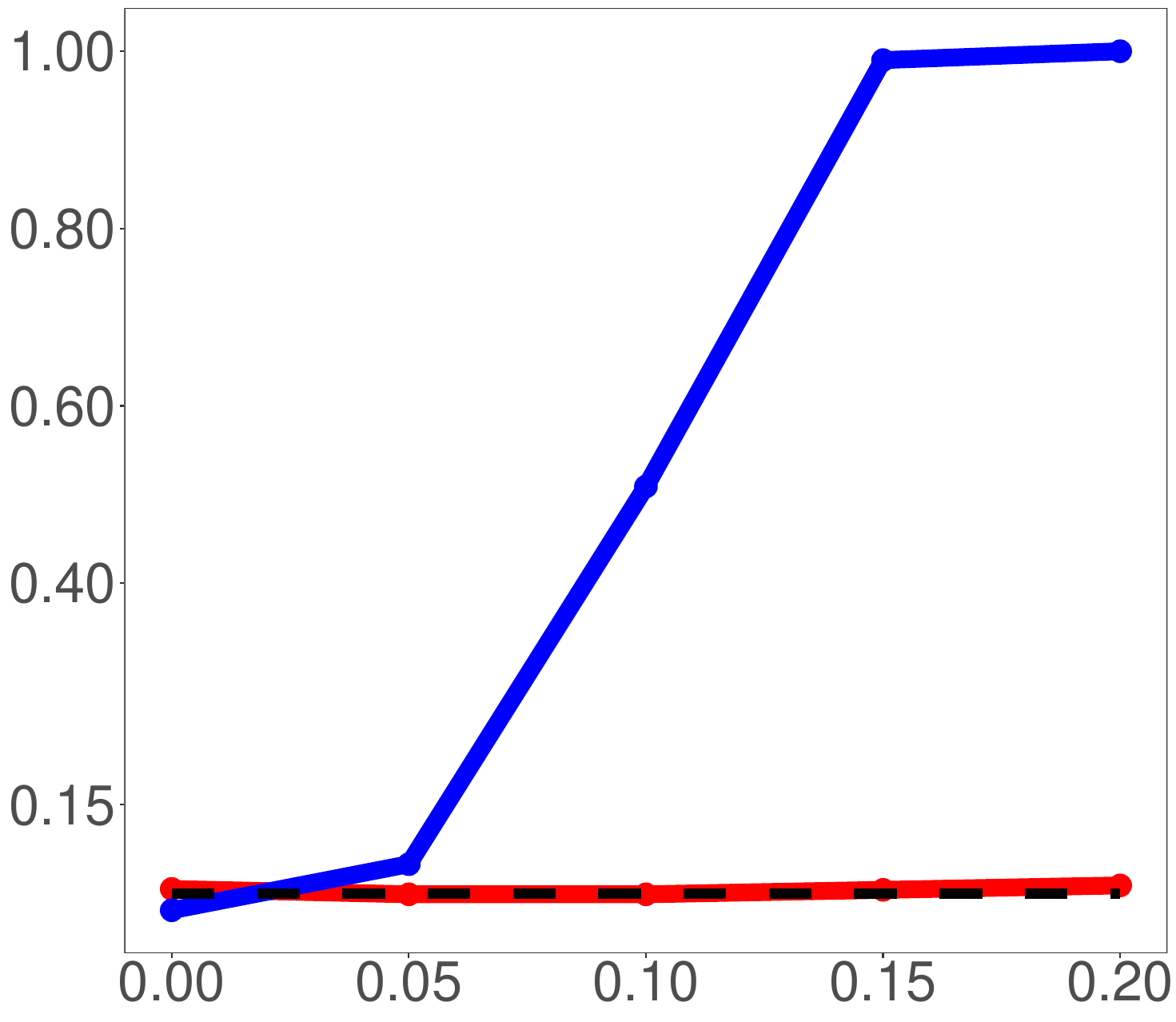}
		\put(-10,-4){\rotatebox{90}{ {  \ \ \ \footnotesize rejection rate  \ \ }}}
\put(50,-9){ \footnotesize $h$  }
\put(27,-22){ Model (4) }
	\end{overpic}}
 % \vspace{0.2cm}
\caption{Results for $y_{11}\sim \frac{\textup{Laplace}(0,1)}{\sqrt{2}}$ when $p/n = 1.5$ and $p=600$.}
\end{figure}

\subsection{Stock market data}
Elliptical distributions are commonly used in finance to model investment returns~\citep{Embrechts:2011,gupta2013elliptically}. Here, we apply the proposed test as well as the normality test from~\citep{chen2023normality} to a collection of stock market return datasets with dimensions ranging from 100 to 480. The datasets were constructed from the monthly log returns of 480 stocks in the S\&P 500 that had complete histories between July 2012 and June 2022, based on the records available at https://finance.yahoo.com/lookup. Hence, the initial ``full'' dataset may be regarded as a matrix of size $120\times 480$, based on $n=120$ monthly observations in $p=480$ dimensions. Next, we labeled the stocks according to the ranks of their prices at the beginning of July 2012, with the first stock having the highest price at that time. Lastly, for each number $d\in\{100,110,120,\dots,480\}$, we applied both the proposed test and the normality test to the $n\times d$ data matrix corresponding to the first $d$ columns of the full $n\times p$ matrix.

\begin{figure}[H]
\centering
	\DeclareGraphicsExtensions{.eps}
	\begin{overpic}[width=0.4\textwidth]{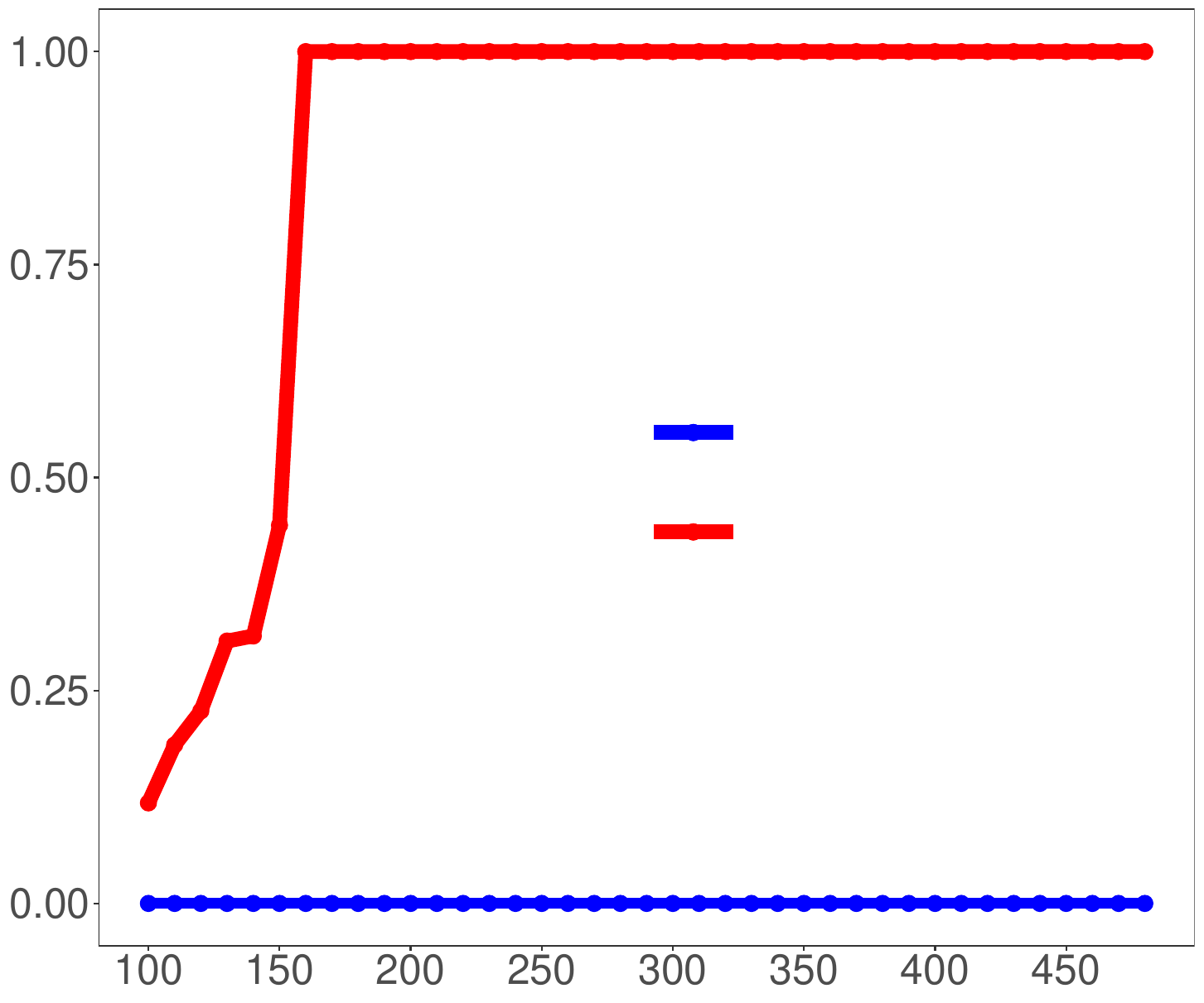}
		\put(-9,24){\rotatebox{90}{ {  \ \ \ p-value  \ \ }}}
\put(35,-7){  dimension $d$ }
\put(61,38){ \footnotesize normality test}
\put(61,47){ \footnotesize proposed test}
	\end{overpic}
 	\caption{p-values for the stock market data using the proposed test and normality test.} 
	\label{fig:stock}
\end{figure}
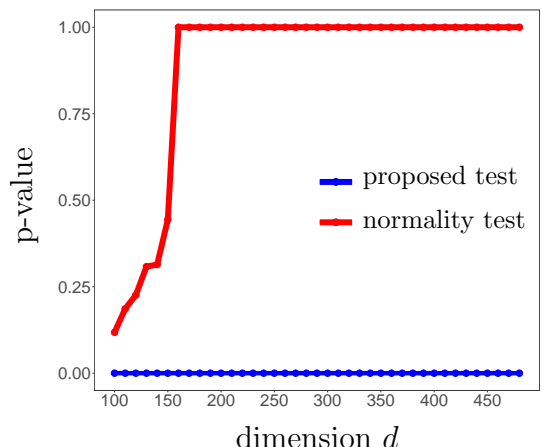

In Figure~\ref{fig:stock}, we plot the resulting p-values as a function of $d$, with the proposed test shown in blue, and the normality test shown in red. For all choices of $d$, the proposed test produced p-values that are effectively 0. On the other hand, none of the p-values of the normality test fall below 0.05, and its p-values corresponding to $d\geq 160$ are all equal to 1. Thus, the results for this example conform with the power advantage of the proposed test that was demonstrated earlier using simulated data.
\vspace{0.2cm}

\subsection{Breast cancer data}

For our second example, we look at a breast cancer dataset that was studied in the influential paper~\cite{van2002gene1} and is available in the R package %{\tt{seventyGeneData}} \citep{seventyGeneData} and 
{\tt{cancerdata}}~\citep{cancerdata}. The observations correspond to $n=295$ tumor samples with $p=24481$ gene expression measures. We excluded the gene expression measures with missing values, resulting in $p = 22223$ gene expression measures. After centering these observations with the sample mean, we constructed 300 datasets in the following way. For each number $d\in\{200,300,500\}$, we randomly sampled $d$ columns without replacement from the full $n\times p$ matrix and repeated this 100 times. So, in other words, for each choice of $d$, we generated 100 datasets of size $n\times d$ corresponding to random subsets of genes. 

We applied the proposed test and the normality test to all the randomly sampled datasets and recorded the p-values. Then, for each choice of $d$, we computed the median of the 100 associated p-values for each test. The proposed test produced median p-values that were effectively 0 for each choice of $d$, whereas the normality test produced median p-values equal to 0.656, 1, and 1 in the respective cases of $d=200, 300, 500$. Once again, the results of this example reflect the power advantage of the proposed test that was observed in the experiments with simulated data.

\bibliographystyle{apalike}

\bibliography{citation}

\newpage

\begin{center}
{\Large\bf Supplementary Material:\\[0.2cm] Testing Elliptical Models in High Dimensions}

Siyao Wang and Miles E. Lopes
\end{center}
\appendix
We first introduce some notation and frequently used facts in Appendix~\ref{app:notation}. In Appendix~\ref{app:T}, we show the asymptotic normality of $T_n$ by analyzing $\tilde \kappa$ and $\check \kappa$ in Appendices~\ref{app:T1} and~\ref{app:T2} respectively. The ratio-consistency of $\hat \sigma_n^2$ is shown in Appendix~\ref{app:estimators}. Theoretical background results are given in Appendix~\ref{app:background}. Additional empirical results are given in Appendix~\ref{app:dim}.

\section{Preliminaries} \label{app:notation}
For two sequences of random variables $\{y_n\}$ and $\{v_n\}$, the relation $y_n=o_{\P}(v_n)$ means that $y_n/v_n \xrightarrow{\P} 0$, and the relation $y_n=\mathcal{O}_{\P}(v_n)$ means that for all $\epsilon > 0$, there is a constant $M$  not depending on $n$ such that $\sup_{n\geq 1} \P(| y_n/v_n| > M) < \epsilon$. The $k\times k$ identity matrix is denoted as $I_k$, and the indicator function for a condition $\cdots$ is denoted as $1\{\cdots\}$. The $j$th canonical basis vector in an Euclidean space is denoted as $\mathbf{e}_{j}$, with the dimension being understood from context. For two matrices $A$ and $B$ of the same size, their Hadamard product $A\circ B$ is the matrix of that size whose $ij$ entry is given by $(A\circ B)_{ij} = A_{ij} B_{ij}$, and $A^{\circ k}$ denotes the $k$-fold Hadamard product of $A$ with itself. The stable rank of the covariance matrix $\SIGMA$ is defined as  ${\tt{r}}(\SIGMA) = \frac{\tr(\SIGMA)^2}{\tr(\SIGMA^2)}$. The correlation matrix associated to $\Sigma$ is denoted by $R$.

Next, let $\z_1, \ldots,\z_n$, denote i.i.d.~standard normal vectors in $\R^p$ independent of $\xi_1,\dots,\xi_n$, and for each $i=1,\dots,n$ let 
\begin{align}
\label{eq:definitionofzeta}
    \boldsymbol\zeta_i  \ = \  \SIGMA^{1/2}\z_i.
\end{align}
This definition allows us to express the $i$th observation $\x_i$ as
\begin{align}
\label{eq:xi}
    \x_i  \ = \  \frac{\xi_i}{\|\z_i\|_2}\boldsymbol\zeta_i  \ = \  \frac{\xi_i}{\|\z_i\|_2} \SIGMA^{1/2}\z_i
\end{align}
for $i = 1, \ldots, n$.
For any $k\geq 0$, we define the parameter
\begin{equation}\label{eq:defofrk}
    \begin{split}
        r_k \ = \ &\frac{\E(\|\mathbf{x}_{1}\|_2^{2k})}{\E((\z_1\ttop \SIGMA \z_1 )^{k})},
    \end{split}
\end{equation}
which can also be expressed as
\begin{equation}\label{eqn:extrark}
          r_k \ = \  \frac{\E(\xi_1^{2k})}{\E(\|\z_1\|_2^{2k})}.
\end{equation}
In particular, we will often use the fact that Assumption \ref{Data generating model} implies
\begin{equation}\label{eq:relationofr2}
    r_2  \ = \  1 + \frac{\tau-2}{p} +o\big(\ts \frac{1}{p}\big).
\end{equation}
By noting that $\boldsymbol\zeta_1 \sim N(0, \SIGMA)$ and $\zeta_{1j} \sim N(0,\Sigma_{jj})$, we have the following moment formula for any integer $k\geq 0$ and $j=1,\dots,p$,
\begin{align*}
    \E(x_{1j}^{2k})  \ = \  &\E\Big(\ts\frac{\xi_1^{2k}}{\|\z_1\|_2^{2k}}  \zeta_{1j}^{2k}\Big)  \\[0.2cm]
\ = \  &\frac{\E(\xi_1^{2k})}{\E(\|\z_1\|_2^{2k})}\E(\zeta_{1j}^{2k}) \\[0.2cm]
\ = \   & (2k-1)!! r_k \Sigma_{jj}^{k},
\end{align*}
where $m!!=m(m-2)(m-4)\cdots 1$ for an odd integer $m$.
In the special case when the diagonal entries of $\SIGMA$ are equal to 1, we will often use the fact that
\begin{equation}\label{eq:momentboundofxij}
   \max_{1\leq j\leq p} \E(x_{1j}^{2k}) \ \lesssim \ 1
\end{equation}
holds for $k=1,\ldots,8$, which is established in Lemma \ref{lem:moment}.

\setcounter{lemma}{0}
\renewcommand{\thelemma}{B.\arabic{lemma}}
\section{Asymptotic normality of $T_n$} \label{app:T}
The following proposition is a counterpart to Theorem~\ref{thm:main} in which $\hat\sigma_n$ is replaced by $\sigma_n$.
\begin{proposition}
    \label{prop:T}
    If Assumption \ref{Data generating model} holds, then as $n\to\infty$,
    \begin{equation}
    \label{eq:T}
    \frac{T_n}{\sigma_n} \ \Rightarrow \ N(0,1).
    \end{equation}
    \end{proposition}
\proof We decompose the statistic $T_n$ as
\begin{equation*}
    \begin{split}
        T_n \ = \ &\sqrt{\frac{pn}{2}}\Big(\frac{\tilde\kappa - \kappa}{3}  + \frac{4}{n} \Big) - \sqrt{\frac{pn}{2}}\Big( \frac{\check{\kappa}- \kappa}{3} \Big) \\
        \ =: \ & T_{n,1} - T_{n,2}.
    \end{split}
\end{equation*}
Also, the formula for $\kappa$ given in \eqref{eq:definitionofkappa1} allows us to write
\begin{align}
    T_{n,1} & \ = \ \sqrt{\frac{pn}{2}}\Big(\frac{\TILDEKAPPA}{3}-1\Big) - \sqrt{\frac{pn}{2}}\Big(\frac{ \kappa}{3}-1\Big) +2\sqrt{\frac{2p}{n}}\label{eqn:T1equiv}\\[0.2cm]
    T_{n,2} & \ = \ \sqrt{\frac{pn}{2}} \Big( \frac{ \check{\varsigma}^2 - 2(\check{\nu}_2 - \frac{2}{n} \check{\nu}_1^2)}{\check{\nu}_1^2 + 2(\check{\nu}_2 - \frac{2}{n} \check{\nu}_1^2)} - \frac{  {\varsigma}^2 - 2 {\nu}_2  }{ {\nu}_1^2 + 2 {\nu}_2  }\Big).\label{eq:defofT2}
\end{align}
We analyze the asymptotic properties of $T_{n,1}$ and $T_{n,2}$ separately. Proposition \ref{prop:T1} shows that as $n\to\infty$,
\begin{equation}\label{eqn:T1normal}
     \frac{T_{n,1}}{\sigma_{n,1}}  \ \Rightarrow \ N(0,1).
\end{equation} 

Regarding $T_{n,2}$,  Proposition \ref{prop:T2smallr} shows that if ${\tt{r}}(\SIGMA) \lesssim \sqrt{p}$, then as $n \to \infty$,
    \begin{align*}
        \frac{T_{n,2}}{\sigma_{n,2}}  \ \Rightarrow \ N(0,1).
    \end{align*} 
Since $T_{n,1}$ and $T_{n,2}$ are constructed from independent sets of observations, it follows that $\big(\frac{1}{\sigma_{n,1}} T_{n,1}, \frac{1}{\sigma_{n,2}} T_{n,2} \big)$ converges weakly to a standard bivariate normal distribution, which implies the stated result~\eqref{eq:T} when ${\tt{r}}(\SIGMA) \lesssim \sqrt{p}$.

When  ${\tt{r}}(\SIGMA)/\sqrt p\to\infty$, Proposition~\ref{prop:T2larger} shows that
    \begin{align}\label{eqn:T2lim0}
        \frac{T_{n,2}}{\sigma_n} \ \xrightarrow{\P} \ 0.
    \end{align}
Also, when ${\tt{r}}(\SIGMA)/\sqrt p\to\infty$, it can be checked that $\frac{\sigma_{n,1}^2}{\sigma_n^2} \to 1$. Therefore, the limits~\eqref{eqn:T1normal} and~\eqref{eqn:T2lim0} imply~\eqref{eq:T} because
    \begin{align*}
    \frac{T_{n,1} -T_{n,2}}{\sigma_n} & \ = \ \frac{T_{n,1}}{\sigma_{n,1}} + o_{\P}(1) \ \Rightarrow \ N(0,1). 
    \end{align*}
\qed

\subsection{Asymptotic normality of $T_{n,1}$}

\label{app:T1}
\begin{proposition}
\label{prop:T1}
 If Assumption \ref{Data generating model} holds, then as $n\to\infty$,
    \begin{equation*}
    \frac{T_{n,1}}{\sigma_{n,1}}\ \Rightarrow \ N(0,1).
    \end{equation*}
\end{proposition}
\proof
Recall from equation~\eqref{eqn:T1equiv} that $T_{n,1}$ may be expressed as
\begin{align*}
    T_{n,1} & \ = \  \sqrt{\frac{pn}{2}}\Big( \frac{\TILDEKAPPA}{3}-1\Big) - \sqrt{\frac{pn}{2}}\Big(\frac{ \kappa}{3}-1\Big) +2\sqrt{\frac{2p}{n}}.
\end{align*}
We begin by analyzing the first term on the right side.
Since each $\tilde \kappa_j$ in the definition of $\TILDEKAPPA$ is invariant to rescaling the coordinates of each observation, we may assume that all the diagonal entries of $\SIGMA$ are equal to 1. Let $\tilde \SIGMA$ denote the sample covariance matrix of the first $n/2$ observations so that
\begin{equation*}
    \tilde \Sigma_{jj}  \ = \  \frac{1}{n/2}\sum_{i=1}^{n/2} x_{ij}^2.
\end{equation*}
Also, for each $j=1,\dots,p$, define the event
\begin{equation}\label{eqn:Ajdef}
A_j \ = \ \big\{ \big|\tilde \Sigma_{jj}  - 1 \big|\leq \e_n \big\},
\end{equation}
where we put $\e_n = n^{-1/20}$.
It is shown in Lemma \ref{lem:A_j} that 
\begin{equation}\label{eqn:Acneg}
    \sqrt{\frac{pn}{2}}\Big(\frac{\TILDEKAPPA}{3}-1\Big) \ = \ \sqrt{\frac{n}{2p}}\sum_{j=1}^p\Big(\frac{\tilde \kappa_j}{3}-1\Big)1\{A_j\} \ + \ o_{\P}(1).
\end{equation}
Next, define the random variables
\begin{align}
\label{eq:definitionofLj}
    \ell_j &  \ = \  \bigg(\frac{1/3}{n/2} \sum_{i=1}^{n/2}\Big( x_{ij}^4 -  \E(x_{1j}^4)\Big)\bigg)
- 2 (\tilde \Sigma_{jj} - 1) \\[0.3cm]
    \Delta_j &  \ = \  -(\tilde \Sigma_{jj} - 1) ^2+\frac{\E(x_{ij}^4)}{3}-1
    \label{eq:definitionofdeltaj}
    \\[0.3cm]
    \ve_j & \ = \ \Big(\frac{1\{A_j\}}{\tilde\Sigma_{jj}^2}-1\Big)1\{\tilde\Sigma_{jj}\neq 0\},\label{eqn:vejef}
\end{align}
which satisfy the algebraic relation
\begin{equation*}
% \label{eqn:veabove}
\Big(\frac{\tilde \kappa_j}{3}-1\Big)1\{A_j\}  \ = \ (1+\ve_j)\big(\ell_j+\Delta_j\big).
\end{equation*}
 It follows from Lemma \ref{lem:epilson} that
\begin{equation*}
    \sqrt{\frac{n}{2p\sigma_{n,1}^2}}\sum_{j=1}^p \ve_j(\ell_j+\Delta_j)  \ = \  o_{\P}(1),
\end{equation*}
and so after combining with~\eqref{eqn:Acneg}, we have
\begin{equation*}
% \label{eqn:midstep}
    \sqrt{\frac{pn}{2\sigma_{n,1}^2}}\Big(\frac{\TILDEKAPPA}{3}-1\Big) \ = \ \sqrt{\frac{n}{2p\sigma_{n,1}^2}}\sum_{j=1}^p\big(\ell_j+\Delta_j\big) \ + \ o_{\P}(1).
\end{equation*}
In addition, Lemma \ref{lem:NullII_n} shows that
\begin{equation*}
    \sqrt{\frac{n}{2p\sigma_{n,1}^2}}\sum_{j=1}^p\Big(\Delta_j - \big[ \ts\frac{ \kappa}{3}-1 -\ts \frac{4}{n}\big]\Big) \ \xrightarrow{\P} \ 0,
\end{equation*}
which implies that the normalized statistic $T_{n,1}/\sigma_{n,1}$ satisfies
\begin{align}
\label{eqn:limitofmaintermofkappa}
    \frac{T_{n,1}}{\sigma_{n,1}} 
     \ = \ &\sqrt{\frac{n}{2p\sigma_{n,1}^2}}\sum_{j=1}^p \ell_j +   o_{\P}(1).
\end{align}
Finally,
Lemma \ref{lem:NullI_n} shows that
\begin{equation*}
    \sqrt{\frac{n}{2p\sigma_{n,1}^2}}\sum_{j=1}^p \ell_j \ \Rightarrow \ N(0,1),
\end{equation*}
and combining this with~\eqref{eqn:limitofmaintermofkappa} completes the proof.
\qed

\subsubsection{Negligible terms associated with $T_{n,1}$}

Since the proof of Proposition~\ref{prop:T1} is based on a reduction to the case where all the diagonal entries of $\Sigma$ are 1, this condition will be assumed in all of the supporting results for Proposition~\ref{prop:T1} in this and the next subsection. Also, for the following lemma, recall that the event $A_j$ is defined in equation~\eqref{eqn:Ajdef}.

\begin{lemma}
\label{lem:A_j}
    If Assumption \ref{Data generating model} holds 
    and all the diagonal entries of $\SIGMA$ are 1, then
    $$\sum_{j=1}^p\Big(\frac{\tilde \kappa_j}{3}-1\Big)1\{A_j^c\}  \ = \  o_{\P}(1).$$
\end{lemma}

\proof
For any $j=1,\ldots,p$, it can be checked that the inequalities $1 \leq \tilde \kappa_j \leq  n$ hold almost surely, yielding
$$ \sum_{j=1}^p \Big|\frac{\tilde \kappa_j}{3}-1\Big| 1\{A_j^c\} \ \leq \  n \sum_{j=1}^p 1\{A_j^c\}.$$
Hence, it is enough to show $\sum_{j=1}^p \P(A_j^c) =o(\ts \frac{1}{n})$. 
We first apply Chebyshev's inequality to each $\P(A_j^c)$, which leads to
$$\P(A_j^c)  \ \leq \  n^{2/5}\Big\|\ts \frac{2}{n}\sum_{i=1}^{n/2} (\ts x_{ij}^2- 1) \Big\|_{L^8}^{8}.$$
The $L^8$ norm of a sum of centered independent random variables can be bounded with Rosenthal's inequality (Lemma \ref{lem:Rosenthal}) to obtain
\begin{equation}
\begin{split}
\P(A_j^c)
& \ \lesssim \ n^{2/5}  \max\left\lbrace\Big(\ts\frac{2}{n}\var(x_{1j}^2) \Big)^{4}, \Big(\ts\frac{2}{n}\Big)^{7} \|x_{1j}^2- 1\|_{L^8}^{8}  \right\rbrace \\
& \ \lesssim \ n^{-18/5},
\end{split}
\label{eq:P(A_jc)}
\end{equation}
where the last step relies on the condition $\Sigma_{jj}=1$ for all $j=1,\dots,p$, as well as the moment bounds given in \eqref{eq:momentboundofxij}. Thus, we conclude that $\sum_{j=1}^p\P(A_j^c)=o(\ts \frac{1}{n})$, as needed.\qed

For the next lemma, recall that the random variables $\ell_j$, $\Delta_j$, and $\ve_j$ are defined in equations~\eqref{eq:definitionofLj},  \eqref{eq:definitionofdeltaj} and \eqref{eqn:vejef}.

\begin{lemma}
\label{lem:epilson}
    If Assumption \ref{Data generating model} holds and the diagonal entries of $\SIGMA$ are equal to 1, then as $n\to\infty$,
    \begin{align*}
        \frac{1}{\sigma_{n,1}}\sum_{j=1}^p \ve_j\big(\ell_j+\Delta_j\big)  \ \xrightarrow{\P} \ 0.
    \end{align*}
\end{lemma}

\proof
For any $a>0$, consider the following algebraic identity, which arises from a Taylor  expansion of the inverse square function at $a=1$,
$$\frac{1}{a^2} - 1 \ = \ -2(a-1) + (\ts\frac{2}{a} +\frac{1}{a^2})(a-1)^2.$$ 
Letting $\tilde\Sigma_{jj}$ play the role of $a$ leads to the following almost-sure expansion for $\ve_j$,
\begin{equation*}
% \label{eqn:1stexpansion}
\begin{split}
\ve_j & \ = \ 1\{A_j\}\Big(-2(\tilde\Sigma_{jj}-1)+ \big(\ts\frac{2}{\tilde\Sigma_{jj}}+\frac{1}{\tilde\Sigma_{jj}^2}\big)(\tilde\Sigma_{jj}-1)^2\Big)-1\{A_j^c\}1\{\tilde\Sigma_{jj}\neq 0\}.
\end{split}
\end{equation*}
Since the key term in this expansion is $-2(\tilde\Sigma_{jj}-1)$, we define the random variable $\eta_j$ to be the remainder that satisfies
\begin{equation}\label{eqn:etadef}
\ve_j \ = \ -2(\tilde\Sigma_{jj}-1)+\eta_j,
\end{equation}
and therefore, the sum in the stated result may be expressed as
\begin{equation*}
    \frac{1}{\sigma_{n,1}} \sum_{j=1}^p \ve_j \big( \ell_j + \Delta_j \big)  \ = \ -\frac{2}{\sigma_{n,1}}\sum_{j=1}^p (\tilde\Sigma_{jj}-1)(\ell_j+\Delta_j) \ + \ \frac{1}{\sigma_{n,1}}\sum_{j=1}^p \eta_j(\ell_j+\Delta_j).
\end{equation*}
The problem of showing that the first sum on the right is $o_{\P}(1)$ is quite involved, and so we defer the details to Lemma \ref{lem:mainterm}. The rest of this proof will show that $\sum_{j=1}^p \eta_j(\ell_j+\Delta_j)$ is $o_{\P}(1)$, which is sufficient because $\sigma_{n,1}^2 \gtrsim 1$. To do this, we start from the observation that the $L^1$ norm of this sum is upper bounded by $\sum_{j=1}^p \|\eta_j\|_{L^2}\|\ell_j+\Delta_j\|_{L^2}$. Next, from the definition of $\eta_j$ in~\eqref{eqn:etadef}, it follows that
\begin{equation*}
\begin{split}
    \|\eta_j\|_{L^2} \ \lesssim \ \|\tilde\Sigma_{jj}-1\|_{L^4}^2 \ + \ (\P(A_j^c))^{1/4}\|\tilde\Sigma_{jj}-1\|_{L^4} \ + \ (\P(A_j^c))^{1/2}.
\end{split}
\end{equation*}
From (\ref{eq:P(A_jc)}), the order of the third term on the right is $(\P(A_j^c))^{1/2}
\lesssim \ n^{-9/5}$. With regard to the first and second terms, Rosenthal's inequality (recorded in Lemma \ref{lem:Rosenthal}) gives
\begin{equation}\label{eq:4thcentralmomentofsigmatilde}   
\begin{split}
   \|\tilde\Sigma_{jj}-1\|_{L^4} &= \Big\| \ts \frac{2}{n}  \sum_{i=1}^{n/2} (x_{ij}^2 -1)\Big\|_{L^4}  \\
    & \ \lesssim \   \max\Big\{\big(\ts \frac{2}{n} \var(x_{1j}^2)\big)^{1/2},  (\frac{2}{n })^{3/4}  \| x_{1j}^2- 1\|_{L^4} \Big\} \\
    & \ \lesssim \ \frac{1}{\sqrt{n}},
\end{split}
\end{equation}
where we have used the condition that the diagonal entries of $\Sigma$ are equal to 1, as well as moment bounds in \eqref{eq:momentboundofxij}. It follows that  $\|\eta_j\|_{L^2}\lesssim \frac{1}{n}$.

Turning to $\|\ell_j + \Delta_j\|_{L^2}$, the definitions of $\ell_j$ and $\Delta_j$ given in~\eqref{eq:definitionofLj} and \eqref{eq:definitionofdeltaj}, as well as the triangle inequality for the $L^2$ norm yield
\begin{align*}
    \|\ell_j+\Delta_j\|_{L^2} 
& \ \leq \ \sqrt{\ts\frac{2}{n}\var\big( \ts\frac{1}{3} x_{1j}^4  -2 x_{1j}^2\big)}
\ + \ \|\tilde\Sigma_{jj}-1\|_{L^4}^2 \ + \ |r_2 -1|.
\end{align*}
The first term on the right can be handled with \eqref{eq:momentboundofxij}, which shows that if the diagonal entries of $\SIGMA$ are equal to 1, then $\var\big( \ts\frac{1}{3} x_{1j}^4  -2 x_{1j}^2 \big) \lesssim 1.$ Combining with equations \eqref{eq:4thcentralmomentofsigmatilde} and     \eqref{eq:relationofr2} yields $\|\ell_j+\Delta_j\|_{L^2} \lesssim \frac{1}{\sqrt{n}}$.
Finally, the previous bounds on $\|\eta_j\|_{L^2}$ and $\|\ell_j+\Delta_j\|_{L^2}$ show that the norm $\big\|\sum_{j=1}^p \eta_j(\ell_j+\Delta_j)\big\|_{L^1}$ is of order $\frac{1}{\sqrt{n}}$, which completes the proof.\qed

\begin{lemma}
\label{lem:mainterm}
     If Assumption \ref{Data generating model} holds 
    and the diagonal entries of $\SIGMA$ are equal to 1, then as $n \to \infty$
    \begin{align}\label{eqn:lhstwopiece}
         \frac{1}{\sigma_{n,1}}\sum_{j=1}^p \big(\tilde\Sigma_{jj}-1\big)\big(\ell_j+\Delta_j\big) \ \xrightarrow{\P} \ 0.
    \end{align}
\end{lemma}

\proof
It suffices to show that the expectation and variance of the left side of~\eqref{eqn:lhstwopiece} are both $o(1)$. To handle the expectation, observe that
\begin{align*}
    \big(\tilde\Sigma_{jj}-1\big) \big(\ell_j+\Delta_j\big) \ = \ &  \Big(\ts \frac{2}{n}\sum_{i=1}^{n/2} (x_{ij}^2-1)\Big)\Big(\ts \frac{2}{3n}\sum_{i=1}^{n/2} x_{ij}^4 -\big(\ts\frac{2}{n}\sum_{i=1}^{n/2} x_{ij}^2\big)^2\Big),
\end{align*}
and so
\footnotesize
\begin{align*}
    \E\Big( \big(\tilde\Sigma_{jj}-1\big)\big(\ell_j+\Delta_j\big)\Big)  \ = \   & \frac{1 }{3(n/2)^2}\sum_{1\leq i_1, i_2 \leq n/2} \E ( x_{i_1 j}^2x_{i_2 j}^4 ) \ - \ \frac{1}{3(n/2)}\sum_{ i_1 = 1}^{n/2} \E ( x_{i_1 j}^4 ) \\[0.1cm]
&  - \ \frac{1}{(n/2)^3} \sum_{1\leq i_1, i_2, i_3 \leq n/2} \E( x_{i_1 j}^2x_{i_2 j}^2x_{i_3 j}^2 )  \ + 
 \ \frac{1}{(n/2)^2}\sum_{1\leq i_1, i_2 \leq n/2}  \E( x_{i_1 j}^2x_{i_2 j}^2).
\end{align*}
\normalsize
When the diagonal entries of $\SIGMA$ are equal to 1, the expressions for $\x_i$ and $r_k$ in \eqref{eq:xi} and \eqref{eqn:extrark} imply $\E ( x_{1 j}^4 ) = 3 r_2$, $\E ( x_{1 j}^6 ) = 15 r_3$, as well as
\small
\begin{align*}
    \sum_{1\leq i_1, i_2 \leq n/2} \E ( x_{i_1 j}^2x_{i_2 j}^4 ) \ = \  &\ts \frac{n}{2} \E( x_{1 j}^6 ) +\frac{n}{2} (\frac{n}{2}-1)\E( x_{1 j}^4 )\E( x_{1 j}^2 ) \\ \ = \  &\ts \frac{15n}{2}r_3 + \frac{3n}{2}(\frac{n}{2}-1)r_2,\\[0.2cm]
    \sum_{1\leq i_1, i_2, i_3 \leq n/2} \E( x_{i_1 j}^2x_{i_2 j}^2x_{i_3 j}^2 ) \ = \  & \ts \frac{n}{2} \E( x_{1 j }^6 ) +  \frac{3n}{2} (\frac{n}{2}-1)\E( x_{1 j }^4 )\E( x_{1 j }^2 )+ \frac{n}{2} (\frac{n}{2}-1)(\frac{n}{2}-2) \E( x_{1 j }^2 )^3 \\  \ = \  &\ts \frac{15n}{2}r_3+ \frac{9n}{2}(\frac{n}{2}-1)r_2 + \frac{n}{2}(\frac{n}{2}-1)(\frac{n}{2}-2),\\[0.2cm]
   \sum_{1\leq i_1, i_2 \leq n/2}  \E( x_{i_1 j}^2x_{i_2 j}^2 )  \ = \   & \ts \frac{n}{2} \E( x_{1 j }^4 ) +\frac{n}{2} (\frac{n}{2}-1)\E( x_{1 j }^2 ) ^2  \\ 
     \ = \   &\ts \frac{3n}{2}r_2 + \frac{n}{2}(\frac{n}{2}-1).
\end{align*}
\normalsize
Combining with Lemma \ref{lem:moment}, which shows that $r_2$ and $r_3$ are both equal to $1+o(p^{-3/4})$, we have
\footnotesize
\begin{align*}
    \E\Big( \big( \tilde\Sigma_{jj}-1\big) \big(\ell_j+\Delta_j\big)\Big) \ = \ & \Big(r_2 + \frac{5r_3- r_2}{n/2}\Big) -r_2- \Big( 1 + \frac{9r_2-3}{n/2} + \frac{15 r_3 - 9 r_2 +2}{(n/2)^2}\Big)   + \Big( 1  +\frac{3 r_2-1}{n/2} \Big)
     \\= \  &o(\ts\frac{1}{n}).
\end{align*}
\normalsize
So, by summing over $j$ and noting that $\sigma_{n,1}^2\gtrsim 1$, we conclude
$$\frac{1}{\sigma_{n,1}}\sum_{j=1}^p  \E\Big(\big(\tilde\Sigma_{jj}-1\big)\big(\ell_j+\Delta_j\big)\Big) = o(1).$$

Next, we show that the variance of the left side of~\eqref{eqn:lhstwopiece} is $o(1)$. For this purpose, we begin by noting that
\footnotesize
\begin{align*}
     \sum_{j=1}^p\big(\tilde\Sigma_{jj}-1\big)\big(\ell_j+\Delta_j\big)
     \ = \ &\sum_{j=1}^p \big(\tilde\Sigma_{jj}-1\big)\Big(\ts \frac{2}{3n} \sum_{i=1}^{n/2} (x_{ij}^4 -  \E(x_{1j}^4))-2(\tilde\Sigma_{jj} -1) -(\tilde\Sigma_{jj} -1)^2 +r_2 -1 \Big). 
\end{align*}
\normalsize
Due to the basic inequality $\var\!\big ( \sum_{l=1}^4 V_l \big) \leq 4 \sum_{l=1}^4 \var(V_l)$
that holds for generic random variables $V_1,\dots,V_4$, as well as the relation $\sigma_1^2=\frac{8}{3p}\|R\|_4^4$, it is sufficient to establish the following four limits,
\begin{align}
\label{eq:maintermvaraince1}
    & \var\Bigg(\sum_{j=1}^p \big(\tilde\Sigma_{jj}-1\big) \Big(\ts \frac{1}{3(n/2)} \displaystyle\sum_{i=1}^{n/2}\big(x_{ij}^4 -  \E(x_{1j}^4)\big)\Big)\Bigg) \ = \ o(\|R\|_4^4/p),\\[0.1cm]
    \label{eq:maintermvaraince2}
    & \var\bigg(\sum_{j=1}^p \big(\tilde\Sigma_{jj}-1\big) ^2\bigg)  \ = \ o(\|R\|_4^4/p),\\[0.1cm]
    \label{eq:maintermvaraince3}
    & \var\bigg(\sum_{j=1}^p \big(\tilde\Sigma_{jj}-1\big) ^3\bigg)  \ = \ o(\|R\|_4^4/p),\\[0.1cm]
    \label{eq:maintermvaraince4}
    & (r_2-1)^2 \var\bigg(\sum_{j=1}^p\big(\tilde\Sigma_{jj}-1\big) \bigg)  \ = \  o(\|R\|_4^4/p).
\end{align}
\noindent\emph{Proof of~\eqref{eq:maintermvaraince1}.} The definition of $\tilde \Sigma_{jj}$ implies that the left hand side of~\eqref{eq:maintermvaraince1} is equal to
\begin{equation}\label{eq:var1}
    \begin{split} 
     \frac{1}{9(n/2)^4}\sum_{1 \leq i_1,\ldots,i_4 \leq n/2} \sum_{j,k = 1}^p \cov\bigg(\big(x_{i_1 j}^4 -  \E(x_{1 j}^4)\big)\big(x_{i_2 j}^2 -  1\big)\,,\,\big(x_{i_3 k}^4 -  \E(x_{1 k}^4)\big)\big(x_{i_4k}^2 -  1\big) \bigg).         
    \end{split}
\end{equation}
The covariance terms in \eqref{eq:var1} are non-zero only when $\big|\{i_1, i_2\} \bigcap \{i_3, i_4\}\big|$ is equal to 1 or 2. When $\big|\{i_1, i_2\} \bigcap \{i_3, i_4\}\big|  =2$, there are two possible types of tuples $(i_1,\dots,i_4)$, corresponding to $i_1 = i_3 \neq i_2 = i_4$ and $i_1 = i_4 \neq i_2 = i_3$, with the number of tuples of each type being $\mathcal{O}(n^2)$. When the covariance terms corresponding to the first type are summed over $(j,k)$, Lemmas \ref{lem:covariances} and \ref{lem:moment} show that
%imply when the diagonal entries of $\SIGMA$ are equal to 1, the corresponding summation of the covariance in \eqref{eq:var1} is \red{I got rid of $=$ and put $\lesssim$} 
\begin{align*}
    \sum_{j,k = 1}^p \big|\cov( x_{1j}^4 ,x_{1k}^4) \cov(x_{1j}^2 ,x_{1k}^2)\big| \ \lesssim \ & \,  p^2|(r_4 - r_2^2)(r_2-1)|+ \|R\|_2^2\\
     \ = \ & o(p\|R\|_4^4),
\end{align*}
where the last step use the general fact that $\|R\|_2^2 \lesssim p^{1/2}\|R\|_4^4$. Likewise, these lemmas show that when the covariance terms corresponding to the second type of tuples $(i_1,\dots,i_4)$ are summed over $(j,k)$ we have
\begin{align*}
    \sum_{j,k = 1}^p \big|\cov(x_{1j}^4 ,x_{1k}^2) \cov(x_{1j}^2 ,x_{1k}^4)\big| \ \lesssim \ & \, p^2(r_3-r_2)^2+ \|R\|_2^2\\
     \ = \ & o(p\|R\|_4^4).
\end{align*}
When $\big|\{i_1, i_2\} \bigcap \{i_3, i_4\}\big| = 1$, there are three possible types of tuples $(i_1,\dots,i_4)$, in which there are 1, 2, or 3 distinct indices. If there are 2 or 3 distinct indices among $(i_1,\dots,i_4)$, then the associated covariance terms in~\eqref{eq:var1} are zero. With regard to tuples $(i_1,\dots,i_4)$ in which all indices are equal, there are $\mathcal{O}(n)$ such tuples. Also, equation \eqref{eq:momentboundofxij} implies that the associated covariance terms in \eqref{eq:var1} are $\mathcal{O}(1)$, uniformly over $(j,k)$. So, the sum of these covariance terms over $(j,k)$ is $\mathcal{O}(p^2)$. Combining the last several steps shows that the total sum in~\eqref{eq:var1} is $\mathcal{O}(n^{-4}[n^2o(p\|R\|_4^4)+ n p^2])=o(\|R\|_4^4/p)$, which proves (\ref{eq:maintermvaraince1}).

\noindent \emph{Proof of~\eqref{eq:maintermvaraince2}.} We proceed analogously to the work above. The left hand side of~\eqref{eq:maintermvaraince2} is equal to
\begin{align}
\label{eq:var2}
    \frac{1}{(n/2)^4} \sum_{1 \leq i_1, i_2, i_3,i_4 \leq n/2} \, \sum_{j,k = 1}^p  \,\cov\Big((x_{i_1 j}^2 -  1)(x_{i_2 j}^2 -  1),(x_{i_3 k}^2 -  1)(x_{i_4 k}^2 -  1) \Big) .
\end{align}
As before, it is enough to restrict our attention to the tuples $(i_1,\dots,i_4)$ for which the cardinality $\big|\{i_1, i_2\} \bigcap \{i_3, i_4\}\big|$ is equal to 1 or 2.
When $\big|\{i_1, i_2\} \bigcap \{i_3, i_4\}\big|  =2$, it follows from equation~\eqref{eq:relationofr2} and Lemma~\ref{lem:covariances} that the sum of the corresponding covariance terms over $(j,k)$ is equal to
\begin{align*}
    \sum_{j,k = 1}^p \Big(\cov(x_{1j}^2 ,x_{1k}^2)\Big)^2 \ \lesssim \  & \, p^2(r_2-1)^2+ \|R\|_2^2\\
     \ = \ & o(p\|R\|_4^4).
\end{align*}
When $\big|\{i_1, i_2\} \bigcap \{i_3, i_4\}\big|  =1$, the associated tuples $(i_1,\dots,i_4)$ may have 1, 2, or 3 distinct indices. If there are 2 or 3 distinct indices, then the corresponding covariance terms in~\eqref{eq:var2} are zero, and so it is enough to consider the tuples for which all indices are equal. Due to equation \eqref{eq:momentboundofxij}, the covariance terms in~\eqref{eq:var2} corresponding to these tuples are $\mathcal{O}(1)$, unformly over $(j,k)$, and so the sum of these terms over $(j,k)$ is $\mathcal{O}(p^2)$. Altogether, it follows that the total sum in~\eqref{eq:var2} is $\mathcal{O}(n^{-4}[n^2o(p\|R\|_4^4)+n p^2])=o(\|R\|_4^4/p)$, which proves~\eqref{eq:maintermvaraince2}.

\noindent \emph{Proof of \eqref{eq:maintermvaraince3}.} Using the triangle inequality for the $L^2$ norm, followed by Rosenthal's inequality (Lemma \ref{lem:Rosenthal}) and equation \eqref{eq:momentboundofxij}, the left hand side of (\ref{eq:maintermvaraince3}) can be bounded as
\begin{align*}
\var\bigg(\sum_{j=1}^p \big(\tilde\Sigma_{jj}-1\big) ^3\bigg) & \ \leq \ \bigg(\sum_{j=1}^p \|\tilde\Sigma_{jj}-1\|_{L^6}^3\bigg)^2\\
%
    %\Big\| \ts \frac{2}{n} \sum_{i=1}^{n/2} (x_{ij}^2 -1)\Big\|_{L^6} \ 
    & \ \lesssim \  p^2\bigg(\max_{1\leq j\leq p}\max \Big\{\big(\ts \frac{1}{n/2} \var(x_{1j}^2)\big)^{1/2},  (\frac{1}{n/2})^{-5/6} \| x_{1j}^2- 1\|_{L^6} \Big\}\bigg)^6\\
    & \ \lesssim \  \frac{1}{n},
\end{align*}
which implies~\eqref{eq:maintermvaraince3}.

\noindent \emph{Proof of \eqref{eq:maintermvaraince4}.} Similarly to the previous argument, the left hand side of and (\ref{eq:maintermvaraince4}) can be bounded as
\begin{align*}   
% \label{eq:var3}
(r_2-1)^2 \var\bigg(\sum_{j=1}^p\big(\tilde\Sigma_{jj}-1\big) \bigg) & \ \lesssim  \ \frac{1}{p^2}
\bigg(\sum_{j=1}^p \|\tilde\Sigma_{jj}-1\|_{L_2}\bigg)^2\\
& \ = \  \frac{1}{p^2}
\bigg(\sum_{j=1}^p \sqrt{\ts\frac{2}{n}\var(x_{1j}^2)}\bigg)^2\\
    &\ \lesssim \ \frac{1}{n},
\end{align*}
which implies~\eqref{eq:maintermvaraince4} and completes the proof of the lemma.\qed

\subsubsection{Main terms of $T_{n,1}$}
    
\begin{lemma}
\label{lem:NullI_n}
    If Assumption \ref{Data generating model} holds and the diagonal entries of $\SIGMA$ are equal to 1, then as $n\to\infty$,
    \begin{align*}
        \sqrt{\frac{n}{2p\sigma_{n,1}^2}}\sum_{j=1}^p \ell_j \ \Rightarrow \ N(0,1).
    \end{align*} 
\end{lemma}
\proof
   We begin the proof by expressing $\sum_{j=1}^p\ell_j$ in terms of a sum of centered i.i.d.~random variables $w_1,\dots,w_n$ defined according to
    \begin{align*} 
    % \label{eq:defofwi}
        w_i  \ = \ 
    \frac{1}{\sqrt{p}}\sum_{j=1}^p \ts\frac{1}{3}\big(x_{ij}^4 - \E(x_{1j}^4)\big) \ - \ 2 (x_{ij}^2-1).
    \end{align*}
   In this notation, we have
    $$\sqrt{\frac{n}{2p}}\sum_{j=1}^p \ell_j  \ = \  \sqrt{\frac{2}{n}}   \sum_{i=1}^{n/2} w_i.$$
    To determine the variance of $w_1$, a direct calculation gives
    \begin{align*}
        \var\left(w_{1}\right)
        & \ = \  \frac{1}{p} \sum_{j,k = 1}^p  \Big(\ts\frac{1}{9}\cov(x_{1j}^4,x_{1k}^4) + 4\cov(x_{1j}^2,x_{1k}^2) - \ts\frac{2}{3}\big(\cov(x_{1k}^4,x_{1j}^2) + \cov(x_{1k}^2,x_{1j}^4)\big)\Big).
        \end{align*}
Each of the covariances appearing on the right is calculated explicitly in Lemma \ref{lem:covariances}, which leads to
        \begin{align}\label{eqn:correction_motivation}
        \var\left(w_{1}\right) & \ = \  \frac{1}{p}\Big(\ts \frac{8r_4}{3} \|R\|_4^4+ 8(r_4+r_2-2r_3) \|R\|_2^2\Big) +p(r_4  - 4r_3 - r_2^2+8r_2-4).
    \end{align}
    Furthermore, Lemma \ref{lem:moment} shows that
    $r_k = 1 + o(p^{-3/4})$ for $k=1,\ldots,4,$
    as well as
    \begin{equation*}
                p(r_4  - 4r_3 - r_2^2+8r_2-4) \ = \ o(1).
    \end{equation*}
    So, combining with the general fact that $\|R\|_2^2 \lesssim p^{1/2}\|R\|_4^4$, and recalling that $\sigma_{n,1}^2=\frac{8}{3p}\|R\|_4^4$, we have
\begin{equation}\label{eqn:sig1equiv}
    \begin{split}
        % \label{eq:rho1negligibleterm}
   \frac{\var(w_{1})}{\sigma_{n,1}^2}
    & \ = \ 1+o(1).
    \end{split}
\end{equation}
To complete the proof, we apply the Lindeberg CLT for triangular arrays given in \cite[Prop.~2.27]{van2000asymptotic}. To verify the Lindeberg condition, it is sufficient to show that the following Lyapunov-type condition holds,
\begin{align}
\label{eq:lindeberg}
    \E (w_1^4) \ = \  o(n\sigma_{n,1}^4).
\end{align}
Verifying this condition is by far the most challenging aspect of the proof and depends on delicate cancellations of high-order terms in an expansion of the fourth moment of $w_1$. To proceed, define the random variable
\begin{equation}\label{eqn:yjdef}
y_j = \ts\frac 13 x_{1j}^4-2x_{1j}^2
\end{equation} 
for each $j=1,\dots,p$. In this notation, we have
\begin{align*}
    p^2\E(w_1^4)  
    \ = \ &   \bigg\|\sum_{j=1}^p (y_j+2 - r_2)\bigg\|_{L^4}^4\\
    \ = \ & \sum_{j_1,j_2,j_3,j_4=1}^p   \!\!\!\!\!\!\E(y_{j_1}y_{j_2}y_{j_3}y_{j_4} )
    \ + \  4(2-r_2)p\!\!\!\!\!\!\sum_{j_1,j_2,j_3=1}^p\!\!\!\!\E(y_{j_1}y_{j_2}y_{j_3} )\\
     & \ + \ 6(2-r_2)^2p^2\!\!\!\sum_{j_1,j_2=1}^p\E(y_{j_1}y_{j_2})
     \ - \  3(2-r_2)^4p^4.
\end{align*}
The last term on the right side can be simplified with respect to its dependence on $r_2$ by using the expansion in \eqref{eq:relationofr2}, which yields
\begin{align*}
    - 3(2-r_2)^4p^4   \ = \ &(4r_2-7)p^4 + 8(\tau-2) p^3 +o(p^3).
\end{align*}
We will use Isserlis' theorem \citep[Theorem 1.28]{janson1997gaussian} to evaluate the three sums in the previous formula for $p^2\E(w_1^4)$, which is a substantial undertaking that we defer to Proposition \ref{prop:combinationofhighmoments}. This proposition establishes the following three equations
\footnotesize
 \begin{align*}
     \sum_{j_1,j_2,j_3,j_4=1}^p   \E(y_{j_1}y_{j_2}y_{j_3}y_{j_4} ) & \ = \ p^4(r_8-8r_7+24r_6-32r_5+16r_4) + 16p^2\|R\|_4^4+ o(p\|R\|_4^8),\\
     4(2-r_2)p\sum_{j_1,j_2,j_3=1}^p\E(y_{j_1}y_{j_2}y_{j_3} ) & \ = \  p^4(4r_6-24r_5+48r_4-32r_3) + 4 (\tau-2) p^3 - 32 p^2\|R\|_4^4+ o(p\|R\|_4^8),\\
     6(2-r_2)^2p^2\sum_{j_1,j_2=1}^p\E(y_{j_1}y_{j_2}) & \ = \  p^4 (6r_4-24r_3 +24r_2)-12(\tau-2) p^3 +  16 p^2\|R\|_4^4+ o(p\|R\|_4^8).
\end{align*}
\normalsize
Based on the relations among $r_1,\dots,r_8$ established in Lemma \ref{lem:moment}, the last several equations can be combined to obtain
\begin{align*}
    p^2\E(w_1^4) \ = \ & o(p\|R\|_4^8),
\end{align*}
which implies desired condition~\eqref{eq:lindeberg}, since $\sigma_{n,1}^2=\frac{8}{3p}\|R\|_4^4$. \qed

\begin{lemma}\label{lem:NullII_n}
    If Assumption \ref{Data generating model} holds and the diagonal entries of $\SIGMA$ are 1, then 
    $$\frac{1}{\sigma_{n,1}}\sum_{j=1}^p\Big(\Delta_j - \big[ \ts\frac{ \kappa}{3}-1 -\ts \frac{4}{n}\big]\Big) \ \xrightarrow{\P} \ 0.$$ 
\end{lemma}
\proof
    First, define the random variable
        $$m_n \ = \ \sum_{j=1}^p\Big(\ts\frac{2}{n} \sum_{i=1}^{n/2} (x_{ij}^2-1)\Big)^2,$$
    and observe that 
    \begin{align*} 
         \sum_{j=1}^p\Big(\Delta_j - \big[ \ts \frac{ \kappa}{3}-1 - \frac{4}{n}\big]\Big) \ = \ \frac{4p}{n}-m_n.
    \end{align*}
Since $\sigma_{n,1}^2 \gtrsim 1$, it is enough to show that $\E(m_n)=4p/n+o(1)$ and $\var(m_n)=o(\sigma_{n,1}^2)$. Noting that $\sigma_{n,1}^2=(8/3)\|R\|_4^4/p$, the condition $\var(m_n)=o(\sigma_{n,1}^2)$ follows directly from~\eqref{eq:maintermvaraince2}. To calculate the expectation of $m_n$, note that the relations $\E(x_{1j}^2)=1$  and $\var(x_{1j}^2) = 3 r_2 - 1$ hold for all $j=1,\dots,p$ when the diagonal entries of $\SIGMA$ are 1, and equation \eqref{eq:relationofr2} implies $3r_2-1=2 + o(1)$. Therefore, we have
\begin{align*}
    \E(m_n)&\ = \  \frac{2}{n} \displaystyle\sum_{j=1}^p\var(x_{1j}^2)\\
    & \ = \ \frac{4p}{n} + o(1),
\end{align*}
which completes the proof.
\qed

\begin{proposition}
\label{prop:combinationofhighmoments}
Suppose that Assumption \ref{Data generating model} holds and the diagonal entries of $\SIGMA$ are equal to 1. Then, the following equations hold,
\footnotesize
\begin{align}
\label{eq:4terms}
     \sum_{j_1,j_2,j_3,j_4=1}^p   \E(y_{j_1}y_{j_2}y_{j_3}y_{j_4} ) & \ = \ p^4 (r_8-8r_7+24r_6-32r_5+16r_4 ) + 16p^2\|R\|_4^4+ o(p\|R\|_4^8),\\
     \label{eq:3terms}
     (2-r_2)p\sum_{j_1,j_2,j_3=1}^p\E(y_{j_1}y_{j_2}y_{j_3} ) & \ = \  p^4 (r_6-6r_5+12r_4-8r_3 ) +(\tau-2) p^3 - 8 p^2\|R\|_4^4+ o(p\|R\|_4^8),\\
     \label{eq:2terms}
     (2-r_2)^2p^2\sum_{j_1,j_2=1}^p\E(y_{j_1}y_{j_2}) & \ = \  p^4 (r_4-4r_3 +4r_2 )-2(\tau-2) p^3 +  \ts\frac{8}{3} p^2\|R\|_4^4+ o(p\|R\|_4^8).
\end{align}
\normalsize
\end{proposition}

\subsubsection*{Proof of Proposition~\ref{prop:combinationofhighmoments}, equation~\eqref{eq:4terms}.}
Recall that $y_j = \ts\frac 13 x_{1j}^4-2x_{1j}^2$ for $j=1,\dots,p$, which leads to the following algebraic relation
\small
\begin{equation}
\label{eq:decomposeof4terms}
\begin{split}
    \sum_{j_1,j_2,j_3,j_4=1}^p   \E(y_{j_1}y_{j_2}y_{j_3}y_{j_4}) \ = \ \  \ \ \  \ \frac{1}{3 ^4} &\sum_{j_1,j_2,j_3,j_4=1}^p \E\big( x_{1j_1}^{4} x_{1j_2}^{4} x_{1j_3}^{4}  x_{1j_4}^{4}  \big)\\
 - \frac{4\cdot2}{3 ^3} &\sum_{j_1,j_2,j_3,j_4=1}^p \E\big( x_{1j_1}^{4} x_{1j_2}^{4} x_{1j_3}^{4}  x_{1j_4}^{2}  \big)   \\
    + \frac{6\cdot2^2}{3 ^2} &\sum_{j_1,j_2,j_3,j_4=1}^p \E\big( x_{1j_1}^{4} x_{1j_2}^{4} x_{1j_3}^{2}  x_{1j_4}^{2}  \big)\\
    - \frac{4\cdot2^3}{3 } &\sum_{j_1,j_2,j_3,j_4=1}^p \E\big( x_{1j_1}^{4} x_{1j_2}^{2} x_{1j_3}^{2}  x_{1j_4}^{2}  \big) \\
    + 2^4 &\sum_{j_1,j_2,j_3,j_4=1}^p \E\big( x_{1j_1}^{2} x_{1j_2}^{2} x_{1j_3}^{2}  x_{1j_4}^{2}  \big).
\end{split}
\end{equation}
\normalsize

The proof consists in developing a way to systematically calculate the sums on the right side of~\eqref{eq:decomposeof4terms}.
When the diagonal entries of $\SIGMA$ are 1, the random vector
$\boldsymbol\zeta_1=\SIGMA^{1/2}\z_1$ defined in (\ref{eq:definitionofzeta}) is equivalent
to $\boldsymbol \zeta_1=R^{1/2}\z_1$, where $\z_1$ is a standard normal vector independent of $\xi_1$.
Likewise, we may write
$\x_1 = \frac{\xi_1}{\|\z_1\|_2} \boldsymbol\zeta_1,$
and so for any positive integers $k_1,\dots,k_4$ we have
\begin{align*} 
    \E\big( x_{1j_1}^{k_1} x_{1j_2}^{k_2} x_{1j_3}^{k_2}  x_{1j_4}^{k_4}  \big) &\  
    =  \ r_{\frac{k_1+ k_2 +k_3+ k_4}{2}}
    \E\big( \zeta_{1j_1}^{k_1} \zeta_{1j_2}^{k_2} \zeta_{1j_3}^{k_3}  \zeta_{1j_4}^{k_4}  \big). 
\end{align*}
To compute the expectation on the right, we will apply Isserlis' theorem \citep[Theorem 1.28]{janson1997gaussian}, which states that the following relation holds for any centered $d$-dimensional normal vector $(g_1,\dots,g_d)$,
\begin{align}\label{eqn:isserlis}
    \E(g_1\cdots g_d) \ = \ \sum_{\pi \in \mathcal{P}_d} \prod_{\{j,j'\}\in \pi} \E(g_jg_{j'}),
\end{align}
where $\mathcal{P}_d$ is the collection of all partitions $\pi$ of $\{1,\dots,d\}$ into distinct pairs $\{j,j'\}$. (Such partitions are also called ``pair partitions''.) More specifically, we will apply Isserlis' theorem in the case when $d=\sum_{l=1}^4 k_l$ and
\begin{equation}\label{eqn:16indices}
(g_1,\dots,g_d) \ = \ (\underbrace{\zeta_{1j_1},\ldots, \zeta_{1j_1}}_{k_1}, \underbrace{\zeta_{1j_2},\ldots, \zeta_{1j_2}}_{k_2}, \underbrace{\zeta_{1j_3},\ldots, \zeta_{1j_3}}_{k_3},\underbrace{\zeta_{1j_4},\ldots, \zeta_{1j_4}}_{k_4}).
\end{equation}
Fortunately, in this context, if a factor $\E(g_jg_{j'})$ does not correspond to one of the diagonal entries of $R$, then it can only take one of six possible forms: $R_{j_1,j_2}$, $R_{j_1,j_3}$, $R_{j_1,j_4}$, $R_{j_2,j_3}$, $R_{j_2,j_4}$, $R_{j_3,j_4}$. Therefore, Isserlis' theorem gives
\begin{equation}\label{eqn:applyIsserlis}
     \E\big( x_{1j_1}^{k_1} x_{1j_2}^{k_2} x_{1j_3}^{k_2}  x_{1j_4}^{k_4}  \big) 
     \  =  \ r_{\frac{k_1+ k_2 +k_3+ k_4}{2}}  \sum_{\mathbf{e} \in \mathcal{E}} c_{\mathbf{e}} \,\!R_{j_1,j_2}^{e_1}R_{j_1,j_3}^{e_2}R_{j_1,j_4}^{e_3}R_{j_2,j_3}^{e_4}R_{j_2,j_4}^{e_5}R_{j_3,j_4}^{e_6},
\end{equation}
where $\mathcal{E}$ is the set of all possible tuples of exponents $\mathbf{e}=(e_1,\dots,e_6)$ for the six factors $R_{j_1,j_2},\dots, R_{j_3,j_4}$ and $c_{\mathbf{e}}$ is the number of pair partitions of $\{1,\dots,k_1+k_2+k_3+k_4\}$ corresponding to $\mathbf{e}$ in Isserlis' theorem.

To facilitate the summation over $\mathcal{E}$, we define an equivalence relation on its elements. We say that $\mathbf{e}$ is equivalent to another tuple $\mathbf{e}'=(e_1',\dots,e_6')$ if the product $R_{j_1,j_2}^{e_1}\cdots R_{j_3,j_4}^{e_6}$ can be obtained from $R_{j_1,j_2}^{e_1'}\cdots R_{j_3,j_4}^{e_6'}$ by permuting $j_1,\dots,j_4$. Hence, if $\mathbf{e}$ and $\mathbf{e}'$ are equivalent, then their corresponding terms in~\eqref{eqn:applyIsserlis} will produce the same value when they are summed over $j_1,\dots, j_4$. Also, let $[\mathbf{e}]$ denote the equivalence class of $\mathbf{e}$, and define 
$$C([\mathbf{e}]) \ = \ \sum_{\mathbf{e}\in[\mathbf{e}]} c_{\mathbf{e}}$$ 
as well as
\begin{equation*}
% \label{eq:ai}
   \begin{split}
       \sps([\mathbf{e}]) \ = \ \sum_{j_1,j_2,j_3,j_4=1}^p R_{j_1,j_2}^{e_1} R_{j_1,j_3}^{e_2} R_{j_1,j_4}^{e_3} R_{j_2,j_3}^{e_4} R_{j_2,j_4}^{e_5} R_{j_3,j_4}^{e_6}.
   \end{split}
\end{equation*} 
 Then, the following formula results from summing~\eqref{eqn:applyIsserlis} over $j_1,\dots,j_4$,
\begin{align*}
    \sum_{j_1,j_2,j_3,j_4=1}^p \E\big( x_{1j_1}^{k_1} x_{1j_2}^{k_2} x_{1j_3}^{k_2}  x_{1j_4}^{k_4}  \big) 
    &\  =  \ r_{\frac{k_1+ k_2 +k_3+ k_4}{2}}  \sum_{[\mathbf{e}]\in [\mathcal{E}]} C([\mathbf{e}]) \sps([\mathbf{e}]),
\end{align*}
where $[\mathcal{E}]$ denotes the collection of equivalence classes.

To choose a unique representative for each equivalence class, we use the largest element with respect to the lexicographical ordering on tuples. A complete system of representatives is displayed in the second column of Table \ref{table:coefficient}. (More precisely, since $\mathcal{E}$ depends on $(k_1,k_2,k_3,k_4)$, each possible choice of $(k_1,k_2,k_3,k_4)$ induces a separate system of representatives. For a column of Table~\ref{table:coefficient} that is headed by a particular choice of $(k_1,k_2,k_3,k_4)$, the representatives correspond to the values of $\mathbf{e}$ associated with the non-empty entries in that column.)

For the equivalence class $[\mathbf{e}]$ represented by the tuple $\mathbf{e}$ in the $i^{\text{th}}$ row of Table~\ref{table:coefficient}, we define the shorthand
\begin{equation}\label{eqn:sidef}
\sps_i \ = \ \sps([\mathbf{e}]).
\end{equation}
In the case when $(k_1,k_2,k_3,k_4)=(4,4,4,4)$, Lemma~\ref{lem:neglibileterms} shows that the sums $\sps_i$ for $i=7$, 8, 10, 11, 12, 13, 15, 17, 19, 20,
are of order $o(p\|R\|_4^8)$. As will be seen below, it follows that these sums play a negligible role in the calculation of~\eqref{eq:4terms}. Hence, it is not necessary to calculate the associated coefficients $C([\mathbf{e}])$, and these are marked by \ding{56} in Table~\ref{table:coefficient}.

To illustrate how the other coefficients $C([\mathbf{e}])$ can be calculated, we give two detailed examples in the case when $(k_1,k_2,k_3,k_4)=(4,4,4,4)$. First, the choice $\mathbf{e}=(0,0,0,0,0,0)$ corresponds to terms in~\eqref{eqn:applyIsserlis} that are equal to 1, which arise from partitioning the 16 indices in~\eqref{eqn:16indices} as
$$\{j_1,j_1\}, \{j_1,j_1\}, \{j_2,j_2\}, \{j_2,j_2\}, \{j_3,j_3\}, \{j_3,j_3\}, \{j_4,j_4\}, \{j_4,j_4\}.$$
(Note that each pair $\{j_l,j_l\}$ corresponds to a factor $R_{j_l,j_l}=1$.)
The number of such partitions is $c_{\mathbf{e}}=3^4$, because there are $k_l=4$ copies of each ``type'' $j_l$ that must be partitioned into pairs,  and any set of 4 elements has exactly 3 pair partitions. Furthermore, the equivalence class of $\mathbf{e}=(0,0,0,0,0,0)$ is a singleton, and so $C([\mathbf{e}])= c_{\mathbf{e}}=81$.

As a second example of how to calculate $C([\mathbf{e}])$, consider the choice $\mathbf{e} = (1,1,0,0,1,1)$, whose equivalence class is equal to $[\mathbf{e}] = \{(1,1,0,0,1,1), (1,0,1,1,0,1), (0,1,1,1,1,0)\}$. This choice of $\mathbf{e}$ corresponds to terms in~\eqref{eqn:applyIsserlis} of the form $R_{j_1,j_2}R_{j_1,j_3}R_{j_2,j_4}R_{j_3,j_4}$, which arise from partitioning the 16 indices in~\eqref{eqn:16indices} as
\begin{align*}
    \{j_{1},j_{1}\}, \{j_{2},j_{2}\}, \{j_{3},j_{3}\}, \{j_{4},j_{4}\}, \{j_{1},j_{2}\}, \{j_{1},j_{3}\}, \{j_{2},j_{4}\}, \{j_{3},j_{4}\}.
\end{align*}
(Note that the first four pairs above lead to unit factors, $R_{j_l,j_l}=1$ for $l=1,\dots,4$, and so these pairs are ``invisible'' in the product $R_{j_1,j_2}R_{j_1,j_3}R_{j_2,j_4}R_{j_3,j_4}$.)
The number of ways to form such a partition is $c_{\mathbf{e}}={4 \choose 2} {4 \choose 2} {4 \choose 2} {4 \choose 2} 2^4$, 
and the equivalence class $[\mathbf{e}]$ has 3 elements, so we conclude that $C([\mathbf{e}]) = 3c_{\mathbf{e}} = 62208$.  Following similar reasoning, we can calculate the coefficients $C([\mathbf{e}])$ associated with the other choices of $\mathbf{e}$, as given in Table~\ref{table:coefficient}.

\begin{table}[H]
\centering
\begin{tabular}{cc|ccccc}
                  &             & \multicolumn{5}{c}{$(k_1,k_2,k_3,k_4)$}               \\ \cline{3-7} 
                  &       $\mathbf{e}$           &  $(4,4,4,4)$  & $(4,4,4,2)$  &  $(4,4,2,2)$  & $(4,2,2,2)$ & $(2,2,2,2)$ \\
                  \hline
           $1$ & $(0,0,0,0,0,0)$ &    81  & 27 & 9 & 3 &  1                \\
                  \hline 
           $2$  & $(1,1,0,0,1,1)$ &    62208 & 10368 &  1728  &  288   & 48 \\
                     \hline
            $3$  & $(1,1,0,1,0,0)$ &    20736  & 4320 & 864 & 168   &  32  \\
              \hline    
            $4$ & $(2,0,0,0,0,0)$ &    3888 & 972 & 234 & 54 &  12  \\
 \hline
                 
            $5$ & $(2,0,0,0,0,2)$ &    15552 &  2592 & 432 & 72 &   12  \\
                  \hline
            $6$ & $(2,1,1,0,0,1)$ &    248832 & 31104 & 3456 & 288   &   \\
             \hline     
            $7$ & $(2,1,1,1,1,0)$ &    \ding{56} &    \ding{56} &    \ding{56} &  &    \\
                  \hline
            $8$ & $(2,1,1,1,1,2)$ &     \ding{56} &  &    &    &                          \\
                  \hline
            $9$ & $(2,2,0,0,0,0)$ &    31104 & 5184 & 720 &  72     &                       \\
                     \hline
                 
            $10$ & $(2,2,0,0,2,0)$ &     \ding{56} &     \ding{56} &     \ding{56} &    &                          \\
                  \hline
            $11$ & $(2,2,0,0,2,2)$ &     \ding{56} &  &          &     &                    \\
                  \hline
            $12$ & $(2,2,0,1,1,1)$ &     \ding{56} &     \ding{56} &   &     &                        \\
                  \hline
            $13$ & $(2,2,0,2,0,0)$ &     \ding{56} &     \ding{56} &  &     &                        \\
                  \hline
            $14$ & $(3,1,0,0,1,1)$ &    165888 &  13824 & 768 &     &                         \\
                  \hline
            $15$ & $(3,1,0,0,1,3)$ &     \ding{56} &  &   &      &                      \\
                  \hline
            $16$ & $(3,1,0,1,0,0)$ &    41472 & 5184 &  384 &     &                        \\
                \hline
            $17$ & $(3,1,0,1,0,2)$ &      \ding{56} &  \ding{56} &  &     &                            \\
                  \hline
            $18$ & $(4,0,0,0,0,0)$ &   1296 & 216 & 24   &  &                            \\
                  \hline
            $19$ & $(4,0,0,0,0,2)$ &       \ding{56} &  \ding{56} &     \ding{56} &    &    \\
                  \hline
            $20$ & $(4,0,0,0,0,4)$ &    \ding{56} &   & &      &                        \\
                  \hline 
\end{tabular}
\caption{Coefficient values $C([\mathbf{e}])$. For a column headed by a particular choice of $(k_1,k_2,k_3,k_4)$, the non-empty entries in that column represent the values of $C([\mathbf{e}])$. If an entry in, say, the $i^{\text{th}}$ row is marked with \ding{56}, then the associated sum $\sps_i$ defined in~\eqref{eqn:sidef} is of negligible size, and so the value of the coefficient is not needed. }
\label{table:coefficient}
\end{table}

Using the coefficients $C([\mathbf{e}])$ from the table yields the equations below. (Note that the first equation has coefficients from the column headed by $(4,4,4,4)$, the second equation has coefficients from the column headed by $(4,4,4,2)$, and so on.)
\footnotesize
\begin{align*}
     \frac{1}{3 ^4} \sum_{j_1,j_2,j_3,j_4=1}^p \E\big( x_{1j_1}^{4} x_{1j_2}^{4} x_{1j_3}^{4}  x_{1j_4}^{4}  \big)
     \ = \ & \frac{r_8}{3 ^4}\big(1296 \sps_{18} +  41472 \sps_{16}  + 165888 \sps_{14} + 31104 \sps_{9} + 248832  \sps_{6} \\&+ 15552 \sps_{5}  + 3888 \sps_{4} + 20736 \sps_{3} +  62208 \sps_{2} +  81\sps_{1}  \big)+ o(p\|R\|_4^8)\\
      \frac{4\cdot2}{3 ^3} \sum_{j_1,j_2,j_3,j_4=1}^p \E\big( x_{1j_1}^{4} x_{1j_2}^{4} x_{1j_3}^{4}  x_{1j_4}^{2}  \big)  \ = \ &\frac{4\cdot2}{3 ^3}r_7\big(216 \sps_{18} +  5184 \sps_{16}  + 13824 \sps_{14} + 5184 \sps_{9} + 31104  \sps_{6} \\&+ 2592 \sps_{5} + 972 \sps_{4}  + 4320 \sps_{3} +  10368 \sps_{2} +  27\sps_{1}  \big) + o(p\|R\|_4^8)\\
      \frac{6\cdot2^2}{3 ^2} \sum_{j_1,j_2,j_3,j_4=1}^p \E\big( x_{1j_1}^{4} x_{1j_2}^{4} x_{1j_3}^{2}  x_{1j_4}^{2}  \big)  \ = \ &\frac{6\cdot2^2}{3 ^2}r_6 \big(24 \sps_{18} +  384 \sps_{16}  + 768 \sps_{14} + 720 \sps_{9} + 3456  \sps_{6} + 432 \sps_{5} \\&+ 234 \sps_{4}  + 864 \sps_{3} +  1728 \sps_{2} +  9 \sps_{1}  \big) + o(p\|R\|_4^8)\\
      \frac{4 \cdot 2^3}{3 } \sum_{j_1,j_2,j_3,j_4=1}^p \E\big( x_{1j_1}^{4} x_{1j_2}^{2} x_{1j_3}^{2}  x_{1j_4}^{2}  \big) 
      \ = \ &\frac{4 \cdot 2^3}{3 }r_5 \big( 72 \sps_{9} + 288  \sps_{6} + 72 \sps_{5} + 54 \sps_{4}  + 168 \sps_{3} +  288 \sps_{2} +  3\sps_{1}  \big) \\&+ o(p\|R\|_4^8)\\
      2^4 \sum_{j_1,j_2,j_3,j_4=1}^p \E\big( x_{1j_1}^{2} x_{1j_2}^{2} x_{1j_3}^{2}  x_{1j_4}^{2}  \big) \ = \ & 2^4 r_4 \big(12 \sps_{5} + 12 \sps_{4}  + 32 \sps_{3} +  48 \sps_{2} +  \sps_{1}  \big)+ o(p\|R\|_4^8).
\end{align*}
\normalsize
Applying these results to (\ref{eq:decomposeof4terms}), it follows from algebraic simplification that
\footnotesize
\begin{align*}
    \sum_{j_1,j_2,j_3,j_4=1}^p   \E(y_{j_1}y_{j_2}y_{j_3}y_{j_4}) \ = \ & 
     16\big(r_8-4r_7+4r_6\big)\sps_{18} + 512 \big(r_8-3r_7+2r_6\big)\sps_{16}+ 2048 \big(r_8-2r_7+r_6\big)\sps_{14} \\&
    + 384 \big(r_8-4r_7+5r_6-2r_5\big)\sps_9 + 3072 \big(r_8-3r_7+3r_6-r_5\big)\sps_6 \\&+ 192 \big(r_8-4r_7+6r_6-4r_5+r_4\big)\sps_5 +48\big(r_8-6r_7+13r_6-12r_5+4r_4\big)  \sps_4\\&
    +256\big(r_8-5r_7+9r_6-7r_5+2r_4\big)  \sps_3  + 768\big(r_8-4r_7+6r_6-4r_5+r_4\big)  \sps_2 \\&+ \big(r_8-8r_7+24r_6-32r_5+16r_4\big)\sps_1+ o(p\|R\|_4^8).
\end{align*}
\normalsize
Lemma \ref{lem:moment} gives $r_{k} = 1 + o(p^{-3/4})$ for $k=1,\ldots,8$ and Lemma \ref{lem:neglibileterms} shows that all of the $\sps_i$ appearing in the equation displayed above, except $\sps_1$ and $\sps_{18}$, 
are $\mathcal{O}(p^{13/8}\|R\|_4^8)$. Also, straightforward calculations show that $\sps_1=p^4$ and $\sps_{18}=p^2\|R\|_4^4$. Accounting for these facts in the equation displayed above leads to~(\ref{eq:4terms}), as needed.\qed

\subsubsection*{Proof of Proposition~\ref{prop:combinationofhighmoments}, equation~\eqref{eq:3terms}.}
We apply the approach developed above. In particular, we start by observing that
\begin{equation}
\label{eq:decomposeof3terms}
\begin{split}
    p\sum_{j_1,j_2,j_3=1}^p\E(y_{j_1}y_{j_2}y_{j_3} ) \ = \ \  \ \ \  \  \frac{1}{3 ^3} &\sum_{j_1,j_2,j_3,j_4=1}^p  \E\big( x_{1j_1}^{4} x_{1j_2}^{4} x_{1j_3}^{4}\big) \\
    -  \frac{3 \cdot 2}{3 ^2} &\sum_{j_1,j_2,j_3,j_4=1}^p \E\big( x_{1j_1}^{4} x_{1j_2}^{4} x_{1j_3}^{2}\big)\\
    +\frac{3\cdot 2^2}{3 } &\sum_{j_1,j_2,j_3,j_4=1}^p \E\big( x_{1j_1}^{4} x_{1j_2}^{2} x_{1j_3}^{2}\big)\\
    - 2 ^3 &\sum_{j_1,j_2,j_3,j_4=1}^p \E\big( x_{1j_1}^{2} x_{1j_2}^{2} x_{1j_3}^{2}\big).
\end{split}
\end{equation}
We need to analyze $\E\big( x_{1j_1}^{k_1} x_{1j_2}^{k_2} x_{1j_3}^{k_2}  x_{1j_4}^{k_4}  \big)$ when the integers $k_1 \geq k_2 \geq k_3$ reside in $\{2,4\}$ and $k_4 = 0$. In this case, the calculation for~\eqref{eq:3terms} may be completed using the information provided in Table~\ref{table:coefficientofthreeterms}, just as Table~\ref{table:coefficient} was used to complete the calculation for~\eqref{eq:4terms}.
\begin{table}[H]
\centering
\begin{tabular}{cc|cccc}
                  &             & \multicolumn{4}{c}{$(k_1,k_2,k_3,k_4)$}               \\ \cline{3-6} 
                  &       $\mathbf{e}$           &  $(4,4,4,0)$  & $(4,4,2,0)$ & $(4,2,2,0)$ & $(2,2,2,0)$\\
                  \hline
           $1$ & $(0,0,0,0,0,0)$ &    27 & 9 & 3 &  1                 \\
                  \hline  
            $3$  & $(1,1,0,1,0,0)$ &    1728 &    288  & 48 & 8\\
              \hline    
            $4$ & $(2,0,0,0,0,0)$ &     648 & 144 & 30& 6\\
             \hline
            $9$ & $(2,2,0,0,0,0)$ &     2592  & 288   & 24&                     \\
                  \hline
            $13$ & $(2,2,0,2,0,0)$ &     \ding{56}  &     & &                       \\
                  \hline
            $16$ & $(3,1,0,1,0,0)$ &    3456 & 192    & &                   \\
                \hline
            $18$ & $(4,0,0,0,0,0)$ &    216 &  24  & &                           \\
                  \hline
\end{tabular}
\caption{Coefficient values $C([\mathbf{e}])$. The table is organized analogously to Table~\ref{table:coefficient}.}
\label{table:coefficientofthreeterms}
\end{table}

Using the coefficients $C([\mathbf{e}])$ from the table then yields 
\begin{align*}
     \frac{1}{3 ^3} \sum_{j_1,j_2,j_3,j_4=1}^p \E\big( x_{1j_1}^{4} x_{1j_2}^{4} x_{1j_3}^{4}   \big)
     \ = \ & \frac{r_6}{3^3}\big(216 \sps_{18} +  3456 \sps_{16}  + 2592 \sps_{9} + 648 \sps_{4} + 1728 \sps_{3} \\&+  27 \sps_{1}  \big)+ o(p\|R\|_4^8)\\
      \frac{3\cdot2}{3 ^2} \sum_{j_1,j_2,j_3,j_4=1}^p \E\big( x_{1j_1}^{4} x_{1j_2}^{4} x_{1j_3}^{2}    \big)  \ = \ &\frac{3\cdot2}{3^2}r_5\big(24 \sps_{18} +  192 \sps_{16}  + 288 \sps_{9} + 144 \sps_{4} + 288 \sps_{3} \\&+  9 \sps_{1}  \big)+ o(p\|R\|_4^8) \\
      \frac{3\cdot2^2}{3 } \sum_{j_1,j_2,j_3,j_4=1}^p \E\big( x_{1j_1}^{4} x_{1j_2}^{2} x_{1j_3}^{2}   \big)  \ = \ &\frac{3\cdot2^2}{3 }r_4 \big( 24 \sps_{9} + 30  \sps_{4} + 48 \sps_{3} +  3 \sps_{1}  \big) + o(p\|R\|_4^8)\\
      2^3 \sum_{j_1,j_2,j_3,j_4=1}^p  \E\big( x_{1j_1}^{2} x_{1j_2}^{2} x_{1j_3}^{2}  \big) \ = \ &2^3 r_3 \big(6  \sps_{4} + 8 \sps_{3} +   \sps_{1}  \big)+ o(p\|R\|_4^8).
\end{align*}
Applying these equations to (\ref{eq:decomposeof3terms}) yields
\begin{align*}
    p \sum_ {j_1,j_2,j_3=1}^p   \E(y_{j_1}y_{j_2}y_{j_3}) \ = \ & 
    8( r_6-2r_5) \sps_{18} +128( r_6-r_5) \sps_{16}+96( r_6-2r_5 +r_4) \sps_9 \\&+24(r_6-4r_5+5r_4-2r_3) \sps_4 +64(r_6-3r_5+3r_4-r_3) \sps_3  \\&+  (r_6-6r_5+12r_4-8r_3)\sps_1+ o(p\|R\|_4^8)\\
     \ = \ & -8p^2\|R\|_{4}^4 + (r_6-6r_5+12r_4-8r_3)p^4 + o(p\|R\|_4^8),
\end{align*}
where the last step is based on $\sps_1=p^4$, $\sps_{18} = p^2\|R\|_{4}^4$, and the formulas for $r_3,\dots,r_6$ and $\sps_3,\dots,\sps_{16}$ given respectively in Lemmas \ref{lem:moment} and \ref{lem:neglibileterms}.
Next, we multiply both sides of the last equation by $2-r_2$. Noting the relation in \eqref{eq:relationofr2} gives $2-r_2=1-\frac{\tau-2}{p} + o(\textstyle \frac{1}{p})$, we have
\begin{align*}
    (2-r_2)p\sum_{j_1,j_2,j_3=1}^p\E(y_{j_1}y_{j_2}y_{j_3} ) \ = \ & - 8 p^2\|R\|_4^4 + (r_6-6r_5+12r_4-8r_3)(p^4 - (\tau-2) p^3 \\&+ o(p^3)) + o(p\|R\|_4^8).
\end{align*}
Finally, the middle term on the right side can be handled by Lemma \ref{lem:moment}, which ensures $r_6-6r_5+12r_4-8r_3 = -1 + o(p^{-3/4})$ and leads to~\eqref{eq:3terms}.\qed

\subsubsection*{Proof of Proposition~\ref{prop:combinationofhighmoments}, equation~\eqref{eq:2terms}.}
Observe that the definition of $y_{j}$ in equation \eqref{eqn:yjdef} gives
\small
\begin{equation*}
% \label{eq:decomposeof2terms}
\begin{split}
    \sum_{j_1,j_2 =1}^p\E(y_{j_1}y_{j_2} ) \ = \  \frac{1}{3^2}\sum_{j_1,j_2=1}^p  \E\big( x_{1j_1}^{4} x_{1j_2}^{4}\big)-\frac{2\cdot2}{3 } \sum_{j_1,j_2=1}^p  \E\big( x_{1j_1}^{4} x_{1j_2}^{2}\big) + 2^2 \sum_{j_1,j_2=1}^p  \E\big( x_{1j_1}^{2} x_{1j_2}^{2}\big).
\end{split}
\end{equation*}
\normalsize
The proof of Lemma \ref{lem:covariances} shows that when the diagonal entries of $\SIGMA$ are equal to 1, the following equations hold
\begin{align*}
     \frac{1}{3 ^2} \sum_{j_1,j_2=1}^p  \E\big( x_{1j_1}^{4} x_{1j_2}^{4}\big)  & \ = \  r_4 \Big(p^2 + 8 \|R\|_2^2+\ts \frac{8}{3}\|R\|_4^4 \Big)\\
     \frac{2\cdot2}{3 } \sum_{j_1,j_2=1}^p  \E\big( x_{1j_1}^{4} x_{1j_2}^{2}\big)  & \ = \  r_3 \Big(4p^2 +  16 \|R\|_2^2 \Big) \\
     2^2 \sum_{j_1,j_2=1}^p  \E\big( x_{1j_1}^{2} x_{1j_2}^{2}\big)  & \ = \  r_2 \Big(4p^2 +  8 \|R\|_2^2 \Big) .
\end{align*}
Combining these calculations, we have
\begin{align*}
    p^2 \sum_{j_1,j_2=1}^p \E(y_{j_1}y_{j_2})  
     \ = \ & (r_4-4r_3 +4r_2) p^4 + 8(r_4- 2r_3 + r_2) p^2\|R\|_2^2 +\frac{8r_4}{3}p^2\|R\|_4^4\\
     \ = \ &(r_4-4r_3 +4r_2) p^4 +\frac{8}{3}p^2\|R\|_4^4 + o(p\|R\|_4^8 ),
\end{align*}
where the last step uses Lemma \ref{lem:moment} and the fact that $\|R\|_2^2 \lesssim p^{1/2}\|R\|_4^4$.
Combining the relation $r_4-4r_3 +4r_2 =1 + o(p^{-3/4})$ given by Lemma \ref{lem:moment} with $(2-r_2)^2 =1-\frac{2(\tau-2)}{p} + o(\textstyle \frac{1}{p})$ implied by Assumption \ref{Data generating model}, we have
\small
\begin{align*}
    (2-r_2)^2 p^2\sum_{j_1,j_2=1}^p\E(y_{j_1}y_{j_2}) 
     \ = \ &(r_4-4r_3 +4r_2)( p^4- 2(\tau-2) p^3 + o(p^3)) +\frac{8}{3}p^2\|R\|_4^4 + o(p\|R\|_4^8 )\\
     \ = \ &p^4 (r_4-4r_3 +4r_2)-2(\tau -2)p^3 + \frac{8}{3} p^2\|R\|_4^4+ o(p\|R\|_4^8),
\end{align*}
\normalsize
which proves~\eqref{eq:2terms}.
\qed

\begin{lemma}
\label{lem:neglibileterms} 
    The sums $\sps_i$ defined in \eqref{eqn:sidef} satisfy
        \begin{subnumcases}{\sps_i=}o(p\|R\|_4^8) \textup{ \ \ \ \ \ \ \ \ \  for \ \ \  } i=7, 8, 10, 11, 12, 13, 15, 17, 19, 20\label{eqn:case1}\\[0.2cm]
        \mathcal{O}(p^{13/8}\|R\|_4^8) \textup{ \ \ \  for \ \ \ } i=2,3, 4, 5, 6, 9, 14, 16.\label{eqn:case2}
        \end{subnumcases}
\end{lemma}
\proof We begin by analyzing the sums $\sps_i$ in the first case~\eqref{eqn:case1}. For these calculations, we will often use the fact that $\frac{\|R\|_2}{\|R\|_4^2} = \mathcal{O}(p^{1/4})$, which holds for a generic correlation matrix $R$.
\begin{itemize} 
    \item $i = 19,$ 20: The sums $\sps_i$ are bounded by 
     $$\sum_{j_1,j_2,j_3,j_4=1}^p R_{j_1,j_2}^4R_{j_3,j_4}^2 \ = \ \|R\|_4^4\|R\|_2^2 \ = \ o(p\|R\|_4^8).$$
    \item $i = 17$: Due to the basic inequality $|ab|\leq \frac{1}{2}(a^2+b^2)$, the absolute value $|\sps_i|$ is bounded by
    \begin{equation}
    \label{eq:lin1}
    \begin{split}
        \sum_{j_1,j_2,j_3,j_4=1}^p \Big|R_{j_1,j_2}^3R_{j_1,j_3}R_{j_2,j_3}R_{j_3,j_4}^2 \Big|& \ \leq \  \sum_{j_1,j_2,j_3,j_4=1}^p \ts \frac{1}{2}R_{j_1,j_2}^2R_{j_3,j_4}^2\big(R_{j_1,j_3}^2 +R_{j_2,j_3}^2\big)\\
        & \ = \  \mathbf{1}\ttop R^{\circ 2}R^{\circ 2}R^{\circ 2}  \mathbf{1}\\
        & \ \leq \   p \tr( R^{\circ 2}R^{\circ 2}R^{\circ 2})\\
        &  \ \leq \  p \tr( R^{\circ 2}R^{\circ 2})^{3/2}\\
        &  \ = \  p (\|R\|_4^4)^{3/2}\\
        &  \ = \  o(p\|R\|_4^8),
    \end{split}
    \end{equation}
    where $\mathbf{1}$ denotes the $p$-dimensional vector of ones, and the third inequality is based on the fact that the Schur product theorem guarantees that $R^{\circ k}$ is positive semidefinite for any integer $k\geq 1$.
    \item $i = 15$: The absolute value $|\sps_i|$ is equal to 
    \begin{align*}
        \Bigg|\sum_{j_1,j_2,j_3,j_4=1}^p R_{j_1,j_2}^3R_{j_1,j_3}R_{j_2,j_4}R_{j_3,j_4}^3\Bigg|
        & \ = \  |\tr(R^{\circ 3}RR^{\circ 3}R)|\\[-.3cm]
& \ \leq \ \|R^{\circ 3} R\|_2^2\\
& \ \leq \ \|R\|_2^2\|R^{\circ 3}\|_2^2\\
        &  \ \leq \  \|R\|_2^2 \|R\|_4^4\\
        &  \ = \  o(p\|R\|_4^8),
    \end{align*}
    where the step going from the third to the fourth line uses the fact that $R_{jk}^6\leq R_{jk}^4$ for all $j,k$.
    \item $i = 13$: The sum $\sps_i$ is equal to
    \begin{align*}
        \sum_{j_1,j_2,j_3,j_4=1}^p R_{j_1,j_2}^2R_{j_1,j_3}^2R_{j_2,j_3}^2  \ = \  &p\,\tr(R^{\circ 2}R^{\circ 2}R^{\circ 2}) 
        \ = \  o(p\|R\|_4^8),
    \end{align*}
    where the last step uses the calculation in (\ref{eq:lin1}).
    \label{item:(2,2,0,2,0,0)}
    \item $i = 12$: The absolute value $|\sps_i|$ is equal to  
     \begin{align*}
        \Bigg|\sum_{j_1,j_2,j_3,j_4=1}^p R_{j_1,j_2}^2R_{j_1,j_3}^2R_{j_2,j_3}R_{j_2,j_4}R_{j_3,j_4}\Bigg| 
        &  \ = \  \big|\tr(R^{\circ 2}( R \circ R^2) R^{\circ 2})\big|\\
        & \ \leq \  \tr(R^{\circ 2} R^{\circ 2})\tr( R \circ R^2)\\
        &  \ \leq \ \|R\|_4^4 \|R\|_2^2      \ \ \  \\
        %Theorem 5.5.12 in \cite{roger1994topics}} \\
        &  \ = \  o(p\|R\|_4^8),
    \end{align*}
    where the second inequality is based on the fact that the Schur product theorem guarantees that $ R \circ R^2$ is positive semidefinite, and the second-to-last step uses $\tr(R\circ R^2) = \tr(R^2)$, which holds since $R$ is a correlation matrix.
    \item $i = 11$: The sum $\sps_i$ is equal to
    \begin{align*}
        \sum_{j_1,j_2,j_3,j_4=1}^p R_{j_1,j_2}^2R_{j_1,j_3}^2 R_{j_2,j_4}^2R_{j_3,j_4}^2 
        &  \ = \  \tr( R^{\circ 2}R^{\circ 2}R^{\circ 2}R^{\circ 2})\\[-.3cm]
        &  \ \leq \  \tr( R^{\circ 2}R^{\circ 2})^2 \\
        &  \ = \  \|R\|_4^8 \\
        &  \ = \  o(p\|R\|_4^8).
    \end{align*}
    \item 
    $i = 10$: The calculation in (\ref{eq:lin1}) implies that the sum $\sps_i$ is equal to
    \begin{align*}
        \sum_{j_1,j_2,j_3,j_4=1}^p R_{j_1,j_2}^2R_{j_1,j_3}^2 R_{j_2,j_4}^2 & \ = \  \mathbf{1}\ttop R^{\circ 2}R^{\circ 2}R^{\circ 2}  \mathbf{1}    \ = \  o(p\|R\|_4^8).
    \end{align*}
    \item $i = 7$, 8: Using the basic inequality $|ab|\leq \frac{1}{2}(a^2+b^2)$ twice, the absolute values $|\sps_i|$ are bounded by 
    \footnotesize
     \begin{align*}
     \sum_{j_1,j_2,j_3,j_4=1}^p \Big|R_{j_1,j_2}^2R_{j_1,j_3}R_{j_1,j_4}R_{j_2,j_3}R_{j_2,j_4}\Big|  \ \leq \  & \sum_{j_1,j_2,j_3,j_4=1}^p \ts\frac{1}{4}R_{j_1,j_2}^2 (R_{j_1,j_3}^2 + R_{j_1,j_4}^2) (R_{j_2,j_3}^2 + R_{j_2,j_4}^2)\\
          \ \lesssim \ & \sum_{j_1,j_2,j_3,j_4=1}^p R_{j_1,j_2}^2R_{j_1,j_3}^2R_{j_2,j_3}^2 \\&+ \sum_{j_1,j_2,j_3,j_4=1}^p R_{j_1,j_2}^2R_{j_1,j_3}^2 R_{j_2,j_4}^2 \\[0.1cm]
          \ = \  &o(p\|R\|_4^8),
    \end{align*}
    \normalsize
    where the last step uses the previous calculations for $\sps_{10}$ and $\sps_{13}$.
\end{itemize}
~\\
In the remainder of the proof, we analyze the sums $\sps_i$ in the second case~\eqref{eqn:case2}. 
\begin{itemize}
    \item $i = 2$: The sum $\sps_i$ is equal to 
    \begin{align*}
        \sum_{j_1,j_2,j_3,j_4=1}^p R_{j_1,j_2}R_{j_1,j_3}R_{j_2,j_4}R_{j_3,j_4}& \ = \  \tr(R^4)\\[-0.4cm] 
& \ \leq \  \|R\|_2^4 \\ 
& \ = \  \mathcal{O}(p \|R\|_4^8).
    \end{align*}
    \item $i = 3$: The sum $\sps_i$ is equal to 
    \begin{align*}
        \sum_{j_1, j_2, j_3, j_4 = 1}^p R_{j_1,j_2} R_{j_1,j_3} R_{j_2,j_3} & \ = \  p\tr(R^3) \\[-0.4cm]
        & \ \leq \  p \|R\|_2^3\\ 
        & \ = \  \mathcal{O}(p^{7/4}\|R\|_4^6)\\ 
        & \ = \  \mathcal{O}(p^{5/4}\|R\|_4^8),
    \end{align*}
    where the last step uses the inequality $\|R\|_4^2 \geq p^{1/2}$ that holds for a generic correlation matrix $R$.
    \item $i = 4$, 5, 6, 9, 14, 16: The absolute value $|\sps_i|$ is bounded by 
    \begin{align*}
        \sum_{j_1,j_2,j_3,j_4=1}^p R_{j_1,j_2}^2  \ = \ \, & p^2\|R\|_2^2 \ = \   \mathcal{O}(p^{3/2}\|R\|_4^8).
    \end{align*}
\end{itemize}
\qed

\begin{lemma}
\label{lem:moment}
    If Assumption \ref{Data generating model} holds, then for
    %r_1=1$, $r_2=1+\frac{\tau-2}{p}+o(\frac{1}{p})$, and for 
    $k = 1,\ldots,8$,
    $$r_{k} = 1 + o(p^{-3/4}).$$ 
    In addition, the following relations hold
    \begin{align}
         r_4  - 4 r_3 -  r_2^2+8 r_2-4 & \ = \ \sum_{l=0}^4 (-1)^{l}\textstyle\binom{4}{l}  r_{l-1} +o(\textstyle \frac{1}{p}) \label{eqn:firstrk}\\
         \sum_{l=0}^4 (-1)^{l}\textstyle\binom{4}{l}  r_{l} & \ = \  o(\textstyle \frac{1}{p})\label{eqn:secondrk}\\
        \sum_{l=0}^8 (-1)^{l}\textstyle\binom{8}{l}  r_{l} & \ = \  o(\textstyle \frac{1}{p}).\label{eqn:thirdrk}
    \end{align} 
\end{lemma}

\proof
  Assumption 1 implies $\E\big((\frac{\xi_1^2}{p}-1)^{k}\big) = o(p^{-3/4})$ for $k=1,\ldots,8$, and so
  \begin{align*}
      \frac{ \E(\xi_1^{2k}) }{p^k} \ = \ & \E\,\Big(\big(\ts\frac{\xi_1^2}{p}-1 + 1\big)^{k}\Big) \\ \ = \ &  1 + o(p^{-3/4}).
  \end{align*}
 Also, the definition of $r_k$ in \eqref{eqn:extrark} implies
\begin{align}
    \begin{split}
        r_k& \ = \  \frac{\E(\xi_1^{2k})}{p^{k}} -\frac{\E(\xi_{1}^{2k})\Big(k(k-1)p^{k-2}+\mathcal{O}(p^{k-3})\Big)}{p^{k-1}\E(\|\z_1\|_2^{2k})}\\
    & \ = \ \frac{\E(\xi_1^{2k})}{p^{k}}  -\frac{ k(k-1) }{p} + o(\textstyle \frac{1}{p})\label{eq:decomposerk}\\
    &  \ = \   1 + o(p^{-3/4}).
    \end{split}
\end{align}
    Towards proving \eqref{eqn:firstrk}, note that the relation $r_2 = 1+ \frac{\tau - 2}{p} + o(\textstyle \frac{1}{p})$ in \eqref{eq:relationofr2} implies
    \begin{align*}
        r_2^2 & \ = \  2r_2 - 1+ o(\textstyle \frac{1}{p}).
    \end{align*}
    Combining $r_0=r_1 = 1$ and the condition $\E\Big(\big| \ts \frac{\xi_1^{2}}{p} -1\big|^{4}\Big)= o(\textstyle \frac{1}{p})$ implied by Assumption \ref{Data generating model} with the second equality in \eqref{eq:decomposerk}, we have
    \begin{align*}
        r_4 - 4r_3- r_2^2+8r_2 -4 & \ = \ r_4- 4r_3+6r_2 -4r_1+r_0 + o(\textstyle \frac{1}{p})\\ 
        & \ = \ \sum_{l=0}^4 (-1)^{l}\binom{4}{l} \frac{\E(\xi_1^{2l})}{p^l} +o(\textstyle \frac{1}{p})\\
        &  \ = \  \E\Big(\big| \ts \frac{\xi_1^{2}}{p} -1\big|^{4}\Big)+ o(\textstyle \frac{1}{p}) \\& \ = \  o(\textstyle \frac{1}{p}),
    \end{align*}
    which proves \eqref{eqn:firstrk} and \eqref{eqn:secondrk}.
    Lastly, to show~\eqref{eqn:thirdrk}, observe that the second equality in \eqref{eq:decomposerk} implies
    \begin{align*}
        \sum_{l=0}^8 (-1)^{l}\binom{8}{l} r_{l}& \ = \  \sum_{l=0}^8 (-1)^{l}\binom{8}{l} \frac{\E(\xi_1^{2l})}{p^l} +o(\textstyle \frac{1}{p})\\
     &  \ = \  \E\Big(\big| \ts \frac{\xi_1^2}{p}-1\big|^8
     \Big) +o(\textstyle \frac{1}{p})\\&  \ = \  o(\textstyle \frac{1}{p}),
    \end{align*}
    where the last step follows from Assumption~\ref{Data generating model}.
\qed

\subsection{Asymptotic normality of $T_{n,2}$}\label{app:T2}

\begin{proposition}
\label{prop:T2smallr}
    Suppose Assumption \ref{Data generating model} holds. If ${\tt{r}}(\SIGMA) \lesssim \sqrt{p}$, then as $n \to \infty$,
    \begin{align*}
        \frac{T_{n,2}}{\sigma_{n,2}}  \ \Rightarrow \ N(0,1).
    \end{align*}
\end{proposition}
\proof
Define the statistics
\begin{align*}
    U &\ = \ \check{\varsigma}^2 - 2(\check{\nu}_2 - \ts\frac{2}{n} \check{\nu}_1^2)\\[0.2cm]
    V &\ = \ \check{\nu}_1^2 + 2(\check{\nu}_2 - \ts\frac{2}{n} \check{\nu}_1^2),
\end{align*}
which allow $T_{n,2}$ to be rewritten as
\begin{equation} \label{eq:definitionofT2useU}
    T_{n,2} \ = \  \sqrt{\frac{pn}{2}} \Big( \frac{U}{V} - \frac{  {\varsigma}^2 - 2 {\nu}_2  }{ {\nu}_1^2 + 2 {\nu}_2  }\Big).
\end{equation}
When ${\tt{r}}(\SIGMA) \lesssim \sqrt{p}$, it will be shown in Lemma~\ref{lem:UCLT} that the statistic $U$ satisfies 
\begin{align}\label{eq:consistencyofU}
    \frac{\sqrt{np/2}(U - ({\varsigma}^2 - 2 {\nu}_2 ))}{\sigma_{n,2}(\nu_1^2 + 2\nu_2)}  
    &\ \Rightarrow \ \ N(0,1)
\end{align}
as $n\to\infty$.
Next, we establish a limit for the statistic $V$. The proofs of Lemmas A.1 and A.2 in \citep{wang2022bootstrap} show that if Assumption \ref{Data generating model} holds, then as $n\to\infty$
\begin{align}\label{eq:consistencyofsigma1}
    \frac{\check\nu_2-\ts \frac{2}{n} \check\nu_1^2}{\nu_2} \  \xrightarrow{\P} \ 1 \ \ \ \ \ \text{ and } \ \ \ \ \
    \frac{\check{\nu}_1^2}{\nu_1^2}\ \xrightarrow{\P} \ 1 ,
\end{align}
and hence
\begin{equation} \label{eq:asymofV}
    \frac{V}{\nu_1^2+ 2\nu_2} \ \xrightarrow{\P} \ 1.
\end{equation}
Now we will combine the previous limits for $U$ and $V$ to reach the stated result. Based on the representation of $T_{n,2}$ given in \eqref{eq:definitionofT2useU}, the limit for $V$ implies
\begin{equation*}
    \frac{T_{n,2}}{\sigma_{n,2}} \ = \  \sqrt{\frac{pn}{2}} \Big( \frac{U}{\sigma_{n,2}(\nu_1^2+ 2\nu_2)(1+o_{\P}(1))}- \frac{  {\varsigma}^2 - 2 {\nu}_2  }{ \sigma_{n,2}({\nu}_1^2 + 2 {\nu}_2)  }\Big)  .
\end{equation*}
Due to the limit for $U$ given in \eqref{eq:consistencyofU}, the desired limit for $T_{n,2}$ will follow from Slutsky's lemma if we can show that 
\begin{equation*}
% \label{eqn:UbigOP}
\frac{\sqrt{pn}\,U}{\sigma_{n,2}(\nu_1^2+2\nu_2)}\ = \ \mathcal{O}_{\P}(1).
\end{equation*}
Noting that the definition of $\sigma_{n,2}^2$ in \eqref{eq:definitionofsigma1} gives
\begin{align}\label{eq:transofsigma2}
    \sqrt{16\nu_4+8\nu_2^2} \ = \  \frac{\sigma_{n,2}(\nu_1^2+2\nu_2)}{\sqrt{p}},
\end{align}
it follows from \eqref{eq:consistencyofU} that
\begin{align*}
    \frac{\sqrt{pn}\,U}{\sigma_{n,2}(\nu_1^2+2\nu_2)} \ = \  \frac{\sqrt{n}(\varsigma^2-2\nu_2)}{\sqrt{16\nu_4+8\nu_2^2}} \ + \ \mathcal{O}_{\P}(1).
\end{align*}
So, to complete the proof, it is enough to show $\varsigma^2-2\nu_2 = \mathcal{O}(\ts \frac{\nu_2}{\sqrt{p}})$.
To see this, observe that when ${\tt{r}}(\SIGMA) \lesssim \sqrt{p}$, Lemma \ref{lem:quadform} and equation~\eqref{eq:relationofr2} imply
\begin{align*}
% \label{eq:varsigma^2-2nu_2}
    \varsigma^2-2\nu_2 \ = \ &\mathcal{O}(\ts \frac{1}{p}) (\nu_1^2 + 2\nu_2)\\
    \ = \ &\mathcal{O}(\ts \frac{\nu_2}{\sqrt{p}}),
\end{align*}
as needed. \qed

\begin{lemma}\label{lem:UCLT}
    Suppose Assumption \ref{Data generating model} holds. If ${\tt{r}}(\SIGMA) \lesssim \sqrt{p}$, then as $n \to \infty$,
    \begin{align*}
        % \label{eq:consistencyofUagain}
        \frac{\sqrt{np/2}(U - ({\varsigma}^2 - 2 {\nu}_2 ))}{\sigma_{n,2}(\nu_1^2 + 2\nu_2)}  
        &\ \Rightarrow \ \ N(0,1).
    \end{align*}
\end{lemma}
\proof To show that $U$ has the structure of a U-statistic, define the kernel
\begin{align}\label{eqn:definitionofh}
    h(\x_1,\x_2) \ = \  \Big(\ts \frac{1}{2}-\frac{2(n-2)}{n^2}\Big)  \Big(\|\mathbf{x}_1\|_2^2 - \|\mathbf{x}_2\|_2^2\Big)^2 - \frac{2(n-2)}{n} (\mathbf{x}_1\ttop \mathbf{x}_2)^2,
\end{align}
which yields
\begin{align*}
% \label{eqn:taudecomposition}
     U \ = \ & \frac{1}{\binom{n/2}{2}} \sum_{n/2 < i < j\leq n}   h(\x_i,\x_j).
\end{align*}
Furthermore, by letting
\begin{equation}\label{eq:definitionofh1}
\begin{split}
\theta &\ = \  \E(h(\x_1,\x_2))\\
    h_1(\x_1) &\ = \  \E(h(\x_1,\x_2) |\x_1) - \theta,
    \end{split}
\end{equation}
the H\'ajek projection of $U-\theta$ can be expressed as 
\begin{equation*}
    \hat U 
    \ = \  \frac{2}{n/2}\sum_{i=n/2+1 }^n h_1(\x_i).
\end{equation*}
It is a classical fact~\citep[Theorem 11.2]{van2000asymptotic} that if 
\begin{equation}\label{eqn:varcond}
    \frac{\var(\hat U)}{\var(U)} \ \to \ 1
\end{equation}
as $n\to\infty$, then
\begin{equation}\label{eqn:CLTequiv}
     \frac{U-\theta}{\sqrt{\var(U)}}  \ = \    \frac{\hat U}{\sqrt{\var(\hat U)}}+ o_{\P}(1).
\end{equation}
This shows that if $U$ and $\hat U$ are suitably normalized, then the asymptotic normality of one implies that of the other. To check the condition~\eqref{eqn:varcond}, we first introduce the parameters
\begin{align*}
    \psi_1 &  \ = \ \cov(h(\x_1,\x_2), h(\x_1,\x_3)) \\
\psi_2 &  \ = \  \var(h(\x_1,\x_2)).
\end{align*}
It is known from the theory of U-statistics~\cite[][Theorem 12.3]{van2000asymptotic} that the variances of $\hat U$ and $U$ have exact expressions in terms of $\psi_1$ and $\psi_2$, namely
\begin{align}
    \var(\hat{U}) \ = \ & \frac{4}{n/2} \psi_1\label{eqn:hatUvar}\\
    \var(U) \ = \ & \frac{2(n/2-2)}{\binom{n/2}{2}} \psi_1 + \frac{1}{\binom{n/2}{2}} \psi_2.\label{eq:varu}
\end{align}
Therefore, the condition~\eqref{eqn:varcond} will hold, if we can show
\begin{equation}\label{eq:varandcovofu}
    \psi_2 \ = \ o(n\,\psi_1).
\end{equation}
The sizes of $\psi_1$ and $\psi_2$ can be compared with Lemmas \ref{lem:psi2} and \ref{lem:psi1}, which show that under the condition ${\tt{r}}(\SIGMA) \lesssim \sqrt{p}$ we have
\begin{align}
    \psi_2 \  \lesssim \ &\tr(\SIGMA^2)^2\label{eq:covarianceofT2smallrankfirst}\\ 
    \psi_1 
    \ = \ &4 \tr(\SIGMA^4) + \ts 2\tr(\SIGMA^2) ^2 + o(\tr(\SIGMA^2) ^2)\label{eq:covarianceofT2smallrank}.
\end{align}
From these relations, it is straightforward to check that~\eqref{eq:varandcovofu} holds, and hence~\eqref{eqn:CLTequiv} holds as well. 

Having shown~\eqref{eqn:CLTequiv}, the proof will be complete if we can establish two more items, which are 
\begin{equation}\label{eqn:hatUCLT}
 \frac{\hat U}{\sqrt{\var(\hat U)}} \ \Rightarrow \ N(0,1) 
\end{equation}
and
\begin{equation}\label{eqn:betterUCLT}
     \frac{U-\theta}{\sqrt{\var(U)}} \ = \    \frac{\sqrt{np/2}(U - ({\varsigma}^2 - 2 {\nu}_2 ))}{\sigma_{n,2}(\nu_1^2 + 2\nu_2)}  + o_{\P}(1).
\end{equation}
Since $\hat U$ is a sum of centered independent random variables, the limit~\eqref{eqn:hatUCLT} can be established by checking the following Lyapunov-type condition, which implies the Lindeberg central limit theorem in the setting of triangular arrays~\cite[Proposition 2.27]{van2000asymptotic} (and so is applicable in the high-dimensional setting),
\begin{equation} \label{eq:lindebergofT2}
    \frac{\E(h_1(\x_1)^{4})}{n \var(h_1(\x_1))^2}  \ \to \ 0.
\end{equation}
This will be established in Lemma~\ref{lem:4thmomentofh1}. 

With regard to~\eqref{eqn:betterUCLT}, it follows from equations~\eqref{eqn:varcond},~\eqref{eqn:hatUvar}, and \eqref{eq:covarianceofT2smallrank} that the standard deviation of $U$ satisfies
\begin{equation*}
% \label{eqn:firstvar}
    \frac{\sqrt{\var(U)}}{\sqrt{\frac{1}{n/2}(16 \nu_4+  8\nu_2 ^2)}} \ \to \ 1.
\end{equation*}
Consequently, equation~\eqref{eq:transofsigma2} gives
\begin{equation}\label{eqn:bettersd}
    \frac{\sqrt{\var(U)}}{\sqrt{\frac{2}{np}}\sigma_{n,2}(\nu_1^2+2\nu_2)} \ \to \ 1.
\end{equation}
Finally, to consider the the mean parameter $\theta$, it can be checked using Lemma \ref{lem:quadform} as well as equation \eqref{eq:relationofr2} that
\begin{equation} \label{eqn:taumean}
    \begin{split}
    \theta & \ = \ \Big(1- \ts \frac{4(n-2)}{n^2}\Big) \varsigma^2 - \frac{2(n-2)}{n}\nu_2\\
    & \ = \ \varsigma^2 - 2 \nu_2 +  \mathcal{O}(\ts \frac{\nu_2}{n}).
    \end{split}
\end{equation}
Combining this with~\eqref{eqn:bettersd} implies~\eqref{eqn:betterUCLT}, which completes the proof.\qed

\begin{lemma}\label{lem:4thmomentofh1}
    Suppose Assumption \ref{Data generating model} holds. If ${\tt{r}}(\SIGMA) \lesssim \sqrt{p}$, then $h_1(\x_1)$ defined in \eqref{eq:definitionofh1} satisfies
     \begin{equation} 
     \E(h_1(\x_1)^{4})  \ = \  o(n \var(h_1(\x_1))^2).
 \end{equation} 
\end{lemma}
\proof
From the definition of $h_1(\x_1)$, we have
\footnotesize
\begin{align*}
    h_1(\x_1) \ = \  \Big(\ts \frac{1}{2}-\frac{2(n-2)}{n^2}\Big)  \Big(\big(\mathbf{x}_1\ttop \mathbf{x}_1- \tr(\SIGMA)\big)^2+ (1-r_2) \tr(\SIGMA)^2 -2 r_2\tr(\SIGMA^2)\Big) - \frac{2(n-2)}{n} \Big(\mathbf{x}_1\ttop \SIGMA \mathbf{x}_1-\tr(\SIGMA^2) \Big). 
\end{align*}
\normalsize
The variance of $h_1(\x_1)$ can be computed via equations~\eqref{eqn:hatUvar} and~\eqref{eq:covarianceofT2smallrank} to obtain
\begin{equation*}
    \var(h_1(\x_1)) = 4\tr(\SIGMA^4) + 2\tr(\SIGMA^2) ^2 + o(\tr(\SIGMA^2)^2).
\end{equation*}
So, to complete the proof, it is sufficient to show that
\begin{equation*}
 \E(h_1(\x_1)^4) \ = \ o(n\,\tr(\SIGMA^2)^4).   
\end{equation*}
Due to equation \eqref{eq:relationofr2}, we have $|1-r_2| \tr(\SIGMA)^2\lesssim \tr(\SIGMA^2)$, and so 
\begin{equation}
    \label{eq:linoftauhat}
\begin{split}
    \E(h_1(\x_1)^{4}) \  &\lesssim \  \|\mathbf{x}_1\ttop \mathbf{x}_1- \tr(\SIGMA)\|_{L^8}^ {8} \ + \   \tr(\SIGMA^2) ^ {4} \ + \ \big\|\mathbf{x}_1\ttop \SIGMA \mathbf{x}_1-\tr(\SIGMA^2) \big\|_{L^4}^ {4}.
\end{split}
\end{equation}
Lemma \ref{lem:8thquadformofu} implies
\begin{align*}
    \big\|\mathbf{x}_1\ttop \SIGMA \mathbf{x}_1-\tr(\SIGMA^2) \big\|_{L^4}^ {4}
    & \ \lesssim \   \tr(\SIGMA ^4)^2 +  o(\ts \frac{1}{p}) \tr(\SIGMA^2) ^4\\
    & \ = \ o(n\tr(\SIGMA^2)^4).
\end{align*}
In addition, when ${\tt{r}}(\SIGMA) \lesssim \sqrt{p}$, Lemma \ref{lem:8thquadformofu} gives
\begin{align*}
    \|\mathbf{x}_1\ttop \mathbf{x}_1- \tr(\SIGMA)\|_{L^8}^ {8}  \ &= \ o(n \tr(\SIGMA^2) ^4).
\end{align*}
Applying the above results to (\ref{eq:linoftauhat}) completes the proof.
\qed

\begin{proposition}
    \label{prop:T2larger}
    Suppose Assumption \ref{Data generating model} holds. If ${\tt{r}}(\SIGMA)/\sqrt p\to\infty$ as $n \to \infty$, then 
    \begin{align*}
        \frac{T_{n,2}}{\sigma_n}  \ \xrightarrow{\P} \ 0.
    \end{align*}
\end{proposition}
\proof
The definition of $T_{n,2}$ in \eqref{eq:definitionofT2useU} gives
\begin{equation}\label{eqn:beginT2}
       \frac{T_{n,2}}{\sigma_n} \ = \ \frac{\sqrt{np/2}U}{\sigma_n(\nu_1^2 + 2\nu_2)} \Big( \frac{\nu_1^2 + 2\nu_2}{V}-1 \Big) + \frac{\sqrt{np/2}(U - ({\varsigma}^2 - 2 {\nu}_2 ))}{\sigma_n(\nu_1^2 + 2\nu_2)}.
\end{equation}
To address the second term, due to the fact $\sigma_n^2 \geq \sigma_{n,1}^2 \gtrsim 1$, we will show
\begin{equation}\label{eq:asymptoticofU}
    \frac{\sqrt{np/2}}{\nu_1^2 + 2\nu_2}\big(U - ({\varsigma}^2 - 2 {\nu}_2)  \big)   \ \xrightarrow{\P} \ 0
\end{equation}
by showing that the mean and variance of the left side are both $o(1)$. When ${\tt{r}}(\SIGMA)/\sqrt p\to\infty$, 
equation (\ref{eqn:taumean}) gives
 \begin{align*}
    \frac{p}{\nu_1^2 + 2\nu_2} \E\big(U - ({\varsigma}^2 - 2 {\nu}_2)  \big) 
    &\ = \  \mathcal{O}\big( \ts \frac{1}{{\tt{r}}(\SIGMA)}\big)\\&\ = \ o(1).
\end{align*}
Regarding the variance of $U$, we will show $\var(\frac{p}{\nu_1^2+2\nu_2} U) = o(1)$.
%$o(\ts \frac{\nu_1^4}{p^2})$. 
Due to equation \eqref{eq:varu}, it is sufficient to show $\psi_1 = o(\ts \frac{\nu_1^4}{p})$ and $\psi_2= o(\nu_1^4)$, and these relations will be established in Lemmas \ref{lem:psi2} and \ref{lem:psi1}.
Consequently, the limit \eqref{eq:asymptoticofU} holds.

Next, to consider the first term in~\eqref{eqn:beginT2}, it is known from \eqref{eq:asymofV} that $\frac{V}{\nu_1^2+ 2\nu_2} \xrightarrow{\P} 1$. Also, it is known from Lemma \ref{lem:quadform},  \eqref{eq:relationofr2} and \eqref{eq:asymptoticofU} that
\begin{equation*}
    \frac{\sqrt{np/2}U}{\sigma_n(\nu_1^2 + 2\nu_2)} \ = \ \mathcal{O}_{\P}(1),
\end{equation*}
which implies that the first term in~\eqref{eqn:beginT2} is $o_{\P}(1)$ and completes the proof.
\qed

\begin{lemma}
\label{lem:psi2}
Suppose that Assumption \ref{Data generating model} holds. If ${\tt{r}}(\SIGMA)/\sqrt p\to\infty$ as $n\to\infty$, then
\begin{align*}
    \psi_2 \ = \  o(\tr(\SIGMA)^4).
\end{align*}
Alternatively, if ${\tt{r}}(\SIGMA) \lesssim \sqrt{p}$, then 
\begin{align*}
   \psi_2 \ \lesssim \ &\tr(\SIGMA^2) ^2.
\end{align*}
\end{lemma}
\proof
Using the definition of $h(\x_1,\x_2)$ in \eqref{eqn:definitionofh}, as well as the bound $\var(V_1+V_2)\leq 2\var(V_1)+2\var(V_2)$ for generic random variables $V_1$ and $V_2$, we have
\begin{equation}
\label{eq:varianceofT2}
\begin{split}
    \psi_2 
    \  \lesssim \  & \var\big((\mathbf{x}_1\ttop \mathbf{x}_1 - \mathbf{x}_2\ttop \mathbf{x}_2)^2\big) \ +  \  \var\big((\mathbf{x}_1\ttop \mathbf{x}_2)^2 \big).
\end{split}
\end{equation}
The second term on the right side can be controlled with Lemma D.2 in \cite{wang2022bootstrap} and equation~\eqref{eq:relationofr2}, which gives
\begin{align*}
    \var\big((\mathbf{x}_1\ttop \mathbf{x}_2)^2 \big)
    & \ = \  3r_2^2 (\tr(\SIGMA^2)^2 + 2 \tr(\SIGMA^4)) - \tr(\SIGMA^2)^2 \\
    & \ \lesssim \ \tr(\SIGMA^2)^2  .
\end{align*}
For the first term on the right side of \eqref{eq:varianceofT2}, we have
\begin{equation*}
    \var\big((\mathbf{x}_1\ttop \mathbf{x}_1 - \mathbf{x}_2\ttop \mathbf{x}_2)^2\big) \ \lesssim \  \E\big((\mathbf{x}_1\ttop \mathbf{x}_1 - \tr(\SIGMA))^4\big).
\end{equation*}
By expanding out the fourth moment directly, and using  Lemmas \ref{lem:momentsofz} and~\ref{lem:moment} as well as the definition of $r_k$ in \eqref{eq:defofrk}, we have
\begin{equation}\label{eq:fourthmomentofx1'x1}
    \begin{split}
        \E\Big(\big|\mathbf{x}_1\ttop \mathbf{x}_1 - \tr(\SIGMA)\big|^4 \Big)
    \ = \ &48 r_4 \tr(\SIGMA^4) +  12r_4 \tr(\SIGMA^2) ^2  -32 (r_3-  r_4) \tr(\SIGMA)\tr(\SIGMA^3) \\
    &+ 12( r_2  - 2r_3+  r_4)  \tr(\SIGMA)^2  \tr(\SIGMA^2)-  (3- 6r_2+4 r_3-r_4 )\tr(\SIGMA)^4\\[0.2cm]
    \ = \ &48 \tr(\SIGMA^4) +  12 \tr(\SIGMA^2)^2   + o(p^{-1})\tr(\SIGMA)^4 + o(p^{-3/4})\tr(\SIGMA)^2\tr(\SIGMA^2).
    \end{split}
\end{equation}
Applying the last few steps to \eqref{eq:varianceofT2} yields 
\begin{align*}
   \psi_2 \ \lesssim \ &\tr(\SIGMA^2) ^2 + o(p^{-1})\tr(\SIGMA)^4 + o(p^{-3/4})\tr(\SIGMA)^2\tr(\SIGMA^2).
\end{align*}
When ${\tt{r}}(\SIGMA) \lesssim \sqrt{p}$, all the three terms on the right side are $\mathcal{O}(\tr(\SIGMA^2)^2)$. On the other hand, if ${\tt{r}}(\SIGMA)/\sqrt p\to\infty$ then all three terms are $o(\tr(\SIGMA)^4)$,  which completes the proof.
\qed
\begin{lemma}
\label{lem:psi1}
Suppose that Assumption \ref{Data generating model} holds. If ${\tt{r}}(\SIGMA)/\sqrt p\to\infty$ as $n\to\infty$, then
\begin{align*}
    \psi_1 & \ = \  o\Big(\ts \frac{\tr(\SIGMA)^4}{p}\Big).
\end{align*}
Alternatively, if ${\tt{r}}(\SIGMA) \lesssim \sqrt{p}$, then 
\begin{align*}
   \psi_1 \ = \ &4 \tr(\SIGMA^4) + \ts 2\tr(\SIGMA^2) ^2 + o(\tr(\SIGMA^2) ^2).
\end{align*}
\end{lemma}
\proof
Using the definition of $h(\x_1,\x_2)$ in \eqref{eqn:definitionofh} and noting that $\x_1$, $\x_2$ and $\x_3$ are independent, we have 
\begin{equation}
\label{eq:covarianceofT2}
\begin{split}
    \psi_1 \ = \ & \big(\ts \frac{1}{4}+ \mathcal{O}(\frac{1}{n})\big) \var\big((\mathbf{x}_1\ttop \mathbf{x}_1 - \tr(\SIGMA))^2\big) \, +\, \big(4 + \ts \mathcal{O}(\frac{1}{n})\big) \var(\mathbf{x}_1\ttop\SIGMA \mathbf{x}_1 ) \\
    &-\big(2 + \ts \mathcal{O}(\frac{1}{n})\big) \cov\big((\mathbf{x}_1\ttop \mathbf{x}_1 - \tr(\SIGMA))^2,  \mathbf{x}_1\ttop \SIGMA\mathbf{x}_1 \big).
    \end{split}
\end{equation}
The first term on the right side can be analyzed by \eqref{eq:relationofr2}, \eqref{eq:fourthmomentofx1'x1} and Lemma \ref{lem:quadform}, which imply
\begin{align*}
     \var\big((\mathbf{x}_1\ttop \mathbf{x}_1 - \tr(\SIGMA))^2\big) \ = \ & 48 \tr(\SIGMA^4) + 8\tr(\SIGMA^2) ^2 + o(p^{-1})\tr(\SIGMA)^4 + o(p^{-3/4}) \tr(\SIGMA)^2\tr(\SIGMA^2).
\end{align*}
Lemma \ref{lem:quadform} and \eqref{eq:relationofr2} also give 
\begin{align*}
    \var(\mathbf{x}_1\ttop\SIGMA \mathbf{x}_1 ) 
    \ = \ & 2 \tr(\SIGMA^4) +  \mathcal{O}(p^{-1}) \tr(\SIGMA^2)^2
\end{align*}
and 
\begin{align*}
    \cov( \mathbf{x}_1\ttop \mathbf{x}_1 , \mathbf{x}_1\ttop \SIGMA\mathbf{x}_1 )
    \ = \  & 2 \tr(\SIGMA^3) +  \mathcal{O}(p^{-1}) \tr(\SIGMA) \tr(\SIGMA^2).
\end{align*}
Noting that 
\begin{align*}
    \cov\Big((\mathbf{x}_1\ttop \mathbf{x}_1 - \tr(\SIGMA))^2,  \mathbf{x}_1\ttop \SIGMA\mathbf{x}_1 \Big) \ = \ \cov\big((\mathbf{x}_1\ttop \mathbf{x}_1)^2,  \mathbf{x}_1\ttop \SIGMA\mathbf{x}_1 \big) -2\tr(\SIGMA)\cov( \mathbf{x}_1\ttop \mathbf{x}_1 , \mathbf{x}_1\ttop \SIGMA\mathbf{x}_1 ),
\end{align*}
it remains to analyze $\cov\big((\mathbf{x}_1\ttop \mathbf{x}_1)^2,  \mathbf{x}_1\ttop \SIGMA\mathbf{x}_1 \big)$. To calculate the covariance term, we may assume without loss of generality that $\SIGMA$ is diagonal, since $\x_1$ follows an elliptical distribution. It follows from the expression for $\x_i$ in \eqref{eq:xi} that 
\begin{align}\label{eq:covinpsi2}
     \cov\big((\mathbf{x}_1\ttop \mathbf{x}_1)^2,  \mathbf{x}_1\ttop \SIGMA\mathbf{x}_1 \big) \ = \ \sum_{j,k,l=1}^p \lambda_l(\SIGMA)^2\lambda_j(\SIGMA)\lambda_k(\SIGMA) \omega_{j,k,l}
\end{align}
where we define
\begin{align*}
 \omega_{j,k,l} &\ = \ \cov\Big(\ts\frac{\xi_1^4}{\|\z_1\|_2^4} z_{1j}^2 z_{1k}^2 ,  \frac{\xi_1^2}{\|\z_1\|_2^2}  z_{1l}^2\Big).
\end{align*}
This quantity can be directly calculated as
\begin{align*}
    \omega_{j,k,l} \ = \  \begin{cases}
  r_3 - r_2  &\text{ if \ \ \ } j \neq k \neq l \\
  3 (r_3 - r_2) &\text{ if \ \ \ } j = k \neq l\\
  3r_3 - r_2  &\text{ if \ \ \ } j \neq k = l 
  \text{\ \ \ or \ \  } k \neq j = l\\
  15r_3 - 3r_2  &\text{ if \ \ \ } j = k = l.
\end{cases}
\end{align*}
Applying this result to \eqref{eq:covinpsi2}, we conclude that
\begin{equation*}
% \label{eq:cov(x1^4,x1^2)}
    \begin{split}
       \cov\big((\mathbf{x}_1\ttop \mathbf{x}_1)^2,  \mathbf{x}_1\ttop \SIGMA\mathbf{x}_1 \big)  
      \ = \  &(15r_3 - 3r_2) \sum_{l=1}^p\lambda_l (\SIGMA)^4   + 2(3r_3 - r_2)\sum_{l\neq j}\lambda_l(\SIGMA)^3 \lambda_j(\SIGMA)  \\ &+3(r_3 - r_2)\sum_{l\neq j}\lambda_l(\SIGMA)^2 \lambda_j(\SIGMA)^2+ (r_3 - r_2)\sum_{j\neq k\neq l}\lambda_l(\SIGMA)^2 \lambda_j(\SIGMA)\lambda_k(\SIGMA) \\
      \ = \  &8 \tr(\SIGMA^4)   + 4 \tr(\SIGMA^3)\tr(\SIGMA)+o(p^{-3/4})\tr(\SIGMA)^2\tr(\SIGMA^2) ,
    \end{split}
\end{equation*}
where the last step is based on Lemma \ref{lem:moment}. Applying the above results to \eqref{eq:covarianceofT2} yields 
\begin{align*}
   \psi_1 \ = \ &4 \tr(\SIGMA^4) + \ts 2\tr(\SIGMA^2) ^2 + o(p^{-1})\tr(\SIGMA)^4 + o(p^{-3/4})\tr(\SIGMA)^2\tr(\SIGMA^2).
\end{align*}
When ${\tt{r}}(\SIGMA)/\sqrt p\to\infty$, the terms on the right side are all $o(\tr(\SIGMA)^4)$, and when ${\tt{r}}(\SIGMA) \lesssim \sqrt{p}$, we have $p^{-1}\tr(\SIGMA)^4 \lesssim \tr(\SIGMA^2)^2$ as well as $p^{-3/4}\tr(\SIGMA)^2 \lesssim \tr(\SIGMA^2)$,
which completes the proof of the lemma.
\qed

\section{Ratio consistency of $\hat\sigma_n^2$}\label{app:estimators}
\setcounter{lemma}{0}
\renewcommand{\thelemma}{C.\arabic{lemma}}
In order to state and prove the next result, recall from equation~\eqref{eq:definitionofsigma1} that the parameter $\sigma_n^2 = \sigma_{n,1}^2+\sigma_{n,2}^2$ is defined by
\begin{align*}
       \sigma_{n,1}^2  &\ = \ \frac{8\| R\|_4^4}{3p} \text{ \ \ \ \  \ \ \ \ \ and \ \ \ \ \ \ \ \ \ }
       \sigma_{n,2}^2  \ = \  \frac{8p(2\nu_4 +\nu_2^2)}{(\nu_1^2 +2\nu_2)^2}.
\end{align*}
In addition, recall from equation \eqref{eqn:hatsigma1def} that the estimate $\hat\sigma_n^2=\hat\sigma_{n,1}^2+\hat\sigma_{n,2}^2$ is defined by
\begin{align*}
    \hat\sigma_{n,1}^2 & \ = \ \frac{8}{3p}\Big(\|\hat R\|_4^4 - \hat\beta\|\hat R\|_4^4 + 3 \llbracket1-\hat\beta+\hat\gamma\rrbracket_{t_p}\|\hat R\|_{2}^2\Big) \text{ \ and \ } 
    \hat\sigma_{n,2}^2  \ = \ \displaystyle\frac{8p(2\hat{\nu}_4 +(\hat{\nu}_2- \frac{1}{n} \hat{\nu}_1^2)^2)}{ (\hat \nu_1^2 +2(\hat{\nu}_2- \frac{1}{n} \hat{\nu}_1^2))^2}.
\end{align*}
\begin{proposition}
\label{prop:rho}
    If Assumption \ref{Data generating model} holds, then as $n \to \infty$, 
    \begin{align*}
        \frac{\hat\sigma_n^2}{\sigma_n^2}  \ \xrightarrow{\P} \  1. 
    \end{align*}
\end{proposition}
\proof
It is sufficient to establish the limits $\hat\sigma_{n,1}^2/\sigma_{n,1}^2\xrightarrow{\P}1$ and $(\hat\sigma_{n,2}^2-\sigma_{n,2}^2)/\sigma_n^2\xrightarrow{\P} 0$, which is done in Lemmas~\ref{lem:rho1} and \ref{lem:rho2} below. \qed

\begin{lemma}
\label{lem:rho1}
    If Assumption \ref{Data generating model} holds, then as $n \to \infty$, 
    \begin{align*}
        \frac{\hat\sigma_{n,1}^2}{\sigma_{n,1}^2}  \ \xrightarrow{\P} \  1. 
    \end{align*}
\end{lemma}
\proof
Note that 
\begin{equation*}
% \label{eq:rho1}
    \begin{split}
        \frac{\hat \sigma_{n,1}^2 }{ \sigma_{n,1}^2} - 1  \ = \ & \frac{\|\hat R\|_4^4 -\| R\|_4^4 }{ \| R\|_4^4} - \frac{\hat \beta\|\hat R\|_4^4}{ \| R\|_4^4} + \frac{3\|\hat R\|_2^2 }{\| R\|_4^4 } \llbracket1-\hat\beta+\hat\gamma\rrbracket_{t_p}.
    \end{split}
\end{equation*}
Lemma \ref{lem:ZandhatZ} shows that the first term on the right is $o_{\P}(1)$. Also, Lemmas \ref{lem:ZandhatZ} and \ref{lem:consistencyofalphak} together show that the second term on the right is $o_{\P}(1)$.
For the last term, combining $\llbracket1-\hat\beta+\hat\gamma\rrbracket_{t_p} = \mathcal{O}(\log(p)p^{-3/4})$ with Lemma \ref{lem:ZandhatZ} gives
\begin{align*}
    \llbracket1-\hat\beta+\hat\gamma\rrbracket_{t_p}  \| \hat R\|_2^2 & \ = \ \mathcal{O}(\log(p)p^{-3/4})  \Big(\| R\|_2^2 + \mathcal{O}_{\P}({\|R\|_2^2})\Big)\\
    &\ = \ o_{\P}({\|R\|_4^4}),
\end{align*}
where the last step uses the inequality  $\frac{\|  R\|_2^2}{\| R\|_4^4 } \leq p^{1/2}$ that holds for a generic correlation matrix $R$.
\qed

\begin{lemma}\label{lem:rho2}
If Assumption \ref{Data generating model} holds, then as $n \to \infty$, 
    \begin{align*}
        \frac{\hat \sigma_{n,2}^2-\sigma_{n,2}^2}{\sigma_n^2}  \ \xrightarrow{\P} \  0. 
    \end{align*}    
\end{lemma}
\proof
Observe that $(\hat \sigma_{n,2}^2-\sigma_{n,2}^2)/\sigma_n^2$ can be decomposed as 
\begin{align*} 
    \frac{\hat \sigma_{n,2}^2 - \sigma_{n,2}^2}{\sigma_n^2}&\ = \ \frac{8p}{\sigma_n^2 (\nu_1^2 +2\nu_2)^2 }\bigg(\frac{(\nu_1^2 +2\nu_2)^2}{(\hat \nu_1^2 +2(\hat{\nu}_2- \ts\frac{1}{n} \hat{\nu}_1^2))^2}\Big(2\hat{\nu}_4 +(\hat{\nu}_2- \ts\frac{1}{n} \hat{\nu}_1^2)^2\Big) - (2\nu_4 + \nu_2^2)\bigg).
\end{align*}
Also, we have  
\begin{align*}
    \frac{\sigma_n^2(\nu_1^2 +2\nu_2)^2}{p}  \ = \ &8 (2\nu_4 +
    \nu_2^2) +\frac{8\| R\|_4^4}{3p^2}(\nu_1^2 +2\nu_2)^2\\
    \ \gtrsim \ & \ts \nu_2^2 + \frac{1}{p}\nu_1^4.
\end{align*}
If we can establish the following two conditions
\begin{align}
\label{eq:sigma21}
   \frac{(\nu_1^2 +2\nu_2)^2}{(\hat \nu_1^2 +2(\hat{\nu}_2- \ts\frac{1}{n} \hat{\nu}_1^2))^2}  & \ = \  1+ o_{\P}(1) \\
   \label{eq:sigma22}
    2\hat{\nu}_4 +(\hat{\nu}_2- \ts\frac{1}{n} \hat{\nu}_1^2)^2 & \ = \  2\nu_4 + \nu_2^2 + o_{\P}(\nu_2^2 + \ts\frac{1}{p}\nu_1^4),
\end{align} 
then combining with ${\nu}_4 \lesssim \nu_2^2$ yields
\begin{align*}
    \bigg|\frac{\hat \sigma_{n,2}^2 - \sigma_{n,2}^2}{\sigma_n^2} \bigg| \ \lesssim \ &  \ts \frac{1}{ \nu_2^2 + \frac{1}{p}\nu_1^4} \bigg| (1+o_{\P}(1))\Big(2\nu_4 + \nu_2^2+  o_{\P}(\nu_2^2 + \ts\frac{1}{p}\nu_1^4)\Big) - (2\nu_4 + \nu_2^2)\bigg|\\
    \ = \ & o_{\P}(1).
\end{align*}
Therefore, it is enough to show \eqref{eq:sigma21} and \eqref{eq:sigma22}. The first of these conditions follows from
\begin{align}\label{eq:consistencyofnu12ofallsample}
    \frac{\hat\nu_2-\ts \frac{1}{n} \hat\nu_1^2}{\nu_2} \  \xrightarrow{\P} \ 1 \ \ \ \ \ \text{ and } \ \ \ \ \
    \frac{\hat{\nu}_1^2}{\nu_1^2}\ \xrightarrow{\P} \ 1 ,
\end{align}
which can be proven in the same manner as~\eqref{eq:consistencyofsigma1}. Lastly, the condition~\eqref{eq:sigma22} is established in Lemma \ref{lem:consistencyof4th}.
\qed

\subsection{Norms of the sample correlation matrix}

\begin{lemma}
    \label{lem:ZandhatZ}
    If Assumption \ref{Data generating model} holds,  then 
    \begin{align}
    \label{eq:hatr4}
        \|\hat R\|_4^4 - \|R\|_4^4  & \ = o_{\P}(\|R\|_4^4) \\[0.1cm]
        \|\hat R\|_2^2 - \|  R\|_2^2 & \ = \ \mathcal{O}_{\P}({\|R\|_2^2}).\label{eq:hatr2}
    \end{align}   
\end{lemma}
\proof
The proof of Lemma \ref{lem:A_j} shows that $ \sum_{j=1}^p \P(|\hat\Sigma_{jj}/\Sigma_{jj} - 1| > n^{-1/20}) = o(1/p)$, and so 
$$\P\Big(\max_{1\leq j\leq p}\big|\hat\Sigma_{jj}/\Sigma_{jj}-1\big|>n^{-1/20}\Big)=o(1/p).$$
This implies
   $$ \Big|\|\hat R\|_4^4 - \| R \|_4^4\Big| \ \leq   \  o_{\P}(1)\sum_{j,k = 1}^p \frac{\hat \Sigma_{jk}^4}{ \Sigma_{jj}^2 \Sigma_{kk}^2} \ + \ \bigg|\sum_{j,k = 1}^p \frac{\hat \Sigma_{jk}^4}{ \Sigma_{jj}^2 \Sigma_{kk}^2}- \| R \|_4^4\bigg|$$
and
$$
    \Big|\|\hat R\|_2^2 - \| R \|_2^2\Big| \ \leq  \ o_{\P}(1)\sum_{j,k = 1}^p \frac{\hat \Sigma_{jk}^2}{ \Sigma_{jj} \Sigma_{kk}} \ + \ \bigg|\sum_{j,k = 1}^p \frac{\hat \Sigma_{jk}^2}{ \Sigma_{jj} \Sigma_{kk}}- \| R \|_2^2\bigg|.
$$
As a consequence, we may work with $\hat\SIGMA$ in place of $\hat R$ and assume that $\SIGMA=R$ when proving \eqref{eq:hatr4} and~\eqref{eq:hatr2}. 
To begin, we bound the $L^1$ error of $\|\hat \SIGMA \|_4^4$ by isolating its bias,
\begin{align}\label{eqn:L1decomp}
    \Big\|\|\hat \SIGMA\|_4^4-\|\SIGMA\|_4^4\Big\|_{L^1} \ \leq \  \Big\|\|\hat \SIGMA\|_4^4-\E\big(\|\hat \SIGMA\|_4^4\big)\Big\|_{L^1} 
\, + \, \Big|\E\big(\|\hat \SIGMA\|_4^4\big)-\|\SIGMA\|_4^4\Big|.
\end{align}
To control the bias, we calculate the expectation of $\|\hat\SIGMA\|_4^4$ directly as
\begin{align*}
    \E\big(\|\hat \SIGMA\|_4^4\big) \  = \  & \sum_{j,k = 1}^p \E(\hat\Sigma_{jk}^4)\\
     \ = \ & \frac{1}{n^4}  \sum_{j,k = 1}^p\sum_{i_1, \ldots, i_4=1}^n \E\bigg(\prod_{l=1}^4 x_{i_lj}x_{i_lk}\bigg).
\end{align*}
In the sum over $i_1,i_2,i_3,i_4$ (with $j$ and $k$ held fixed), the number of terms involving 4 distinct indices is $n(n-1)(n-2)(n-3)$, and the common expectation of these terms is $\Sigma_{jk}^4$. The number of terms involving 3 distinct indices is $6n(n-1)(n-2)$, and by the proof of Lemma \ref{lem:covariances} their common expectation is $r_2 \Sigma_{jk}^2(1+2\Sigma_{jk}^2)$.
The number of terms involving 1 or 2 distinct indices is $\mathcal{O}(n^2)$ and by \eqref{eq:momentboundofxij} their expectations are $\mathcal{O}(1)$. Combining these observations with $r_2\lesssim 1$ and the fact that $\|\SIGMA\|_2 \leq p^{1/4}\|\SIGMA\|_4^2$ holds when $\SIGMA=R$, we have 
\begin{equation}\label{eqn:4thpowerbias}
    \begin{split}
        \E\big(\|\hat \SIGMA\|_4^4\big) - \|\SIGMA\|_4^4
    \ = \ &(1+\mathcal{O}(\textstyle \frac{1}{n}))\|\SIGMA\|_4^4 + \mathcal{O}(\textstyle  \frac{1}{n}) (\|\SIGMA\|_2^2 +\|\SIGMA\|_4^4) +\mathcal{O}(\ts\frac{1}{n^2}) -\|\SIGMA\|_4^4 \\
    \ = \ & o(\|\SIGMA\|_4^4).
    \end{split}
\end{equation}
Hence, this handles the bias term in the bound~\eqref{eqn:L1decomp}. Next, we address the $L^1$ deviation term in that bound by introducing an independent copy of $\hat\SIGMA$, denoted by $\hat \SIGMA'$. Using Jensen's inequality, as well as the basic inequality $|a^4 - b^4|\leq 4|a-b|(a^3+b^3)$ for non-negative $a$ and $b$, we have
\begin{align*}
\E\bigg[\Big|\|\hat \SIGMA\|_4^4-\E\big(\|\hat \SIGMA\|_4^4\big)\Big|\bigg]
& \ \lesssim \ \E\bigg[\|\hat \SIGMA-\hat \SIGMA'\|_4\Big(\|\hat \SIGMA\|_4^3+\|\hat \SIGMA'\|_4^3\Big)\bigg]\\
& \ \lesssim \ \Big(\E\big[ \|\hat \SIGMA-\hat \SIGMA'\|_4^4\big]\Big)^{1/4}\bigg(\E\bigg[\Big(\|\hat \SIGMA\|_4^3+\|\hat \SIGMA'\|_4^3\Big)^{4/3}\bigg]\bigg)^{3/4},
\end{align*}
where H\"older's inequality with exponents $4$ and $4/3$ has been used in the last step.
By simplifying the last line with the basic inequality $(a^3+b^3)^{4/3}\leq 2^{4/3}(a^4+b^4)$ and then applying our calculation of the bias of $\|\hat\SIGMA\|_4^4$ in equation~\eqref{eqn:4thpowerbias}, we have
\begin{align*}
% \label{eqn:thirdpower}
\E\bigg[\Big|\|\hat \SIGMA\|_4^4-\E\big(\|\hat \SIGMA\|_4^4\big)\Big|\bigg] & \ \lesssim \  \Big(\E\big[ \|\hat \SIGMA-\hat \SIGMA'\|_4^4\big]\Big)^{1/4}\Big(\E\big[\|\hat \SIGMA\|_4^4\big]\Big)^{3/4}\\
& \ \lesssim \ \Big(\E\big[ \|\hat \SIGMA-\hat \SIGMA'\|_4^4\big]\Big)^{1/4} \| \SIGMA\|_4^3.
\end{align*}
Hence, it remains to show $\E\big[\|\hat \SIGMA-\SIGMA\|_4^4\big] = o(\|\SIGMA\|_4^4)$.  Letting $\mathbf{e}_{j}\in\R^p$ denote the $j$th canonical basis vector, Rosenthal's inequality (Lemma \ref{lem:Rosenthal})  gives 
\begin{align*}
    \E\big[\|\hat \SIGMA-\SIGMA\|_4^4\big]
    & \ = \  \sum_{j,k = 1}^p \bigg\| \sum_{i=1}^n  \frac{1}{n}\big(\mathbf{e}_{j}\ttop\mathbf{x}_{i} \mathbf{x}_{i}\ttop\mathbf{e}_{k}- \mathbf{e}_{j}\ttop\SIGMA\mathbf{e}_{k}\big)\bigg\|_{L^4}^4\\
    & \ \lesssim \ \sum_{j,k = 1}^p  \max\bigg\{\big(\ts \frac{1}{n}\var(\mathbf{e}_{j}\ttop\mathbf{x}_{1} \mathbf{x}_{1}\ttop\mathbf{e}_{k})\big)^2\, , \, \frac{1}{n^3}\big\|\mathbf{e}_{j}\ttop\mathbf{x}_{1} \mathbf{x}_{1}\ttop\mathbf{e}_{k}- \mathbf{e}_{j}\ttop\SIGMA\mathbf{e}_{k}\big\|_{L^4}^4\bigg\}.
\end{align*}
Since we may assume that all the diagonal entries of $\SIGMA$ are equal to 1, it follows from Lemmas \ref{lem:quadform}, \ref{lem:8thquadformofu} and \ref{lem:moment} that the variance and $L^4$ norm inside the maximum are both $\mathcal{O}(1)$.
This leads to $\E\big[\|\hat \SIGMA-\SIGMA\|_4^4\big] \lesssim 1=o(\|\SIGMA\|_4^4)$ and completes the proof of (\ref{eq:hatr4}).

To prove (\ref{eq:hatr2}), the $L^1$ triangle inequality gives
\begin{equation}\label{eqn:applyLSS}
    \Big\|\|\hat\SIGMA\|_2^2-\|\SIGMA\|_2^2\Big\|_{L^1} \ \leq \   \Big\|\big(\|\hat\SIGMA\|_2^2-\ts \frac{1}{n}\tr(\hat\SIGMA)^2\big)-\|\SIGMA\|_2^2\Big\|_{L^1} \ + \ \Big\| \ts \frac{1}{n}\tr(\hat\SIGMA)^2\Big\|_{L^1}.
\end{equation}
Lemma A.1 in \citep{wang2022bootstrap} shows that the first term on the right side of~\eqref{eqn:applyLSS} is $o(\|\SIGMA\|_2^2)$. 
For the second term in~\eqref{eqn:applyLSS}, 
Lemma \ref{lem:quadform} gives 
\begin{align*}
    \|\tr(\hat\SIGMA)^2\|_{L^1} & \ = \  \frac{n-1}{n}\tr(\SIGMA)^2 + \frac{r_2}{n} \Big(\tr(\SIGMA)^2 + 2\tr(\SIGMA^2)\Big)\\
    & \ \lesssim \ p\|\SIGMA\|_2^2,
\end{align*}
where the last step uses Lemma \ref{lem:moment} and the condition that all diagonal entries of $\SIGMA$ are equal to 1.\qed

\subsection{The limit of $\hat \beta$} 
\begin{lemma}
\label{lem:consistencyofalphak}
If Assumption \ref{Data generating model} holds, then as $n\to\infty$,
\begin{align*}
       \hat \beta \ \xrightarrow{\P} \ 0. 
\end{align*} 
\end{lemma}
\proof
Recall the definition of $\hat \beta$ in \eqref{eq:definitionofhatbeta} and let $\Delta(\hat \SIGMA) = - \ts \frac{12}{n}(\hat{\nu}_1^4 +2\hat{\nu}_1^2\hat{\nu}_2 - \ts \frac{1}{n}\hat{\nu}_1^4)$. 
Since the quantity $r_4=\frac{\E(\|\mathbf{x}_{1}\|_2^8)}{g_4(\SIGMA)}$ satisfies $r_4=1+o(1)$ by Lemma \ref{lem:moment}, we have 
\begin{align*}
    \hat \beta \ = \ &\frac{\E(\|\mathbf{x}_{1}\|_2^8)}{g_4(\SIGMA)} - \frac{\frac{1}{n}\sum_{i=1}^n \|\mathbf{x}_{i}\|_2^{8}}{g_4(\hat \SIGMA) + \Delta(\hat \SIGMA)}+ o(1) \\
    \ = \  &  r_4 \Big( 1 - \frac{g_4(  \SIGMA)}{g_4(\hat \SIGMA) + \Delta(\hat \SIGMA)}     \Big) + r_4 \frac{g_4(  \SIGMA)}{g_4(\hat \SIGMA) + \Delta(\hat \SIGMA)} \Big( 1- \frac{\ts \frac{1}{n}\sum_{i=1}^n \|\mathbf{x}_{i}\|_2^8 }{\E(\|\mathbf{x}_{1}\|_2^8)}    \Big) + o(1).
\end{align*}
Since Lemma \ref{lem:consistencyofg4} shows $g_4(  \SIGMA)/(g_4(\hat \SIGMA) + \Delta(\hat \SIGMA)) \xrightarrow{\P}1$, it is sufficient to show  $\frac{1}{n}\sum_{i=1}^n \|\mathbf{x}_{i}\|_2^8$ is a ratio-consistent estimator of $\E(\|\mathbf{x}_{1}\|_2^8)$. To this end, observe that 
%Without loss of generality, we assume that $\SIGMA$ is diagonal in the following analysis, then we have
\begin{equation*}
    \begin{split}
        \bigg\| \frac{\frac{1}{n}\sum_{i=1}^n \|\mathbf{x}_{i}\|_2^8}{\E(\|\mathbf{x}_{1}\|_2^8)} -1 \bigg\|_{L^2}^2 &  \ = \ \frac{\var(\|\x_1\|_2^8)}{n\big(\E(\|\x_1\|_2^8)\big)^2}\\
                    & \ \lesssim  \ \frac{1}{n} \frac{r_8}{r_4^2}\frac{\E((\z_1\ttop \SIGMA \z_1 )^{8})}{\big(\E((\z_1\ttop \SIGMA \z_1 )^{4})\big)^2}.
    \end{split}
\end{equation*}
Combining the fact $\E((\z_1\ttop \SIGMA \z_1 )^{4})  \geq (\E(\z_1\ttop \SIGMA \z_1))^{4}= \tr(\SIGMA)^4$ with $r_8/r_4^2\lesssim 1$ shown in Lemma \ref{lem:moment} and $\E((\z_1\ttop \SIGMA \z_1 )^{8}) \lesssim \tr(\SIGMA)^8$ given by Lemma \ref{lem:8thquadformofz}, we conclude $\frac{1}{n}\sum_{i=1}^n \|\mathbf{x}_{i}\|_2^8/\E(\|\mathbf{x}_{1}\|_2^8) \xrightarrow{\P} 1$, which completes the proof.
\qed

\begin{lemma}\label{lem:consistencyofg4}
If Assumption \ref{Data generating model} holds, then as $n\to\infty$,
    \begin{align*}
        \frac{g_4(\hat \SIGMA) + \Delta(\hat \SIGMA)}{g_4(  \SIGMA)} \ \xrightarrow{\P} \ 1.
    \end{align*}
\end{lemma}
\proof
Noting that equation \eqref{eq:consistencyofnu12ofallsample} and Lemma \ref{lem:momentsofz} imply $\Delta(\hat \SIGMA)/ g_4(  \SIGMA) \xrightarrow{\P} 0$, it is enough to show $g_4(\hat \SIGMA)$ is a ratio-consistent estimator of $g_4(\SIGMA)$.
Let $\z\in\R^p$ denote a standard normal vector independent of $\x_1,\dots,\x_n$, and observe that 
\begin{align*}
     g_4(\hat \SIGMA)-g_4(  \SIGMA) \ = \ \E\Big[(\z\ttop\hat\SIGMA \z)^4-(\z\ttop \SIGMA \z)^4\Big|\hat\SIGMA\Big].
\end{align*}
So, it is sufficient to show that
\begin{align*}
% \label{eqn:tempreduction}
    \Big\|(\z\ttop\hat\SIGMA \z)^4-(\z\ttop \SIGMA \z)^4\Big\|_{L^1} \ = \ o(g_4(\Sigma)).
\end{align*}
Using the basic inequality $|a^4 - b^4|\leq 4|a-b|\sqrt{a^6+b^6}$ that holds for non-negative $a$ and $b$, it follows from the Cauchy-Schwarz inequality that
\begin{align}\label{eqn:CS6}
\Big\|(\z\ttop\hat\SIGMA \z)^4-(\z\ttop \SIGMA \z)^4\Big\|_{L^1} 
\ \lesssim \ \Big\| \z\ttop\hat\SIGMA \z-\z\ttop \SIGMA \z\Big\|_{L^2}\Big\|(\z\ttop\hat\SIGMA \z)^6+(\z\ttop \SIGMA \z)^6\Big\|_{L^1}^{1/2}.
\end{align}
To handle the second factor,  Lemma \ref{lem:8thquadformofz} implies
\begin{equation*}
    \Big\|(\z\ttop\hat\SIGMA \z)^6+(\z\ttop \SIGMA \z)^6\Big\|_{L^1}^{1/2} \ \lesssim \  \Big( \E(\tr(\hat\SIGMA)^6)+\tr(\SIGMA)^6\Big)^{1/2}.
\end{equation*}
By viewing $\tr(\hat\SIGMA)$ as the sample average $\frac{1}{n}\sum_{i=1}^n \|\x_i\|_2^2$, Jensen's inequality implies
\begin{align*}
    \E\big(\tr(\hat \SIGMA)^6\big) & \ \leq \ \frac{1}{n}\sum_{i=1}^n \E(\|\x_1\|_2^{12})\\
    & \ \lesssim \ \tr(\SIGMA)^6,
\end{align*} 
where the last step follows from Lemma \ref{lem:8thquadformofu}. To handle the first factor on the right side of \eqref{eqn:CS6}, observe that Lemma \ref{lem:momentsofz} gives
$$\E\Big[ \big(\z\ttop\hat\SIGMA \z-\z\ttop \SIGMA \z\big)^2 \Big|\hat\SIGMA\Big] \ = \ 2\|\hat\SIGMA-\SIGMA\|_2^2+\big(\tr(\hat\SIGMA)-\tr(\SIGMA)\big)^2.$$
Furthermore, Lemmas \ref{lem:quadform} and \ref{lem:moment} imply
\begin{align*}
    \E\Big(\|\hat \SIGMA-\SIGMA\|_2^2\Big)
    & \ = \   \frac{1}{n}\sum_{j,k = 1}^p\var(\mathbf{x}_{1}\ttop\mathbf{e}_{k}\mathbf{e}_{j}\ttop\mathbf{x}_{1} ) \\
    & \ \lesssim \ \frac{1}{n} \sum_{j,k = 1}^p  ( \Sigma_{jk}^2 + \Sigma_{jj}\Sigma_{kk}) \\[0.2cm]
    & \ \lesssim \  \frac{1}{n}\tr(\SIGMA)^2
\end{align*} 
as well as
\begin{align*}
    \E\big((\tr(\hat\SIGMA)-\tr(\SIGMA))^2\big)
    \ = \ & \frac 1n \var(\x_1\ttop \x_1)
    \ \lesssim \ \frac{1}{n}\tr(\SIGMA)^2.
\end{align*}
The last few steps imply 
\begin{align*}
    \Big\| \z\ttop\hat\SIGMA \z-\z\ttop \SIGMA \z\Big\|_{L^2} \ = \ \frac{1}{\sqrt{n}}\tr(\SIGMA).
\end{align*}
Applying the above results to \eqref{eqn:CS6} and observing $g_4(\Sigma)\geq\tr(\SIGMA)^4$, we have 
\begin{align*}
    \Big\|(\z \ttop \hat\SIGMA \z ) ^4 - (\z\ttop \SIGMA \z)^4\Big\|_{L^1} \ = \ \mathcal{O}(\ts \frac{1}{\sqrt{n}}\tr(\SIGMA)^4) \ = \ o(g_4(\Sigma)),
\end{align*}
which completes the proof.
\qed

\subsection{Consistency of $\hat \nu_4$}
\begin{lemma}
\label{lem:consistencyof4th}
If Assumption \ref{Data generating model} holds, then as $n \to \infty$, 
    \begin{align*}
        \frac{ \hat \nu_4 - \nu_4}{ \nu_2^2 + \frac{1}{p}\nu_1 ^4} \ \xrightarrow{\P} \ 0. 
    \end{align*}
\end{lemma}
\proof
It is enough to show $\E\big(|\tr(\hat\SIGMA^4)-\tr(\SIGMA^4)|\big) = o\big( \tr(\SIGMA^2)^2 + \frac{\tr(\SIGMA)^4}{p}\big)$.
First, observe that the Cauchy-Schwarz inequality implies
\begin{equation*}
    \begin{split}
|\tr(\hat\SIGMA^4)-\tr(\SIGMA^4)|^2 & \ = \ \Big|\sum_{j=1}^p \lambda_j(\hat\SIGMA)^4 -\lambda_j(\SIGMA)^4\Big|^2\\[0.2cm]
& \ \leq \  \sum_{j=1}^p \Big(\lambda_j(\hat\SIGMA)^2 -\lambda_j(\SIGMA)^2\Big)^2 \sum_{j=1}^p\Big(\lambda_j(\hat\SIGMA)^2 +\lambda_j(\SIGMA)^2\Big)^2\\[0.2cm]
& \ \leq \  2\|\hat\SIGMA^2-\SIGMA^2\|_2^2 \Big(\tr(\hat\SIGMA^4)+\tr(\SIGMA^4)\Big),
    \end{split}
\end{equation*}
where Wielandt's inequality~\cite[][Corollary 7.3.5]{horn2012matrix} has been used in the last step.
Taking the square root of both sides, followed by the expectation, the Cauchy-Schwarz inequality implies
\begin{equation}\label{eq:l1oftr}
\small
    \begin{split}
        \E\big(|\tr(\hat\SIGMA^4)-\tr(\SIGMA^4)|\big) &\ \lesssim  \ \sqrt{\E\big(\|\hat\SIGMA^2-\SIGMA^2\|_2^2\big)}\sqrt{\E(\tr(\hat\SIGMA^4))+\tr(\SIGMA^4)}\\%
        & \ = \ \sqrt{\E(\tr(\hat\SIGMA^4))-2\E(\tr(\SIGMA^2\hat\SIGMA^2))+\tr(\SIGMA^4)}\sqrt{\E(\tr(\hat\SIGMA^4))+\tr(\SIGMA^4)}.
    \end{split}
\end{equation}
In the remainder of the proof, we will separately analyze $\E(\tr(\hat\SIGMA^4))$ and $\E(\tr(\SIGMA^2\hat\SIGMA^2))$. For the first of these quantities, note that 
\begin{align}\label{eqn:Etr4}
   \E(\tr(\hat \SIGMA^4)) & \ = \ \frac{1}{n^4} \ \sum_{i_1,\ldots,i_4 =1}^n\E\bigg(\tr\Big(\textstyle\prod_{l=1}^{4} \x_{i_l}\x_{i_l}\ttop\Big)\bigg).
\end{align} 
The number of terms involving 4 distinct indices is $n(n-1)(n-2)(n-3)$ and the corresponding expectation is $\tr(\SIGMA^4)$. For the terms in~\eqref{eqn:Etr4} that involve 3 distinct indices, there are two possible cases. One case occurs when $i_1=i_2\neq i_3 \neq i_4$. The number of corresponding terms is $\mathcal{O}(n^3)$, and by Lemmas \ref{lem:quadform} and \ref{lem:moment}, the expectation of each such term is $r_2(2\tr(\SIGMA^4) + \tr(\SIGMA^3)\tr(\SIGMA)) \lesssim \sqrt{p}\tr(\SIGMA^2)^2$. The second case occurs when $i_1=i_3\neq i_2 \neq i_4$. The number of corresponding terms is $\mathcal{O}(n^3)$, and by Lemmas \ref{lem:quadform} and \ref{lem:moment}, the expectation of each such term is $r_2(2\tr(\SIGMA^4) + \tr(\SIGMA^2)^2) \lesssim \tr(\SIGMA^2)^2$. The number of terms in~\eqref{eqn:Etr4} that involve 1 or 2 distinct indices is $\mathcal{O}(n^2)$, with their expectations being at most $\mathcal{O}(\tr(\SIGMA)^4)$. Consequently, we have 
\begin{align*}
    \E(\tr(\hat \SIGMA^4)) 
    & \ = \ \tr(\SIGMA^4) + \mathcal{O}(\ts \frac{\sqrt{p}}{n})\tr(\SIGMA^2)^2 + \mathcal{O}(\ts \frac{1}{n^2})\tr(\SIGMA)^4.
\end{align*} 
In addition, Lemmas \ref{lem:quadform} and \ref{lem:moment} give
\begin{align*}
    \E(\tr(\SIGMA^2\hat\SIGMA^2)) & \ = \ \frac{1}{n^2}\sum_{i_1,i_2 = 1}^n\E( \x_{i_1}\ttop \x_{i_2} \x_{i_2}\ttop\SIGMA^2 \x_{i_1}) \\& \ = \ \frac{n-1}{n} \tr(\SIGMA^4) + \frac{r_2}{n}\Big( \tr(\SIGMA^3)\tr(\SIGMA)+2\tr(\SIGMA^4)\Big)\\
    & \ = \ \tr(\SIGMA^4) + \mathcal{O}(\ts \frac{\sqrt{p}}{n})\tr(\SIGMA^2)^2. 
\end{align*}
Applying the last two bounds to \eqref{eq:l1oftr}, it follows that 
\begin{align*}
    \E\big(|\tr(\hat\SIGMA^4)-\tr(\SIGMA^4)|\big) &\ \lesssim \ \Big( \ts\frac{\sqrt{p}}{n}\tr(\SIGMA^2)^2 + \frac{1}{n^2}\tr(\SIGMA)^4 \Big)^{1/2} \Big (\tr(\SIGMA^4) + \ts \frac{\sqrt{p}}{n}\tr(\SIGMA^2)^2 + \frac{1}{n^2}\tr(\SIGMA)^4  \Big)^{1/2}\\
    & \ = \ o\Big(\big( \tr(\SIGMA^2)^2 + \ts \frac{\tr(\SIGMA)^4}{p}\big)^{1/2}\Big) \cdot \mathcal{O}\Big(\big( \tr(\SIGMA^2)^2 + \ts \frac{\tr(\SIGMA)^4}{p}\big)^{1/2}\Big)\\
    & \ = \ o\big( \tr(\SIGMA^2)^2 + \ts \frac{\tr(\SIGMA)^4}{p}\big),
\end{align*}
which implies the stated result.\qed

\section{Background results} \label{app:background}
\setcounter{lemma}{0}
\renewcommand{\thelemma}{D.\arabic{lemma}}
\normalsize

% {\blue{please use index notation for a data vector like in the rest of the paper e.g. $\x_1$ rather than $\x$}}}
\begin{lemma}\label{lem:equal_kurt}
     Let $\x_1 =  \xi_1 \Sigma^{1/2} \u_1$, where $\Sigma \in \R^{p\times p}$ is deterministic and positive semidefinite with positive diagonal entries, $\u_1$ is a random vector drawn from the uniform distribution on the unit sphere of $\R^p$, and $\xi_1\geq 0 $ is a random variable independent of $\u_1$. In addition, suppose that $\E(x_{1j}^4)<\infty$ for all $j=1,\dots,p$. Then, $\frac{\E(x_{11}^4)}{(\E(x_{11}^2))^2} = \cdots =\frac{\E(x_{1p}^4)}{(\E(x_{1p}^2))^2}$.
\end{lemma}
\proof
Note that the rotation invariance of the uniform distribution on the surface of the unit sphere implies that $x_{1j} =\mathbf{e}_j \ttop \Sigma^{1/2}\xi_1 \u_1$ is equal in distribution to  $\| \Sigma^{1/2}\mathbf{e}_j\|_{2}\xi_1 u_{11}$, where $\mathbf{e}_j$ denotes the $j$th standard basis vector (with 1 in the $j$th entry and 0 in all other entries). Combining with the fact that kurtosis is a scale-invariant parameter, it follows that the kurtosis of $x_{1j}$  is equal to that of $\xi_1 u_{11}$ for all $j = 1,\ldots, p$.

% Note that $\u$ can be expressed as $\frac{\z}{\|\z\|_2}$, where $\z$ is a standard normal vector in $\R^p$ independent of $\xi$. It follows that $x_i = \frac{\xi}{\|\z\|_2}\sum_{j=1}^p M_{ij}z_j$, then the properties of multivariate normal distribution imply
% \begin{align*}
%     \E(x_i^2) \ = \ &\frac{\E(\xi^2) }{\E(\|\z\|_2^2)}\E\Big(\big(\sum_{j=1}^p M_{ij}z_j\big)^2\Big)\\
%     \ = \ &\frac{\E(\xi^2) }{p}\sum_{j=1}^p M_{ij}^2
% \end{align*}
% and
% \begin{align*}
%     \E(x_i^4) \ = \ &\frac{\E(\xi^4) }{\E(\|\z\|_2^4)}\E\Big(\big(\sum_{j=1}^p M_{ij}z_j\big)^4\Big)\\
%     \ = \ &\frac{\E(\xi^4) }{p(p+2)}\Big(\sum_{j=1}^p 3 M_{ij}^4 + 3\sum_{1\leq j\neq k\leq p} M_{ij}^2 M_{ik}^2\Big) \\
%     \ = \ &\frac{3 \E(\xi^4) }{p(p+2)}\Big(\sum_{j=1}^p M_{ij}^2\Big)^2.
% \end{align*}
% Therefore, we have for any $i=1,2,\ldots,p$, 
% \begin{align*}
%     \frac{\E(x_i^4)}{\E(x_i^2)^2} = \frac{3p \E(\xi^4) }{(p+2)\E(\xi^2)^2}.
% \end{align*}
\qed

\begin{lemma}[\cite{hu2019aos}, Lemma A.1]\label{lem:quadform}
     Let $\xi_1\in\R$ and $\mathbf{u}_1\in\R^p$ satisfy the conditions in Assumption \ref{Data generating model}, and fix any matrices $M, \tilde M \in \R^{p\times p}$. Then,
     \small
    \begin{align*}
       \cov \Big(\xi_1^2\mathbf{u}_1\ttop M \mathbf{u}_1 \,,\, \xi_1^2\mathbf{u}_1\ttop \tilde M \mathbf{u}_1   \Big) \ = \  \frac{\E(\xi_1^4)}{p(p+2)}\Big(\tr(M)\tr(\tilde M)+\tr(M\tilde M) +\tr(M\tilde M\ttop)\Big) - \tr(M)\tr(\tilde M) .
    \end{align*}
    \normalsize
\end{lemma}

\begin{lemma}[\cite{bai2010spectral}, Lemma B.26]\label{lem:8thquadformofz}
     Let $\z_1 \in\R^p$ be a standard normal random vector, and fix any matrix $M\in \R^{p\times p}$. Then, for any $q\geq 1$,
     $$\big\|\z_1\ttop M \z_1-\tr(M)\big\|_{L^q} \ \leq \ C_q \tr(M M\ttop)^{1/2} ,$$
    where $C_q$ is a constant depending on $q$ only.
\end{lemma}
\begin{lemma}\label{lem:8thquadformofu}
If Assumption \ref{Data generating model} holds, then for any fixed matrix $M \in \R^{p\times p}$ and any \smash{$1\leq q\leq 8$,} 
\begin{align*}
    \Big\|\mathbf{x}_1\ttop M \mathbf{x}_1 - \tr(M\SIGMA)\Big\|_{L^q} &  \ \lesssim \ \sqrt{\tr(M\SIGMA M\ttop\SIGMA)} \ + \ o(\ts \frac{1}{p^{1/4}})\tr(M \SIGMA).
\end{align*} 
\end{lemma}
\proof 
By Lyapunov's inequality, it is enough to handle the case when $q=8$. For any fixed matrix $M \in \R^{p\times p}$, the definition of $\x_1$ and Assumption \ref{Data generating model} give
\small
\begin{align*}
    \Big\|\mathbf{x}_1\ttop M \mathbf{x}_1 - \tr(M\SIGMA)\Big\|_{L^8}^8 & \ \lesssim \  \E(\xi_1^{16}) \Big\|\mathbf{u}_1\ttop \SIGMA ^{1/2}  M \SIGMA ^{1/2} \mathbf{u}_1 - \ts \frac{1}{p}\tr(M\SIGMA)\Big\|_{L^8}^8  \ + \ \tr(M\SIGMA)^8\Big\|\ts \frac{1}{p}\xi_1^2-1\Big\|_{L^8}^8 \\[0.2cm]
    &  \ \lesssim \ p^8\Big\|\mathbf{u}_1\ttop \SIGMA ^{1/2}  M \SIGMA ^{1/2} \mathbf{u}_1 - \ts \frac{1}{p}\tr(\SIGMA^{1/2}M\SIGMA^{1/2})\Big\|_{L^8}^8  \ + \ o(\textstyle \frac{1}{p^2})\tr(M\SIGMA)^8.
\end{align*}
\normalsize
If we can show
\begin{align} \label{eq:8thmomentofu}
    p^{8}  \Big\|\mathbf{u}_1\ttop A  \mathbf{u}_1 - \ts\frac{1}{p}\tr(A)\Big\|_{L^8}^8  \ \lesssim \  \tr(A A\ttop)^4
\end{align}
holds for any matrix $A \in \R^{p\times p}$, then the sub-additivity of the eighth-root function implies the stated result.
To show \eqref{eq:8thmomentofu}, observe that
\begin{align*}
    p^{8}  \Big\|\mathbf{u}_1\ttop A  \mathbf{u}_1 - \ts\frac{1}{p}\tr(A)\Big\|_{L^8}^8 \  \lesssim \  \Big\| p \u_1\ttop A \u_1 -\z_1\ttop A \z_1\ttop\Big\|_{L^8}^8 \ + \ \Big\|\z_1\ttop A \z_1-\tr(A)\Big\|_{L^8}^8,
    % \label{eqn:8thmoment}
\end{align*}
where $\z_1$ denotes a standard normal vector in $\R^p$ satisfying $\u_1 = \frac{\z_1}{\|\z_1\|_2}$.
The second term of this bound is $\mathcal{O}(\tr (AA\ttop)^4)$ due to Lemma \ref{lem:8thquadformofz}, while the first term can be rewritten as
$$\Big\| p \u_1\ttop A \u_1 -\z_1\ttop A \z_1\ttop\Big\|_{L^8}^8 \ = \  \big\|p-\|\z_1\|_2^2\big\|_{L^8}^8\big\|\u_1\ttop A\u_1\|_{L^8}^8.$$
Lemma \ref{lem:8thquadformofz} shows that the first factor on the right is $\mathcal{O}(p^4)$,
and combining with the lower bound $\E(\|\mathbf{z}_{1}\|_2^{16}) \geq (\E(\|\mathbf{z}_{1}\|_2^2))^{8}\geq p^{8}$ implies that the second factor satisfies
\begin{equation*}
\begin{split}
\big\|\u_1\ttop A\u_1\|_{L^8}^8  & \ \lesssim \ \frac{\big\|\z_1\ttop A\z_1-\tr(A)\big\|_{L^8}^8}{\E(\|\z_1\|_2^{16})}+ \frac{\tr(A)^8}{\E(\|\z_1\|_2^{16})}\\[0.2cm]
& \ \lesssim \ \frac{\tr(A A\ttop)^4}{p^8}+\frac{\tr(A)^8}{p^8}.
\end{split}
\end{equation*}
This leads to \eqref{eq:8thmomentofu}, since $\tr(A)^8\leq p^4\tr(A A\ttop)^4$.
\qed

\begin{lemma}[\cite{magnus1978moments}, Lemma 2.2]\label{lem:momentsofz}
     If $\z_1\in\R^p$ is a standard normal vector and $M\in \R^{p\times p}$ is a symmetric matrix, then
     \begin{align*}
     \E((\z_1\ttop M \z_1 )^2)  \ = \ & 2 \tr(M^2) + \tr(M)^2,\\
    \E((\z_1\ttop M \z_1 )^3) \ = \ & 8\tr(M^3) +6 \tr(M^2)\tr(M)+ \tr(M)^3, \\
    \E((\z_1\ttop M \z_1 )^4) \ = \ &48 \tr(M^4)  +32 \tr(M^3)\tr(M) + 12\tr(M^2)^2  + 12\tr(M^2)\tr(M)^2 + \tr(M)^4.
\end{align*}
\end{lemma}

\begin{lemma}
\label{lem:covariances}
    If Assumption \ref{Data generating model} holds, then for any $j,k \in \{1,\ldots,p\}$
    \begin{align*}
        \cov(x_{1j}^4,x_{1k}^4)  & \ = \ r_4 (72\Sigma_{jj}\Sigma_{kk}\Sigma_{jk}^2+24\Sigma_{jk}^4)+ 9 (r_4 - r_2^2)\Sigma_{jj}^2 \Sigma_{kk}^2,\\[0.2cm]
        \cov(x_{1j}^2,x_{1k}^2) 
        & \ = \ 2r_2\Sigma_{jk}^2 + ( r_2-1)\Sigma_{jj} \Sigma_{kk},\\[0.2cm]
        \cov(x_{1k}^2,x_{1j}^4) & \ = \ 12 r_3\Sigma_{jj}\Sigma_{jk}^2  + 3 ( r_3- r_2)\Sigma_{jj}^2 \Sigma_{kk} .
    \end{align*}
\end{lemma}

\proof
Recall the relation $\x_1=(\xi_1/\|\z_1\|_2)\boldsymbol\zeta_1$ from~\eqref{eq:xi}, where $\boldsymbol\zeta_1 \sim N(0, \SIGMA)$. Isserlis' theorem (as presented in~\eqref{eqn:isserlis}) ensures that for any $j,k \in \{1,\ldots, p\}$ we have
\begin{align*}
    \cov(x_{1j}^4,x_{1k}^4) & \ = \ \frac{\E(\xi_1^8)}{\E(\|\z_1\|_2^8)}\E(\zeta_{1j}^4\zeta_{1k}^4) - \Big(\frac{\E(\xi_1^4)}{\E(\|\z_1\|_2^4)}\Big)^2\E(\zeta_{1j}^4)\E(\zeta_{1k}^4)\\[0.2cm]
    & \ = \   r_4(9\Sigma_{jj}^2\Sigma_{kk}^2+72\Sigma_{jj}\Sigma_{kk}\Sigma_{jk}^2+24\Sigma_{jk}^4) - 9  r_2^2 \Sigma_{jj}^2 \Sigma_{kk}^2, \\[0.3cm]
    \cov(x_{1j}^2,x_{1k}^2) & \ = \ \frac{\E(\xi_1^4)}{\E(\|\z_1\|_2^4)} \E(  \zeta_{1j}^2\zeta_{1k}^2)- \E(\zeta_{1j}^2)\E(\zeta_{1k}^2)\\[0.2cm]
& \ = \  r_2 (\Sigma_{jj}\Sigma_{kk}+2\Sigma_{jk}^2) - \Sigma_{jj} \Sigma_{kk},\\[0.3cm]
    \cov(x_{1k}^2,x_{1j}^4)& \ = \  \frac{\E(\xi_1^6 )}{\E(\|\z_1\|_2^6)}\E( \zeta_{1j}^4\zeta_{1k}^2) - \frac{\E(\xi_1^4)}{\E(\|\z_1\|_2^4)}\E(\zeta_{1j}^4)\E(\zeta_{1k}^2)\\[0.2cm]
    & \ = \  r_3(3\Sigma_{jj}^2\Sigma_{kk}+12\Sigma_{jj}\Sigma_{jk}^2) - 3  r_2 \Sigma_{jj}^2 \Sigma_{kk}.
\end{align*}
\qed

\begin{lemma}[\cite{johnson1985best}, Rosenthal's inequality]
\label{lem:Rosenthal}
     Fix $r > 2$ and let $v_1, \ldots, v_n$ be centered independent random variables. Then, there is an absolute constant $c>0$ such that
    \begin{align*}
        \bigg\|\sum_{i=1}^n v_i \bigg\|_{L^r} \ \leq \ c \cdot r\cdot \max\bigg\{\bigg \|\sum_{i=1}^n v_i \bigg \|_{L^2}, \Big(\sum_{i=1}^n\|v_i\|_{L^r}^r\Big)^{1/r}\bigg\}.
    \end{align*}
\end{lemma}
\qed

\section{Additional numerical experiments}\label{app:dim}

The following is an extension of Table~\ref{table:level} from Section~\ref{sec:level}, which reports the empirical level of the proposed test for larger values of the ratio $p/n\in\{3,5,10\}$.

\begin{table}[H]
% \small
\centering
\setlength\tabcolsep{6.5pt}
\caption{Results on empirical level for the proposed test (5\% nominal level).}
\begin{tabular}{lllllllllllllll}
\hline
      &\multicolumn{4}{c}{$p/n=3$} &  & \multicolumn{4}{c}{$p/n=5$} &  & \multicolumn{4}{c}{$p/n=10$} \\ \cline{2-5} \cline{7-10} \cline{12-15} 
      \specialrule{0em}{1pt}{1pt}
 & (1)   & (2)   & (3)   & (4) &  & (1)   & (2)   & (3)   & (4) &  & (1)   & (2)   & (3)   & (4) \\
\hline
\specialrule{0em}{1pt}{1pt}
(i) &4.26&4.34 &4.00 &4.29 &  &3.89&3.84 &3.71 &3.83 &  &3.89&3.86 &3.80 &3.81        \\
 \specialrule{0em}{1pt}{1pt}
 \hline
\specialrule{0em}{1pt}{1pt}
(ii) &4.02&4.13 &3.88 &4.17 &  &3.74&3.85 &3.87 &3.89 &  &3.77&3.73 &3.71 &3.81        \\
 \specialrule{0em}{1pt}{1pt}
 \hline
\specialrule{0em}{1pt}{1pt}
(iii) &4.43&4.51 &4.22 &4.47 &  &4.12&4.11 &3.89 &4.06 &  &4.10&4.00 &3.92 &4.05 \\
 \specialrule{0em}{1pt}{1pt}
 \hline
\specialrule{0em}{1pt}{1pt}
(iv) & 3.99&4.18 &3.95 &4.17 &  &3.73&3.80 &3.97 &3.83 &  &3.62&3.61 &3.33 &3.71        \\
 \specialrule{0em}{1pt}{1pt}
 \hline
\specialrule{0em}{1pt}{1pt}
(v) &4.07&4.23 &3.99 &4.12 &  &3.78&3.79 &3.92 &3.87 &  &3.73&3.68 &3.66 &3.68         \\
 \specialrule{0em}{1pt}{1pt}
 \hline
\end{tabular}
\label{table:levelapp}
\end{table}

% \begin{table}[H]
% \centering
% \setlength\tabcolsep{10pt}
% \caption{Results on empirical level for the proposed test (5\% nominal level)}
% \begin{tabular}{lllllllllllll}
% %
% \hline
%       &  & \multicolumn{3}{c}{$p/n=3$} &  & \multicolumn{3}{c}{$p/n=5$} &  & \multicolumn{3}{c}{$p/n=10$} \\ \cline{3-5} \cline{7-9} \cline{11-13} 
%       \specialrule{0em}{1pt}{1pt}
%  &  & (1)   & (2)   & (3) &  & (1)   & (2)   & (3) &  & (1)   & (2)   & (3) \\
% \hline
% \specialrule{0em}{1pt}{1pt}
%  %
% (i) & &4.26&4.34 &4 &  &3.89&3.84 &3.71 &  &3.89&3.86 &3.80         \\
% %
%  \specialrule{0em}{1pt}{1pt}
%  \hline
% \specialrule{0em}{1pt}{1pt}
%  %
% (ii) & &4.02&4.13 &3.88 &  &3.74&3.85 &3.87 &  &3.77&3.73 &3.71        \\
%  %
%  \specialrule{0em}{1pt}{1pt}
%  \hline
% \specialrule{0em}{1pt}{1pt}
%  % 
% (iii) & &4.43&4.51 &4.22 &  &4.12&4.11 &3.89 &  &4.10&4.00 &3.92 \\
%  \specialrule{0em}{1pt}{1pt}
%  \hline
% \specialrule{0em}{1pt}{1pt}
%  %
% (iv) & &3.99&4.18 &3.95 &  &3.73&3.8 &3.97 &  &3.62&3.61 &3.33         \\
%  %
%  \specialrule{0em}{1pt}{1pt}
%  \hline
% \specialrule{0em}{1pt}{1pt}
%  %   
% (v) & &4.07&4.23 &3.99 &  &3.78&3.79 &3.92 &  &3.73&3.68 &3.66          \\
%  %
%  \specialrule{0em}{1pt}{1pt}
%  \hline
% \end{tabular}
% \label{table:level}
% \end{table}
\end{document}